\documentclass{ArticleModern}

\usepackage{SymbolsModern}
\usepackage{StylesModern}

\newcommand{\ForestF}{\mathfrak{f}}
\newcommand{\TrimmedForestF}{f}
\newcommand{\TrimmedTermT}{t}
\newcommand{\MultisetM}{\mathfrak{m}}

\newcommand{\Reduced}{\mathrm{rd}}
\newcommand{\Weight}{\mathrm{w}}
\newcommand{\Charge}{\mathrm{ch}}
\newcommand{\Content}{\mathrm{ct}}
\newcommand{\Infix}{\mathrm{in}}
\newcommand{\Trim}{\mathrm{tr}}
\newcommand{\Pack}{\mathrm{pck}}
\newcommand{\Edge}[2]{\mathrel{\overset{#1}{\rightarrow}}_{#2}}
\newcommand{\SignatureExample}{\Signature_{\mathrm{e}}}

\newcommand{\SetForests}{\mathfrak{F}}

\newcommand{\SetDecoratedLetters}{\mathrm{L}}

\newcommand{\Occurrences}[1]{\mathrm{Pos}_{#1}}
\newcommand{\WordToTensor}{\theta}
\newcommand{\ClassAlphabets}{\mathsf{A}}
\newcommand{\ClassForestLikeAlphabets}{\mathsf{F}}
\newcommand{\ClassTotallyOrderedAlphabets}{\mathsf{O}}
\newcommand{\Alphabet}{\mathbb{A}}
\newcommand{\AlphabetSum}{\mathbin{+\mkern-16mu+}}
\newcommand{\AlphabetProduct}{\mathbin{\times\mkern-16mu\times}}
\newcommand{\Monomial}{\mathrm{m}}
\newcommand{\PositionsEncoding}{\mathrm{pos}}

\newcommand{\RelatedAlphabetSignature}{\mathcal{R}}
\newcommand{\PositionsAlphabet}{\Alphabet_{\mathrm{p}}}
\newcommand{\LengthsAlphabet}{\Alphabet_{\mathrm{l}}}
\newcommand{\ProjectionLength}{\pi_{\mathrm{l}}}

\newcommand{\Letter}{\mathbf{a}}
\newcommand{\NAryRelation}{\mathfrak{R}}
\newcommand{\RootRelation}{\mathrm{R}}
\newcommand{\DecorationRelation}{\mathrm{D}}
\newcommand{\EdgeRelation}[1]{\mathrel{\prec_{#1}}}
\newcommand{\CompatibleWord}[1]{\mathrel{\Vdash^{#1}}}
\newcommand{\AdmissibleSet}{\mathrel{\vdash}}
\newcommand{\Position}{\mathrm{p}}
\newcommand{\Decoration}{\mathrm{d}}
\newcommand{\Value}{\mathrm{v}}

\newcommand{\NaturalHopfAlgebra}{\mathbf{N}}
\newcommand{\Space}{\mathcal{V}}
\newcommand{\HilbertSeries}{\mathcal{F}}

\newcommand{\BasisE}{\mathsf{E}}
\newcommand{\BasisM}{\mathsf{M}}
\newcommand{\Realization}{\mathsf{r}}
\newcommand{\NCSF}{\mathbf{Sym}}
\newcommand{\FQSym}{\mathbf{FQSym}}
\newcommand{\WQSym}{\mathbf{WQSym}}
\newcommand{\PQSym}{\mathbf{PQSym}}
\newcommand{\NCFdB}{\mathbf{FdB}}
\newcommand{\NCK}{\mathbf{NCK}}
\newcommand{\Phr}{\mathbf{Phr}}

\newcommand{\As}{\mathsf{As}}
\newcommand{\MAs}{\mathsf{MAs}}
\newcommand{\Int}{\mathsf{Int}}
\newcommand{\FCat}{\mathsf{FCat}}
\newcommand{\Schr}{\mathsf{Schr}}
\newcommand{\Motz}{\mathsf{Motz}}
\newcommand{\Comp}{\mathsf{Comp}}
\newcommand{\DA}{\mathsf{DA}}
\newcommand{\SComp}{\mathsf{SComp}}
\newcommand{\Dias}{\mathsf{Dias}}
\newcommand{\Dendr}{\mathsf{Dendr}}

\Title{%
    Polynomial realizations of Hopf algebras \\ built from nonsymmetric operads%
    \vspace{-1.5ex}
}
\TitleShort{Realizations of natural Hopf algebras of operads}
\Author{Samuele Giraudo \vspace{-1.5ex}}
\AuthorShort{S. Giraudo}
\Subjects{%
    16T30, 
    18M80, 
    05E05. 
}
\Date{\today}
\Address{%
    Université du Québec à Montréal, LACIM, \\
    Pavillon Président-Kennedy, 201 Avenue du Président-Kennedy, Montréal, H2X~3Y7, Canada.
    \vspace{-1.5ex}
}
\Email{giraudo.samuele@uqam.ca}
\Funding{%
    This research has been partially supported by the projects CARPLO (ANR-20-CE40-0007) and
    ALCOHOL (ANR-19-CE40-0006) of the Agence nationale de la recherche, and by the Natural
    Sciences and Engineering Research Council of Canada (RGPIN-2024-04465).
}
\Abstract{%
    The natural Hopf algebra $\NaturalHopfAlgebra \App \Operad$ of an operad $\Operad$ is a
    Hopf algebra whose bases are indexed by some words on $\Operad$. We construct polynomial
    realizations of $\NaturalHopfAlgebra \App \Operad$ by using alphabets of noncommutative
    variables endowed with unary and binary relations. By using particular alphabets, we
    establish links between $\NaturalHopfAlgebra \App \Operad$ and some other Hopf algebras
    including the Hopf algebra of word quasi-symmetric functions of Hivert, the decorated
    versions of the noncommutative Connes-Kreimer Hopf algebra of Foissy, the noncommutative
    Faà di Bruno Hopf algebra and its deformations, the noncommutative multi-symmetric
    functions Hopf algebras of Novelli and Thibon, and the double tensor Hopf algebra of
    Ebrahimi-Fard and Patras.
}

\begin{document}

\MakeFirstPage

\section{Introduction}
A polynomial realization of a Hopf algebra consists of interpreting its elements as
polynomials, either commutative or not, in such a way that its product translates as
a polynomial multiplication and its coproduct translates as a simple transformation of the
alphabet of variables. A great portion of combinatorial Hopf algebras appearing in
combinatorics admit polynomial realizations~\cite{DHT02,NT06,Mau13,FNT14,Foi20}. It is
striking to note that Hopf algebras involving a variety of different families of
combinatorial objects and operations on them can be translated and understood in a common
manner through adequate polynomial realizations.

Such realizations are crucial for several reasons. First, they allow us to prove easily that
the Hopf algebra axioms (like associativity, coassociativity, and the Hopf compatibility
between the product and the coproduct) are satisfied. Indeed, if a Hopf algebra admits a
polynomial realization, such properties are almost immediate on polynomials~\cite{Hiv07}.
Second, given a polynomial realization of a Hopf algebra, it is in most cases fruitful to
specialize the associated polynomials (for instance by letting the variables commute) in
order to get Hopf algebra morphisms to other Hopf algebras. This leads to the construction
of new Hopf algebras or establishes links between already existing ones~\cite{FNT14}.
Finally, such families of polynomials realizing a Hopf algebra lead to the definition of new
classes of special functions, analogous to Schur or Macdonald functions~\cite{DHT02}.

Under right conditions, an operad $\Operad$ produces a Hopf algebra $\NaturalHopfAlgebra
\App \Operad$, called the \textit{natural Hopf algebra} of $\Operad$. The bases of
$\NaturalHopfAlgebra \App \Operad$ are indexed by some words on $\Operad$, and its coproduct
is inherited from the composition map of $\Operad$. This construction is considered
in~\cite{vdl04,CL07,BG16}, and we focus here on a noncommutative variation for nonsymmetric
operads, appearing first in~\cite{ML14}. In contrast to many examples of Hopf algebras
having polynomial realizations, none are known for $\NaturalHopfAlgebra \App \Operad$. The
main contribution of this work is to provide a polynomial realization for this family of
Hopf algebras built from nonsymmetric operads. The particularity of our approach is that we
consider a polynomial realization based on variables belonging to alphabets endowed with
several unary and binary relations in order to capture the particularities of the coproduct
of $\NaturalHopfAlgebra \App \Operad$. This approach, using what we call \textit{related
alphabets}, generalizes the previous approaches using totally ordered
alphabets~\cite{DHT02,NT06,Hiv07}, quasi-ordered alphabets~\cite{Foi20}, or alphabets
endowed with a single binary relation~\cite{FNT14,Gir11}.

This work is presented as follows.

In Section~\ref{sec:natural_hopf_algebras}, the main notions about natural Hopf algebras of
operads, related alphabets, and polynomial realizations are provided. We also present
elementary but important definitions concerning free operads and terms that constitute their
elements. We conclude with forests, which are finite sequences of terms.

In Section~\ref{sec:natural_hopf_algebras_free_operads}, we provide some properties of the
natural Hopf algebra $\NaturalHopfAlgebra \App \SetTerms \App \Signature$ of the free operad
$\SetTerms \App \Signature$ generated by a signature $\Signature$. This Hopf algebra is
defined on the linear span of forests decorated on $\Signature$. We offer an alternative
expression for the coproduct of $\NaturalHopfAlgebra \App \SetTerms \App \Signature$
(Proposition~\ref{prop:admissible_sets_coproduct}) which is useful for its subsequent
polynomial realization. We also describe an injection of $\NaturalHopfAlgebra \App \Operad$
into $\NaturalHopfAlgebra \App \SetTerms \App \Signature$ when $\Operad$ is an operad
quotient of $\SetTerms \App \Signature$ subject to specific properties
(Theorem~\ref{thm:quotient_operad_hopf_morphism}).

Section~\ref{sec:polynomial_realization} is the principal part of this work. Here, we
introduce forest-like alphabets, a type of alphabet endowed with certain unary and binary
relations. Next, we define a map $\Realization_A$ that sends any element of
$\NaturalHopfAlgebra \App \SetTerms \App \Signature$ to a polynomial on $A$, where $A$ is a
forest-like alphabet, and we show that it forms a polynomial realization of
$\NaturalHopfAlgebra \App \SetTerms \App \Signature$ (Theorem~\ref{thm:realization}). A key
component to establish this property is a particular forest-like alphabet
$\PositionsAlphabet \App \Signature$ that encodes the shape and the decorations of a forest,
ensuring that $\Realization_{\PositionsAlphabet \App \Signature}$ is injective. We also show
that the previous injection of $\NaturalHopfAlgebra \App \Operad$ into $\NaturalHopfAlgebra
\App \SetTerms \App \Signature$ can be used to obtain a polynomial realization
of~$\NaturalHopfAlgebra \App \Operad$.

In the final Section~\ref{sec:links_hopf_algebras}, we establish links between
$\NaturalHopfAlgebra \App \Operad$ and other Hopf algebras by using the previous polynomial
realization. First, we introduce a generalization of word quasi-symmetric
functions~\cite{Hiv99,NT06} on alphabets with decorated letters, similar to what is
considered in~\cite{NT10} for different classes of polynomials derived from polynomial
realizations of Hopf algebras. We then show that $\NaturalHopfAlgebra \App \SetTerms \App
\Signature$ admits, as a quotient, a space of decorated word quasi-symmetric functions
(Theorem~\ref{thm:map_to_wqsym}). Next, we show that this quotient is isomorphic to a Hopf
subalgebra of a decorated version of the noncommutative Connes-Kreimer Hopf
algebra~\cite{CK98,Foi02a,Foi02b} (Theorem~\ref{thm:map_to_nck}). Finally, we consider
polynomial realizations of natural Hopf algebras of two families of not necessarily free
operads. Among these, as an application of the results of the previous section, we propose
two polynomial realizations of the noncommutative Faà di Bruno Hopf
algebra~\cite{FG05,BFK06,Foi08,Bul11} (Theorem~\ref{thm:n_mas_length_realization}). We also
propose a polynomial realization of the double tensor Hopf algebra~\cite{EP15}.

\paragraph{General notations and conventions}
For any ``$X$-$C$'' concept that depends on an $X$ entity, if in some circumstance $X$ is
either known or insignificant, we may simply denote it by ``$C$''. For instance, we shall
write simply ``forest'' instead of ``$\Signature$-forest'' when $\Signature$ is known or
insignificant. If $f$ is an entity parameterized by an input $x$, we write $f \App x$ for
$f(x)$. By definition, $\App$ associates from right to left so that $f \App g \App x$
denotes $f \App (g \App x)$. Moreover, $\App$ is defined as having higher priority than any
other operator. For a statement $P$, the Iverson bracket $\Iverson{P}$ takes $1$ as value if
$P$ is true and $0$ otherwise. For an integer $i$, $[i]$ (resp.\ $\HanL{i}$) denotes the set
$\{1, \dots, i\}$ (resp.\ $\{0, \dots, i\}$). For a set $A$, $A^*$ is the set of words on
$A$. For $w \in A^*$, $\Length \App w$ is the length of $w$, and for $i \in [\Length \App
w]$, $w \App i$ is the $i$-th letter of $w$. The only word of length $0$ is the empty word
$\epsilon$. For a subset $A'$ of $A$, $w_{|A'}$ is the subword of $w$ made of letters of
$A'$. Let moreover $\Occurrences{A'} \App w$ be the set of positions of the letters of $w$
which belong to $A'$. Given two words $w$ and $w'$, the concatenation of $w$ and $w'$ is
denoted by~$w \Conc w'$ or by $ww'$.

\section{Natural Hopf algebras of operads} \label{sec:natural_hopf_algebras}
In this preliminary section, we recall the concept of a natural Hopf algebra of an operad.
We also review the essential concepts about polynomial realizations of Hopf algebras and
introduce the notion of related alphabet. This section concludes with definitions
concerning terms, forests, and free operads.

In the entirety of this article, all algebraic structures defined on a vector space assume
that this vector space is over a field $\K$ of characteristic zero.

\subsection{Natural Hopf algebras of nonsymmetric operads}
\label{subsec:natural_hopf_algebras}
We provide here some basic definitions about operads and natural Hopf algebras of operads.

\subsubsection{Signatures} \label{subsubsec:signatures}
A \Def{signature} is a set $\Signature$ endowed with a map $\Arity : \Signature \to \N$.
Given $\GenS \in \Signature$, $\Arity \App \GenS$ is the \Def{arity} of $\GenS$. For any $n
\in \N$, let $\Signature \App n := \Bra{\GenS \in \Signature : \Arity \App \GenS = n}$. The
signature $\Signature$ is \Def{positive} if $\Signature \App 0 = \emptyset$. When all
$\Signature \App n$ are finite, the \Def{profile} of $\Signature$ is the infinite word $w$
such that for any $i \geq 1$, $w \App i$ is the cardinality of $\Signature \App \Par{i -
1}$. To write profiles, we shall use the notation $a^\omega$, $a \in \N$, to specify an
infinite sequence of letters $a$. For instance, the infinite word $1020^\omega$ is the
profile of a signature $\Signature$ such that $\# \Signature \App 0 = 1$, $\# \Signature
\App 1 = 0$, $\# \Signature \App 2 = 2$, and $\# \Signature \App n = 0$ for all $n \geq 3$.
A signature $\Signature$ is \Def{binary} if its profile is of the form $00r0^\omega$ where
$r \in \N$.

For the examples that will follow throughout the article, we shall consider the signature
$\SignatureExample := \Bra{\GenA, \GenB, \GenC}$ of profile $01110^\omega$ where $\Arity
\App \GenA = 1$, $\Arity \App \GenB = 2$, and $\Arity \App \GenC = 3$.

\subsubsection{Nonsymmetric operads}
We follow the usual notations about nonsymmetric operads~\cite{Gir18} (called simply
\Def{operads} here). An operad $\Operad$ is above all considered to be a signature. We
denote by
\begin{math}
    -\Han{-, \dots, -} :
    \Operad \App n \times \Operad \App m_1 \times \dots \times \Operad \App m_n
    \to \Operad \App \Par{m_1 + \dots + m_n}
\end{math}
the composition map of $\Operad$ and by $\Unit$ the unit of $\Operad$. The partial
composition map of $\Operad$ is denoted by~$\circ_i$.

Let us introduce two properties an operad $\Operad$ can satisfy. When each $x \in \Operad$
admits finitely many factorizations $x = y \Han{y'_1, \dots, y'_{\Arity \App y}}$ where $y,
y'_1, \dots, y'_{\Arity \App y} \in \Operad$, $\Operad$ is \Def{finitely factorizable}. When
there exists a map $\Deg : \Operad \to \N$ such that $\Deg^{-1} \App 0 = \Bra{\Unit}$ and,
for any $y, y'_1, \dots, y'_{\Arity \App y} \in \Operad$, $\Deg \App \Par{y \Han{y'_1,
\dots, y'_{\Arity \App y}}} = \Deg \App y + \Deg \App y'_1 + \dots + \Deg \App y'_{\Arity
\App y}$, the map $\Deg$ is a \Def{grading} of~$\Operad$.

\subsubsection{Natural Hopf algebras} \label{subsubsec:natural_hopf_algebras}
The \Def{natural Hopf algebra}~\cite{vdl04,CL07,ML14,BG16} of a finitely factorizable operad
$\Operad$ admitting a grading $\Deg$ is the Hopf algebra $\NaturalHopfAlgebra \App \Operad$
defined as follows.  Let $\Reduced : \Operad^* \to \Par{\Operad \setminus \Bra{\Unit}}^*$ be
the map such that $\Reduced \App w$ is the subword of $w \in \Operad^*$ consisting of its
letters different from $\Unit$. A word $w$ on $\Operad$ is \Def{reduced} if $\Reduced \App w
= w$. Let $\NaturalHopfAlgebra \App \Operad$ be the $\K$-linear span of the set $\Reduced
\App \Operad^*$. The bases of $\NaturalHopfAlgebra \App \Operad$ are thus indexed by
$\Reduced \App \Operad^*$, and the \Def{elementary basis} (or \Def{$\BasisE$-basis} for
short) of $\NaturalHopfAlgebra \App \Operad$ is the set $\Bra{\BasisE_w : w \in \Reduced
\App \Operad^*}$. This vector space is endowed with an associative algebra structure through
the product $\Product$ satisfying, for any $w_1, w_2 \in \Reduced \App \Operad^*$,
\begin{equation}
    \BasisE_{w_1} \Product \BasisE_{w_2} = \BasisE_{w_1 \Conc w_2}.
\end{equation}
Moreover, $\NaturalHopfAlgebra \App \Operad$ is endowed with the coproduct $\Coproduct$
defined as the unique associative algebra morphism satisfying, for any $x \in \Operad$,
\begin{equation} \label{equ:natural_coproduct}
    \Coproduct \App \BasisE_x =
    \sum_{y \in \Operad}
    \;
    \sum_{w \in \Operad^{\Arity \App y}}
    \;
    \Iverson{x = y \Han{w \App 1, \dots, w \App \Length \App w}}
    \;
    \BasisE_{\Reduced \App y} \otimes \BasisE_{\Reduced \App w},
\end{equation}
where the outer $\Iverson{-}$ denotes the Iverson bracket. Due to the fact that $\Operad$ is
finitely factorizable, \eqref{equ:natural_coproduct} is a finite sum. This coproduct endows
$\NaturalHopfAlgebra \App \Operad$ with the structure of a bialgebra. By extending
additively $\Deg$ on $\Operad^*$, the map $\Deg$ defines a grading of $\NaturalHopfAlgebra
\App \Operad$. Thus, $\NaturalHopfAlgebra \App \Operad$ admits an antipode and becomes a
Hopf algebra.

\subsubsection{Noncommutative Faà di Bruno Hopf algebra}
\label{subsubsec:faa_di_bruno_hopf_algebra}
An important example of a natural Hopf algebra of an operad is the following. Let us
consider the \Def{associative operad} $\As$, defined by $\As := \Bra{\alpha_n : n \in \N}$,
for any $\alpha_n \in \As$, $\Arity \App \alpha_n := n + 1$, for any $\alpha_n,
\alpha_{m_1}, \dots, \alpha_{m_{n + 1}} \in \As$, $\alpha_n \Han{\alpha_{m_1}, \dots,
\alpha_{m_{n + 1}}} := \alpha_{n + m_1 + \dots + m_{n + 1}}$, and $\Unit := \alpha_0$. The
map $\Deg$ satisfying, for any $\alpha_n \in \As$, $\Deg \App \alpha_n = n$ is a grading of
$\As$. The bases of $\NaturalHopfAlgebra \App \As$ are indexed by $\Reduced \App \As^*$. We
have for instance
\begin{equation}
    \Coproduct \App \BasisE_{\alpha_3} =
    \BasisE_\epsilon \otimes \BasisE_{\alpha_3}
    + 2 \BasisE_{\alpha_1} \otimes \BasisE_{\alpha_2}
    + \BasisE_{\alpha_1} \otimes \BasisE_{\alpha_1 \alpha_1}
    + 3 \BasisE_{\alpha_2} \otimes \BasisE_{\alpha_1}
    + \BasisE_{\alpha_3} \otimes \BasisE_\epsilon.
\end{equation}
It is shown in~\cite{BG16} that $\NaturalHopfAlgebra \App \As$ is isomorphic to the
\Def{noncommutative Faà di Bruno Hopf algebra} $\NCFdB$ (see~\cite{FG05,BFK06,Foi08}).

\subsection{Polynomial realizations} \label{subsec:polynomial_realizations}
We establish here our framework to work with polynomial realizations of Hopf algebras.

\subsubsection{Related alphabets}
A \Def{related alphabet signature} is a positive signature $\RelatedAlphabetSignature$. An
\Def{$\RelatedAlphabetSignature$-related alphabet} is a set $A$ endowed with an $n$-ary
relation $\NAryRelation^A$ for each $\NAryRelation \in \RelatedAlphabetSignature$ where $n =
\Arity \App \NAryRelation$. As usual, an $n$-ary relation $\NAryRelation^A$ on $A$ is a
subset of $A^n$. We denote by $\NAryRelation^A\Par{a_1, \dots, a_n}$ the fact that
$\Par{a_1, \dots, a_n} \in \NAryRelation^A$. When $n = 2$, we write $a_1 \, \NAryRelation^A
\, a_2$ instead of $\NAryRelation^A\Par{a_1, a_2}$.

Let $A_1$ and $A_2$ be two $\RelatedAlphabetSignature$-related alphabets. An
\Def{$\RelatedAlphabetSignature$-related alphabet morphism} is a map $\phi : A_1 \to A_2$
such that for any $\NAryRelation \in \RelatedAlphabetSignature$, by denoting by $n$ the
arity of $ \NAryRelation$, for any $a_1, \dots, a_n \in A_1$, $\NAryRelation^{A_1}\Par{a_1,
\dots, a_n}$ implies $\NAryRelation^{A_2}\Par{\phi \App a_1, \dots, \phi \App a_n}$. An
\Def{$\RelatedAlphabetSignature$-related alphabet congruence} of an
$\RelatedAlphabetSignature$-related alphabet $A$ is an equivalence relation $\Equiv$ on $A$.
For any $a \in A$, we denote by $[a]_{\Equiv}$ the $\Equiv$-equivalence class of $a$. The
\Def{quotient} of $A$ by $\Equiv$ is the $\RelatedAlphabetSignature$-related alphabet $A
/_{\Equiv}$ on the set of $\Equiv$-equivalence classes such that, for any $\NAryRelation \in
\RelatedAlphabetSignature$, by denoting by $n$ the arity of $\NAryRelation$, for any $a_1,
\dots, a_n \in A$, if $\NAryRelation^A\Par{a_1, \dots, a_n}$ then $\NAryRelation^{A
/_{\Equiv}}\Par{\Han{a_1}_{\Equiv}, \dots, \Han{a_n}_{\Equiv}}$.

We shall consider classes of $\RelatedAlphabetSignature$-related alphabets satisfying
possibly some additional conditions. For instance, the \Def{class of totally ordered
alphabets} is the class $\ClassTotallyOrderedAlphabets$ of
$\RelatedAlphabetSignature$-related alphabets where $\RelatedAlphabetSignature$ contains a
binary element $\leq$ and such that for any alphabet $A$ of $\ClassTotallyOrderedAlphabets$,
$\leq^A$ is a total order relation on~$A$.

\subsubsection{Noncommutative polynomials}
For any alphabet $A$, $\K \Angle{A}$ is the $\K$-vector space of \Def{$A$-polynomials},
which are noncommutative polynomials with variables in $A$, having a possibly infinite
support but a finite degree. For instance, for $A := \Bra{\Letter_i : i \in \N}$, the
infinite sum
\begin{equation}
    \sum_{i_1, i_2 \in \N} \Iverson{i_1 \leq i_2} \; \Letter_{i_1} \Letter_{i_2}
\end{equation}
is an element of $\K \Angle{A}$, which is also homogeneous of degree~$2$. In contrast, the
infinite sum $\sum_{n \in \N} \Letter_0^n$ has no finite degree and is not in $\K
\Angle{A}$. The vector space $\K \Angle{A}$ is a graded unital associative algebra for the
usual product of polynomials.

Besides, given two alphabets $A_1$ and $A_2$, let
\begin{math}
    \WordToTensor_{A_1, A_2} : \K \Angle{A_1 \sqcup A_2}
    \to \K \Angle{A_1} \otimes \K \Angle{A_2}
\end{math}
be the linear map such that
\begin{math}
    \WordToTensor_{A_1, A_2} \App w = w_{|A_1} \otimes w_{|A_2}
\end{math}
for any $w \in \Par{A_1 \sqcup A_2}^*$.

\subsubsection{Polynomial realizations}
A \Def{polynomial realization} of a Hopf algebra $\HopfAlgebra$ is a quadruple
$\Par{\ClassAlphabets, \AlphabetSum, \Realization_A, \Alphabet}$ such that
\begin{enumerate}[label=({\sf \roman*})]
    \item \label{item:polynomial_realization_1}
    $\ClassAlphabets$ is a class of related alphabets;
    \item \label{item:polynomial_realization_2}
    $\AlphabetSum$ is an associative operation on $\ClassAlphabets$ which is a disjoint sum
    on the underlying sets of the related alphabets;
    \item \label{item:polynomial_realization_3}
    for any related alphabet $A$ of $\ClassAlphabets$, $\Realization_A : \HopfAlgebra \to \K
    \Angle{A}$ is a graded unital associative algebra morphism;
    \item \label{item:polynomial_realization_4}
    for any alphabets $A_1$ and $A_2$ of $\ClassAlphabets$ and any $x \in \HopfAlgebra$,
    \begin{math}
        \WordToTensor_{A_1, A_2} \App \Realization_{A_1 \AlphabetSum A_2} \App x
        = \Par{\Realization_{A_1} \otimes \Realization_{A_2}} \App \Coproduct \App x;
    \end{math}
    \item \label{item:polynomial_realization_5}
    $\Alphabet$ is an alphabet of $\ClassAlphabets$ such that the map
    $\Realization_\Alphabet$ is injective.
\end{enumerate}
By a slight abuse of terminology, we shall write that $\Realization_\Alphabet$ itself is a
\Def{polynomial realization} when all components $\ClassAlphabets$, $\AlphabetSum$,
$\Realization_A$, and $\Alphabet$ are known.

Property~\ref{item:polynomial_realization_4} is known as the \Def{alphabet doubling trick}.
This property, enjoyed by polynomial realizations of a large number of Hopf algebras, allows
us to rephrase their coproduct via such alphabet
transformations~\cite{DHT02,NT06,Hiv07,Gir11,FNT14,Foi20}.

An interesting point about polynomial realizations is based on the elementary fact that if
$\phi : \Space_1 \to \Space_2$ is a linear map between two vector spaces $\Space_1$ and
$\Space_2$, then $\Space_1 / \Ker \App \phi$ is isomorphic to $\Im \App \phi$, where $\Ker
\App \phi$ is the kernel of $\phi$ and $\Im \App \phi$ is the image of $\phi$. Using this,
if $\Realization_\Alphabet$ is a polynomial realization of a Hopf algebra $\HopfAlgebra$,
any alphabet $A$ of $\ClassAlphabets$ gives rise to a quotient $\HopfAlgebra / \Ker \App
\Realization_A$ which is isomorphic to the subspace $\Im \App \Realization_A$ of $\K
\Angle{A}$.

\subsection{Terms, forests, and free operads}
We end this preliminary section with some combinatorial notions about terms and forests, and
also on free operads.

\subsubsection{Terms} \label{subsubsec:terms}
Let $\Signature$ be a signature. An \Def{$\Signature$-term} is either the \Def{leaf} $\Leaf$
or a pair $\Par{\GenS, \Par{\TermT_1, \dots, \TermT_n}}$ where $n \in \N$, $\GenS \in
\Signature \App n$, and $\TermT_1$, \dots, $\TermT_n$ are $\Signature$-terms. For
convenience, we write $\GenS \Par{\TermT_1, \dots, \TermT_n}$ for $\Par{\GenS,
\Par{\TermT_1, \dots, \TermT_n}}$.  By definition, an $\Signature$-term is therefore a
planar rooted tree where each internal node having $n$ children is decorated on $\Signature
\App n$. The set of $\Signature$-terms is denoted by $\SetTerms \App \Signature$. Let
$\TermT \in \SetTerms \App \Signature$. The \Def{degree} $\Deg \App \TermT$ of $\TermT$ is
the number of internal nodes of $\TermT$. The \Def{arity} $\Arity \App \TermT$ of $\TermT$
is the number of occurrences of leaves of $\TermT$.

For instance, $\GenC \Par{\Leaf, \GenB \Par{\Leaf, \GenA \Par{\Leaf}}, \GenB \Par{\Leaf,
\Leaf}}$ in an $\SignatureExample$-term. This $\SignatureExample$-term $\TermT$ writes as
the planar rooted tree
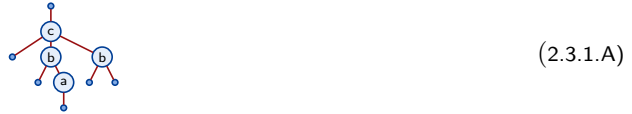
\begin{equation} \label{equ:example_term}
    \scalebox{.75}{
    \begin{tikzpicture}[Centering,xscale=0.45,yscale=0.45]
        \node[Leaf](1)at(1.5,0){};
        \node[Node](2)at(1.5,-1){$\GenC$};
        \node[Leaf](3)at(0,-2){};
        \node[Node](4)at(1.5,-2){$\GenB$};
        \node[Leaf](5)at(1,-3){};
        \node[Node](6)at(2,-3){$\GenA$};
        \node[Leaf](7)at(2,-4){};
        \node[Node](8)at(3.5,-2){$\GenB$};
        \node[Leaf](9)at(3,-3){};
        \node[Leaf](10)at(4,-3){};
        \draw[Edge](1)--(2);
        \draw[Edge](2)--(3);
        \draw[Edge](2)--(4);
        \draw[Edge](2)--(8);
        \draw[Edge](4)--(5);
        \draw[Edge](4)--(6);
        \draw[Edge](6)--(7);
        \draw[Edge](8)--(9);
        \draw[Edge](8)--(10);
    \end{tikzpicture}}
\end{equation}
and is such that $\Deg \App \TermT = 4$ and $\Arity \App \TermT = 5$.

\subsubsection{Free operads}
The \Def{free operad} on a signature $\Signature$ is the set $\SetTerms \App \Signature$
considered as a signature through the arity map $\Arity$, with the composition map such that
for any $\TermT, \TermT'_1, \dots, \TermT'_{\Arity \App \TermT} \in \SetTerms \App
\Signature$, $\TermT \Han{\TermT'_1, \dots, \TermT'_{\Arity \App \TermT}}$ is the
$\Signature$-term obtained by grafting the root of each $\TermT'_i$, $i \in [\Arity \App
\TermT]$, onto the $i$-th leaf of $\TermT$, and with the unit $\Leaf$. Moreover, the partial
composition map of $\SetTerms \App \Signature$ admits the following description. For any
$\TermT, \TermT' \in \SetTerms \App \Signature$ and $i \in [\Arity \App \TermT]$, $\TermT
\circ_i \TermT'$ is the $\Signature$-term obtained by grafting the root of $\TermT'$ onto
the $i$-th leaf of $\TermT$, where the numbering of leaves is from left to right.

For instance, in $\SetTerms \App \SignatureExample$, we have
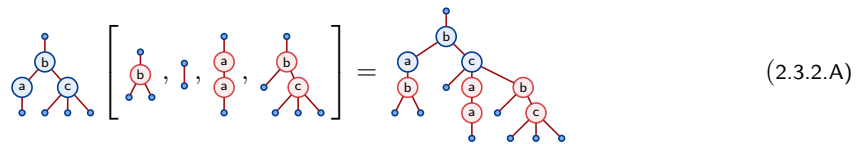
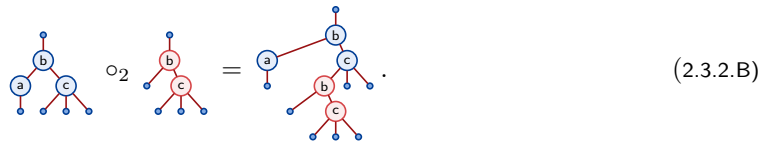
\begin{equation}
    \scalebox{.75}{
        \begin{tikzpicture}[Centering,xscale=0.4,yscale=0.45]
            \node[Leaf](1)at(1,0){};
            \node[Node](2)at(1,-1){$\GenB$};
            \node[Node](3)at(0,-2){$\GenA$};
            \node[Leaf](4)at(0,-3){};
            \node[Node](5)at(2,-2){$\GenC$};
            \node[Leaf](6)at(1,-3){};
            \node[Leaf](7)at(2,-3){};
            \node[Leaf](8)at(3,-3){};
            \draw[Edge](1)--(2);
            \draw[Edge](2)--(3);
            \draw[Edge](2)--(5);
            \draw[Edge](3)--(4);
            \draw[Edge](5)--(6);
            \draw[Edge](5)--(7);
            \draw[Edge](5)--(8);
        \end{tikzpicture}
    }
    \Han{
        \scalebox{.75}{
            \begin{tikzpicture}[Centering,scale=0.4]
                \node[Leaf](1)at(0.5,0){};
                \node[Node,MarkB](2)at(0.5,-1){$\GenB$};
                \node[Leaf](3)at(0,-2){};
                \node[Leaf](4)at(1,-2){};
                \draw[Edge](1)--(2);
                \draw[Edge](2)--(3);
                \draw[Edge](2)--(4);
            \end{tikzpicture}
        },
        \scalebox{.75}{
            \begin{tikzpicture}[Centering,scale=0.4]
                \node[Leaf](1)at(0,0){};
                \node[Leaf](2)at(0,-1){};
                \draw[Edge](1)--(2);
            \end{tikzpicture}
        },
        \scalebox{.75}{
            \begin{tikzpicture}[Centering,scale=0.45]
                \node[Leaf](1)at(0,0){};
                \node[Node,MarkB](2)at(0,-1){$\GenA$};
                \node[Node,MarkB](3)at(0,-2){$\GenA$};
                \node[Leaf](4)at(0,-3){};
                \draw[Edge](1)--(2);
                \draw[Edge](2)--(3);
                \draw[Edge](3)--(4);
            \end{tikzpicture}
        },
        \scalebox{.75}{
            \begin{tikzpicture}[Centering,xscale=0.4,yscale=0.45]
                \node[Leaf](1)at(1.5,0){};
                \node[Node,MarkB](2)at(1.5,-1){$\GenB$};
                \node[Leaf](3)at(0.5,-2){};
                \node[Node,MarkB](4)at(2,-2){$\GenC$};
                \node[Leaf](5)at(1,-3){};
                \node[Leaf](6)at(2,-3){};
                \node[Leaf](7)at(3,-3){};
                \draw[Edge](1)--(2);
                \draw[Edge](2)--(3);
                \draw[Edge](2)--(4);
                \draw[Edge](4)--(5);
                \draw[Edge](4)--(6);
                \draw[Edge](4)--(7);
            \end{tikzpicture}
        }
    }
    =
    \scalebox{.75}{
        \begin{tikzpicture}[Centering,xscale=0.45,yscale=0.45]
            \node[Leaf](1)at(2,0){};
            \node[Node](2)at(2,-1){$\GenB$};
            \node[Node](3)at(0.5,-2){$\GenA$};
            \node[Node,MarkB](4)at(0.5,-3){$\GenB$};
            \node[Leaf](5)at(0,-4){};
            \node[Leaf](6)at(1,-4){};
            \node[Node](7)at(3.0,-2){$\GenC$};
            \node[Leaf](8)at(2,-3){};
            \node[Node,MarkB](9)at(3,-3){$\GenA$};
            \node[Node,MarkB](10)at(3,-4){$\GenA$};
            \node[Leaf](11)at(3,-5){};
            \node[Node,MarkB](12)at(5.0,-3){$\GenB$};
            \node[Leaf](13)at(4.0,-4){};
            \node[Node,MarkB](14)at(5.5,-4){$\GenC$};
            \node[Leaf](15)at(4.5,-5){};
            \node[Leaf](16)at(5.5,-5){};
            \node[Leaf](17)at(6.5,-5){};
            \draw[Edge](1)--(2);
            \draw[Edge](2)--(3);
            \draw[Edge](2)--(7);
            \draw[Edge](3)--(4);
            \draw[Edge](4)--(5);
            \draw[Edge](4)--(6);
            \draw[Edge](7)--(8);
            \draw[Edge](7)--(9);
            \draw[Edge](7)--(12);
            \draw[Edge](9)--(10);
            \draw[Edge](10)--(11);
            \draw[Edge](12)--(13);
            \draw[Edge](12)--(14);
            \draw[Edge](14)--(15);
            \draw[Edge](14)--(16);
            \draw[Edge](14)--(17);
        \end{tikzpicture}
    }
\end{equation}
and
\begin{equation}
    \scalebox{.75}{
        \begin{tikzpicture}[Centering,xscale=0.4,yscale=0.45]
            \node[Leaf](1)at(1,0){};
            \node[Node](2)at(1,-1){$\GenB$};
            \node[Node](3)at(0,-2){$\GenA$};
            \node[Leaf](4)at(0,-3){};
            \node[Node](5)at(2,-2){$\GenC$};
            \node[Leaf](6)at(1,-3){};
            \node[Leaf](7)at(2,-3){};
            \node[Leaf](8)at(3,-3){};
            \draw[Edge](1)--(2);
            \draw[Edge](2)--(3);
            \draw[Edge](2)--(5);
            \draw[Edge](3)--(4);
            \draw[Edge](5)--(6);
            \draw[Edge](5)--(7);
            \draw[Edge](5)--(8);
        \end{tikzpicture}
    }
    \circ_2
    \scalebox{.75}{
        \begin{tikzpicture}[Centering,xscale=0.4,yscale=0.45]
            \node[Leaf](1)at(1.5,0){};
            \node[Node,MarkB](2)at(1.5,-1){$\GenB$};
            \node[Leaf](3)at(0.5,-2){};
            \node[Node,MarkB](4)at(2,-2){$\GenC$};
            \node[Leaf](5)at(1,-3){};
            \node[Leaf](6)at(2,-3){};
            \node[Leaf](7)at(3,-3){};
            \draw[Edge](1)--(2);
            \draw[Edge](2)--(3);
            \draw[Edge](2)--(4);
            \draw[Edge](4)--(5);
            \draw[Edge](4)--(6);
            \draw[Edge](4)--(7);
        \end{tikzpicture}
    }
    =
    \scalebox{.75}{
        \begin{tikzpicture}[Centering,xscale=0.4,yscale=0.45]
            \node[Leaf](1)at(3,0){};
            \node[Node](2)at(3,-1){$\GenB$};
            \node[Node](3)at(0,-2){$\GenA$};
            \node[Leaf](4)at(0,-3){};
            \node[Node](5)at(3.5,-2){$\GenC$};
            \node[Node,MarkB](6)at(2.5,-3){$\GenB$};
            \node[Leaf](7)at(1,-4){};
            \node[Node,MarkB](8)at(3,-4){$\GenC$};
            \node[Leaf](9)at(2,-5){};
            \node[Leaf](10)at(3,-5){};
            \node[Leaf](11)at(4,-5){};
            \node[Leaf](12)at(3.5,-3){};
            \node[Leaf](13)at(4.5,-3){};
            \draw[Edge](1)--(2);
            \draw[Edge](2)--(3);
            \draw[Edge](2)--(5);
            \draw[Edge](3)--(4);
            \draw[Edge](5)--(6);
            \draw[Edge](5)--(12);
            \draw[Edge](5)--(13);
            \draw[Edge](6)--(7);
            \draw[Edge](6)--(8);
            \draw[Edge](8)--(9);
            \draw[Edge](8)--(10);
            \draw[Edge](8)--(11);
        \end{tikzpicture}
    }.
\end{equation}

Observe that $\SetTerms \App \Signature$ is finitely factorizable and that the map $\Deg$ is
a grading of~$\SetTerms \App \Signature$.

\subsubsection{Forests} \label{subsubsec:forests}
Let $\Signature$ be a signature. An \Def{$\Signature$-forest} is any element of $\SetForests
\App \Signature := {\SetTerms \App \Signature}^*$. Let $\ForestF \in \SetForests \App
\Signature$. The \Def{degree} $\Deg \App \ForestF$ of
$\ForestF$ is the sum of the degrees of the terms forming it. The \Def{arity} $\Arity \App
\ForestF$ of $\ForestF$ is the sum of the arities of the terms forming it. We identify each
internal node of an $\Signature$-forest $\ForestF$ with its position, starting by $1$,
according to the left to right preorder traversal of $\ForestF$. The \Def{decoration map} of
$\ForestF$ is the map $\Decoration_\ForestF$ sending each internal node of $\ForestF$ to its
decoration. The \Def{height map} of $\ForestF$ is the map $\Height_\ForestF$ sending each
internal node $i$ of $\ForestF$ to the length of the path connecting the $i$ to the root of
the $\Signature$-term to which $i$ belongs. In particular, if $i$ is a root of $\ForestF$,
then $\Height_\ForestF \App i = 0$. Let moreover for any $j \geq 1$ the binary relation
$\Edge{\ForestF}{j}$ on the set of internal nodes of $\ForestF$ such that $i
\Edge{\ForestF}{j} i'$ holds if $i'$ is the $j$-th child of $i$ in $\ForestF$. Let $i$ be an
internal node of $\ForestF$ and $i_0$ be the root of the $\Signature$-term to which $i$
belongs. The \Def{position} of $i$ in $\ForestF$ is the word $\Position_\ForestF \App i =
j_1 j_2 \dots j_{\ell - 1} j_\ell$ provided that
\begin{math}
    i_0 \Edge{\ForestF}{j_1} i_1 \Edge{\ForestF}{j_2} \dots
    \Edge{\ForestF}{j_{\ell - 1}} i_{\ell - 1}
    \Edge{\ForestF}{j_\ell} i
\end{math}
for some internal nodes $i_1, \dots, i_{\ell - 1}$ of $\ForestF$ and positive integers $j_1,
j_2, \dots, j_{\ell - 1}, j_\ell$.

Let us give some examples of the previous notions. For this, let
\begin{equation} \label{equ:example_forest}
    \ForestF :=
    \scalebox{.75}{
    \begin{tikzpicture}[Centering,xscale=0.55,yscale=0.5]
        \node[Leaf](2)at(0,-1){};
        \node[Leaf](3)at(0,-2){};
        \node[Leaf](4)at(2.5,-1){};
        \node[Node](5)at(2.5,-2){$\GenC$};
        \node[Node](6)at(1.5,-3){$\GenA$};
        \node[Leaf](7)at(1.5,-4){};
        \node[Leaf](8)at(2.5,-3){};
        \node[Node](9)at(3.5,-3){$\GenB$};
        \node[Node](10)at(3,-4){$\GenA$};
        \node[Leaf](11)at(3,-5){};
        \node[Leaf](12)at(4,-4){};
        \node[Leaf](13)at(5,-1){};
        \node[Leaf](14)at(5,-2){};
        \node[Leaf](15)at(6,-1){};
        \node[Leaf](16)at(6,-2){};
        \node[Leaf](17)at(8,-1){};
        \node[Node](18)at(8,-2){$\GenB$};
        \node[Leaf](19)at(7,-3){};
        \node[Node](20)at(8.5,-3){$\GenB$};
        \node[Node](21)at(8,-4){$\GenA$};
        \node[Leaf](22)at(8,-5){};
        \node[Leaf](23)at(9,-4){};
        \draw[Edge](2)--(3);
        \draw[Edge](4)--(5);
        \draw[Edge](5)--(6);
        \draw[Edge](5)--(8);
        \draw[Edge](5)--(9);
        \draw[Edge](6)--(7);
        \draw[Edge](9)--(10);
        \draw[Edge](9)--(12);
        \draw[Edge](10)--(11);
        \draw[Edge](13)--(14);
        \draw[Edge](15)--(16);
        \draw[Edge](17)--(18);
        \draw[Edge](18)--(19);
        \draw[Edge](18)--(20);
        \draw[Edge](20)--(21);
        \draw[Edge](20)--(23);
        \draw[Edge](21)--(22);
        \node[left of=5,node distance=10pt,font=\small,text=ColB]{$1$};
        \node[left of=6,node distance=10pt,font=\small,text=ColB]{$2$};
        \node[right of=9,node distance=10pt,font=\small,text=ColB]{$3$};
        \node[left of=10,node distance=10pt,font=\small,text=ColB]{$4$};
        \node[right of=18,node distance=10pt,font=\small,text=ColB]{$5$};
        \node[right of=20,node distance=10pt,font=\small,text=ColB]{$6$};
        \node[left of=21,node distance=10pt,font=\small,text=ColB]{$7$};
    \end{tikzpicture}}
\end{equation}
be an $\SignatureExample$-forest where internal nodes are identified by the integers close
to them. The degree of $\ForestF$ is $7$, its arity is $10$, and we have for instance
$\Decoration_\ForestF \App 1 = \GenC$, $\Decoration_\ForestF \App 3 = \GenB$, $1
\Edge{\ForestF}{1} 2$, $1 \Edge{\ForestF}{3} 3$, and $5 \Edge{\ForestF}{2} 6$. Moreover, we
have $\Position_\ForestF \App 1 = \epsilon$, $\Position_\ForestF \App 5 = \epsilon$,
$\Position_\ForestF \App 4 = 31$, and $\Position_\ForestF \App 7 = 21$. Notice that the
internal nodes $1$ and $5$ of $\ForestF$ have the same position.

\section{Natural Hopf algebras of free operads}
\label{sec:natural_hopf_algebras_free_operads}
In this section, we begin by focusing specifically on the natural Hopf algebras of free
operads and derive some of their properties. Next, we consider the natural Hopf algebras of
not necessarily free operads and interpret them as Hopf subalgebras of natural Hopf
algebras of free operads.

\subsection{First properties}
Here, we give a description of the basis elements of $\NaturalHopfAlgebra \App \SetTerms
\App \Signature$ where $\Signature$ is a signature, its Hilbert series, and necessary and
sufficient conditions for the commutativity and the cocommutativity of this Hopf algebra. We
end this section by providing an alternative description for its coproduct, useful later to
establish a polynomial realization of $\NaturalHopfAlgebra \App \SetTerms \App \Signature$.

\subsubsection{Basis elements, product, and coproduct}
By construction, for any signature $\Signature$, the Hopf algebra $\NaturalHopfAlgebra \App
\SetTerms \App \Signature$ is graded by $\Deg$ and its bases are indexed by the set
$\Reduced \App \SetForests \App \Signature$ of reduced $\Signature$-forests. By definition,
a reduced $\Signature$-forest is a word on the alphabet~$\SetTerms \App \Signature \setminus
\Bra{\Leaf}$.

On the $\BasisE$-basis, the product of $\NaturalHopfAlgebra \App \SetTerms \App \Signature$
works by concatenating the reduced forests. For instance, in
$\NaturalHopfAlgebra \App \SetTerms \App \SignatureExample$,
\begin{equation}
    \BasisE_{
        \scalebox{0.75}{
            \begin{tikzpicture}[Centering,xscale=0.5,yscale=0.45]
                \node[Leaf](2)at(0.5,-1){};
                \node[Node](3)at(0.5,-2){$\GenA$};
                \node[Node](4)at(0.5,-3){$\GenB$};
                \node[Leaf](5)at(0,-4){};
                \node[Leaf](6)at(1,-4){};
                \node[Leaf](7)at(3,-1){};
                \node[Node](8)at(3,-2){$\GenC$};
                \node[Node](9)at(2,-3){$\GenA$};
                \node[Leaf](10)at(2,-4){};
                \node[Leaf](11)at(3,-3){};
                \node[Leaf](12)at(4,-3){};
                \draw[Edge](2)--(3);
                \draw[Edge](3)--(4);
                \draw[Edge](4)--(5);
                \draw[Edge](4)--(6);
                \draw[Edge](7)--(8);
                \draw[Edge](8)--(9);
                \draw[Edge](8)--(11);
                \draw[Edge](8)--(12);
                \draw[Edge](9)--(10);
            \end{tikzpicture}
        }
    }
    \Product
    \BasisE_{
        \scalebox{0.75}{
        \begin{tikzpicture}[Centering,xscale=0.5,yscale=0.4]
            \node[Leaf](1)at(0.5,0){};
            \node[Node,MarkB](2)at(0.5,-1){$\GenB$};
            \node[Node,MarkB](3)at(0,-2){$\GenA$};
            \node[Leaf](4)at(0,-3){};
            \node[Leaf](5)at(1,-2){};
            \draw[Edge](1)--(2);
            \draw[Edge](2)--(3);
            \draw[Edge](2)--(5);
            \draw[Edge](3)--(4);
        \end{tikzpicture}}}
    =
    \BasisE_{
        \scalebox{0.75}{
        \begin{tikzpicture}[Centering,xscale=0.45,yscale=0.45]
            \node[Leaf](2)at(0.5,-1){};
            \node[Node](3)at(0.5,-2){$\GenA$};
            \node[Node](4)at(0.5,-3){$\GenB$};
            \node[Leaf](5)at(0,-4){};
            \node[Leaf](6)at(1,-4){};
            \node[Leaf](7)at(3,-1){};
            \node[Node](8)at(3,-2){$\GenC$};
            \node[Node](9)at(2,-3){$\GenA$};
            \node[Leaf](10)at(2,-4){};
            \node[Leaf](11)at(3,-3){};
            \node[Leaf](12)at(4,-3){};
            \node[Leaf](13)at(5.5,-1){};
            \node[Node,MarkB](14)at(5.5,-2){$\GenB$};
            \node[Node,MarkB](15)at(5,-3){$\GenA$};
            \node[Leaf](16)at(5,-4){};
            \node[Leaf](17)at(6,-3){};
            \draw[Edge](2)--(3);
            \draw[Edge](3)--(4);
            \draw[Edge](4)--(5);
            \draw[Edge](4)--(6);
            \draw[Edge](7)--(8);
            \draw[Edge](8)--(9);
            \draw[Edge](8)--(11);
            \draw[Edge](8)--(12);
            \draw[Edge](9)--(10);
            \draw[Edge](13)--(14);
            \draw[Edge](14)--(15);
            \draw[Edge](14)--(17);
            \draw[Edge](15)--(16);
        \end{tikzpicture}
    }
}.
\end{equation}
On the other hand, on the $\BasisE$-basis the coproduct of $\NaturalHopfAlgebra \App
\SetTerms \App \Signature$ works by summing on all ways to cut a reduced $\Signature$-forest
into a upper part and a lower part, and then, by removing the obtained $\Signature$-terms
equal to the leaf in both parts. For instance, in $\NaturalHopfAlgebra \App \SetTerms \App
\SignatureExample$,
\begin{align}
    \Coproduct \App
    \BasisE_{
        \scalebox{0.75}{
        \begin{tikzpicture}[Centering,xscale=0.4,yscale=0.45]
            \node[Leaf](2)at(1,-1){};
            \node[Node](3)at(1,-2){$\GenC$};
            \node[Leaf](4)at(0,-3){};
            \node[Node](5)at(1,-3){$\GenA$};
            \node[Leaf](6)at(1,-4){};
            \node[Leaf](7)at(2,-3){};
            \node[Leaf](8)at(3.5,-1){};
            \node[Node](9)at(3.5,-2){$\GenB$}; 
            \node[Leaf](10)at(3,-3){};
            \node[Leaf](11)at(4,-3){};
            \draw[Edge](2)--(3);
            \draw[Edge](3)--(4);
            \draw[Edge](3)--(5);
            \draw[Edge](3)--(7);
            \draw[Edge](5)--(6);
            \draw[Edge](8)--(9);
            \draw[Edge](9)--(10);
            \draw[Edge](9)--(11);
        \end{tikzpicture}}}
    & =
    \BasisE_\epsilon
    \otimes
    \BasisE_{
        \scalebox{0.75}{
        \begin{tikzpicture}[Centering,xscale=0.4,yscale=0.45]
            \node[Leaf](2)at(1,-1){};
            \node[Node](3)at(1,-2){$\GenC$};
            \node[Leaf](4)at(0,-3){};
            \node[Node](5)at(1,-3){$\GenA$};
            \node[Leaf](6)at(1,-4){};
            \node[Leaf](7)at(2,-3){};
            \node[Leaf](8)at(3.5,-1){};
            \node[Node](9)at(3.5,-2){$\GenB$}; 
            \node[Leaf](10)at(3,-3){};
            \node[Leaf](11)at(4,-3){};
            \draw[Edge](2)--(3);
            \draw[Edge](3)--(4);
            \draw[Edge](3)--(5);
            \draw[Edge](3)--(7);
            \draw[Edge](5)--(6);
            \draw[Edge](8)--(9);
            \draw[Edge](9)--(10);
            \draw[Edge](9)--(11);
        \end{tikzpicture}}}
    +
    \BasisE_{
        \scalebox{0.75}{
        \begin{tikzpicture}[Centering,scale=0.4]
            \node[Leaf](1)at(1,0){};
            \node[Node](2)at(1,-1){$\GenC$};
            \node[Leaf](3)at(0,-2){};
            \node[Leaf](4)at(1,-2){};
            \node[Leaf](5)at(2,-2){};
            \draw[Edge](1)--(2);
            \draw[Edge](2)--(3);
            \draw[Edge](2)--(4);
            \draw[Edge](2)--(5);
        \end{tikzpicture}}}
    \otimes
    \BasisE_{
        \scalebox{0.75}{
        \begin{tikzpicture}[Centering,scale=0.4]
            \node[Leaf](2)at(0,-1){};
            \node[Node](3)at(0,-2){$\GenA$};
            \node[Leaf](4)at(0,-3){};
            \node[Leaf](5)at(1.5,-1){};
            \node[Node](6)at(1.5,-2){$\GenB$};
            \node[Leaf](7)at(1,-3){};
            \node[Leaf](8)at(2,-3){};
            \draw[Edge](2)--(3);
            \draw[Edge](3)--(4);
            \draw[Edge](5)--(6);
            \draw[Edge](6)--(7);
            \draw[Edge](6)--(8);
        \end{tikzpicture}}}
    +
    \BasisE_{
        \scalebox{0.75}{
        \begin{tikzpicture}[Centering,scale=0.4]
            \node[Leaf](1)at(0.5,0){};
            \node[Node](2)at(0.5,-1){$\GenB$};
            \node[Leaf](3)at(0,-2){};
            \node[Leaf](4)at(1,-2){};
            \draw[Edge](1)--(2);
            \draw[Edge](2)--(3);
            \draw[Edge](2)--(4);
        \end{tikzpicture}}}
    \otimes
    \BasisE_{
        \scalebox{0.75}{
        \begin{tikzpicture}[Centering,xscale=0.4,yscale=0.45]
            \node[Leaf](1)at(1,0){};
            \node[Node](2)at(1,-1){$\GenC$};
            \node[Leaf](3)at(0,-2){};
            \node[Node](4)at(1,-2){$\GenA$};
            \node[Leaf](5)at(1,-3){};
            \node[Leaf](6)at(2,-2){};
            \draw[Edge](1)--(2);
            \draw[Edge](2)--(3);
            \draw[Edge](2)--(4);
            \draw[Edge](2)--(6);
            \draw[Edge](4)--(5);
        \end{tikzpicture}}}
    \\
    & \quad +
    \BasisE_{
        \scalebox{0.75}{
        \begin{tikzpicture}[Centering,xscale=0.4,yscale=0.45]
            \node[Leaf](1)at(1,0){};
            \node[Node](2)at(1,-1){$\GenC$};
            \node[Leaf](3)at(0,-2){};
            \node[Node](4)at(1,-2){$\GenA$};
            \node[Leaf](5)at(1,-3){};
            \node[Leaf](6)at(2,-2){};
            \draw[Edge](1)--(2);
            \draw[Edge](2)--(3);
            \draw[Edge](2)--(4);
            \draw[Edge](2)--(6);
            \draw[Edge](4)--(5);
        \end{tikzpicture}}}
    \otimes
    \BasisE_{
        \scalebox{0.75}{
        \begin{tikzpicture}[Centering,scale=0.4]
            \node[Leaf](1)at(0.5,0){};
            \node[Node](2)at(0.5,-1){$\GenB$};
            \node[Leaf](3)at(0,-2){};
            \node[Leaf](4)at(1,-2){};
            \draw[Edge](1)--(2);
            \draw[Edge](2)--(3);
            \draw[Edge](2)--(4);
        \end{tikzpicture}}}
    +
    \BasisE_{
        \scalebox{0.75}{
        \begin{tikzpicture}[Centering,scale=0.4]
            \node[Leaf](2)at(1,-1){};
            \node[Node](3)at(1,-2){$\GenC$};
            \node[Leaf](4)at(0,-3){};
            \node[Leaf](5)at(1,-3){};
            \node[Leaf](6)at(2,-3){};
            \node[Leaf](7)at(3.5,-1){};
            \node[Node](8)at(3.5,-2){$\GenB$};
            \node[Leaf](9)at(3,-3){};
            \node[Leaf](10)at(4,-3){};
            \draw[Edge](2)--(3);
            \draw[Edge](3)--(4);
            \draw[Edge](3)--(5);
            \draw[Edge](3)--(6);
            \draw[Edge](7)--(8);
            \draw[Edge](8)--(9);
            \draw[Edge](8)--(10);
        \end{tikzpicture}}}
    \otimes
    \BasisE_{
        \scalebox{0.75}{
        \begin{tikzpicture}[Centering,scale=0.4]
            \node[Leaf](1)at(0,0){};
            \node[Node](2)at(0,-1){$\GenA$};
            \node[Leaf](3)at(0,-2){};
            \draw[Edge](1)--(2);
            \draw[Edge](2)--(3);
        \end{tikzpicture}}}
    +
    \BasisE_{
        \scalebox{0.75}{
        \begin{tikzpicture}[Centering,xscale=0.4,yscale=0.45]
            \node[Leaf](2)at(1,-1){};
            \node[Node](3)at(1,-2){$\GenC$};
            \node[Leaf](4)at(0,-3){};
            \node[Node](5)at(1,-3){$\GenA$};
            \node[Leaf](6)at(1,-4){};
            \node[Leaf](7)at(2,-3){};
            \node[Leaf](8)at(3.5,-1){};
            \node[Node](9)at(3.5,-2){$\GenB$}; 
            \node[Leaf](10)at(3,-3){};
            \node[Leaf](11)at(4,-3){};
            \draw[Edge](2)--(3);
            \draw[Edge](3)--(4);
            \draw[Edge](3)--(5);
            \draw[Edge](3)--(7);
            \draw[Edge](5)--(6);
            \draw[Edge](8)--(9);
            \draw[Edge](9)--(10);
            \draw[Edge](9)--(11);
        \end{tikzpicture}}}
    \otimes
    \BasisE_\epsilon.
    \notag
\end{align}

\subsubsection{Hilbert series}
Here is a description of the Hilbert series $\HilbertSeries_{\NaturalHopfAlgebra \App
\SetTerms \App \Signature}$ of $\NaturalHopfAlgebra \App \SetTerms \App \Signature$ when
$\Signature$ is a finite signature.

\begin{Statement}{Proposition}{prop:hilbert_series_natural_hopf_algebra_free_operad}
    For any finite signature $\Signature$, the Hilbert series of $\NaturalHopfAlgebra \App
    \SetTerms \App \Signature$ satisfies
    \begin{equation}
        \HilbertSeries_{\NaturalHopfAlgebra \App \SetTerms \App \Signature}
        = \frac{1}{2 - T}
    \end{equation}
    where $T$ is the generating series satisfying
    \begin{math}
        T = 1 + z \; S \Han{z := T},
    \end{math}
    and $S$ is the polynomial defined by
    \begin{math}
        S := \sum_{n \in \N} \# \Signature \App n \; z^n.
    \end{math}
\end{Statement}
\begin{Proof}
    From the functional equation defining $T$, this series is the generating series of
    $\Signature$-terms enumerated w.r.t.\ their degrees. By construction, since
    $\NaturalHopfAlgebra \App \SetTerms \App \Signature$ is the linear span of reduced
    $\Signature$-forests, its Hilbert series is $1 / (1 - T')$, where $T'$ is the generating
    series of the $\Signature$-terms different from the leaf. Since $T' = T - 1$, the
    statement of the proposition follows.
\end{Proof}

\subsubsection{Commutativity and cocommutativity}
Here is a necessary and sufficient condition for the commutativity and cocommutativity of
$\NaturalHopfAlgebra \App \SetTerms \App \Signature$ for any signature $\Signature$.

\begin{Statement}{Proposition}
{prop:commutativity_cocommutativity_natural_hopf_algebra_free_operad}
    Let $\Signature$ be a signature of profile $w$. The Hopf algebra $\NaturalHopfAlgebra
    \App \SetTerms \App \Signature$ is
    \begin{enumerate}[label=({\sf \roman*})]
        \item \label{item:commutativity_cocommutativity_natural_hopf_algebra_free_operad_1}
        commutative if and only if $w = 0^\omega$ or $w = 10^\omega$;
        \item \label{item:commutativity_cocommutativity_natural_hopf_algebra_free_operad_2}
        cocommutative if and only if $w = k0^\omega$, $k \in \N$, or $w = 010^\omega$.
    \end{enumerate}
\end{Statement}
\begin{Proof}
    Assume that $w = 0^\omega$. In this case, each reduced $\Signature$-forest is
    necessarily empty. Hence, $\NaturalHopfAlgebra \App \SetTerms \App \Signature$ is the
    linear span of $\Bra{\BasisE_\epsilon}$ and is commutative. When $w = 10^\omega$, each
    reduced $\Signature$-forest is a concatenation of $\Signature$-terms which are
    themselves made of a single node of arity $0$ and decorated by the same $\GenS \in
    \Signature \App 0$. From the definition of the product $\Product$ of
    $\NaturalHopfAlgebra \App \SetTerms \App \Signature$, it follows that
    $\NaturalHopfAlgebra \App \SetTerms \App \Signature$ is commutative. Conversely, if
    $\Signature$ is such that $\# \Signature \App 0 \geq 2$ or $\# \Signature \App n \geq 1$
    for an $n \geq 1$, then it is possible to build two different $\Signature$-terms
    $\TermT_1$ and $\TermT_2$ which are both different from the leaf. Since in this case,
    \begin{math}
        \BasisE_{\TermT_1} \Product \BasisE_{\TermT_2}
        = \BasisE_{\TermT_1 \Conc \TermT_2}
        \ne \BasisE_{\TermT_2 \Conc \TermT_1}
        = \BasisE_{\TermT_2} \Product \BasisE_{\TermT_1},
    \end{math}
    it follows that $\NaturalHopfAlgebra \App \SetTerms \App \Signature$ is not commutative
    and proves~\ref{item:commutativity_cocommutativity_natural_hopf_algebra_free_operad_1}.

    Assume that $w = k0^\omega$, $k \in \N$. When $k = 0$, each reduced $\Signature$-forest
    is necessarily empty. In this case, $\NaturalHopfAlgebra \App \SetTerms \App \Signature$
    is the linear span of $\Bra{\BasisE_\epsilon}$ and is cocommutative. When $w =
    k0^\omega$, $k \geq 1$, each $\Signature$-term $\TermT$ is made of a single node of
    arity $0$ and decorated by an $\GenS \in \Signature \App 0$. From the definition of the
    coproduct $\Coproduct$ of $\NaturalHopfAlgebra \App \SetTerms \App \Signature$,
    $\BasisE_\TermT$ is a primitive element, implying that this Hopf algebra is
    cocommutative. Assume that $w = 010^\omega$. In this case, each $\Signature$-term
    $\TermT$ is of the form $\GenS \Par{\GenS \Par{ \dots \GenS \Par{\Leaf} \dots}}$ where
    $\GenS \in \Signature \App 1$. Moreover, if there are two such $\Signature$-terms
    $\TermT_1$ and $\TermT_2$ such that $\TermT = \TermT_1 \Han{\TermT_2}$, then we have
    also $\TermT = \TermT_2 \Han{\TermT_1}$. From the definition of the coproduct
    $\Coproduct$ of $\NaturalHopfAlgebra \App \SetTerms \App \Signature$, this implies that
    this Hopf algebra is cocommutative. Conversely, assume that $w \ne k0^\omega$, $k \in
    \N$, and $w \ne 010^\omega$. By an elementary logical reasoning, we obtain that the
    signature $\Signature$ can take exactly three different forms, leading to the following
    cases.
    \begin{enumerate}[label=({\sf \arabic*})]
        \item If $\# \Signature \App 0 \geq 1$ and $\# \Signature \App 1 = 1$, let $\GenS_0
        \in \Signature \App 0$ and $\GenS_1 \in \Signature \App 1$. We have
        \begin{equation}
            \Coproduct \App
            \BasisE_{
                \scalebox{0.75}{
                    \begin{tikzpicture}[Centering,xscale=0.55,yscale=0.45]
                        \node[Leaf](1)at(0,0){};
                        \node[Node](2)at(0,-1){$\GenS_1$};
                        \node[Node](3)at(0,-2){$\GenS_0$};
                        \draw[Edge](1)--(2);
                        \draw[Edge](2)--(3);
                    \end{tikzpicture}
                }
            }
            =
            \BasisE_\epsilon
            \otimes
            \BasisE_{
                \scalebox{0.75}{
                    \begin{tikzpicture}[Centering,xscale=0.55,yscale=0.45]
                        \node[Leaf](1)at(0,0){};
                        \node[Node](2)at(0,-1){$\GenS_1$};
                        \node[Node](3)at(0,-2){$\GenS_0$};
                        \draw[Edge](1)--(2);
                        \draw[Edge](2)--(3);
                    \end{tikzpicture}
                }
            }
            +
            \BasisE_{
                \scalebox{0.75}{
                    \begin{tikzpicture}[Centering,xscale=0.55,yscale=0.45]
                        \node[Leaf](1)at(0,0){};
                        \node[Node](2)at(0,-1){$\GenS_1$};
                        \draw[Edge](1)--(2);
                    \end{tikzpicture}
                }
            }
            \otimes
            \BasisE_{
                \scalebox{0.75}{
                    \begin{tikzpicture}[Centering,xscale=0.55,yscale=0.45]
                        \node[Leaf](1)at(0,0){};
                        \node[Node](2)at(0,-1){$\GenS_0$};
                        \draw[Edge](1)--(2);
                    \end{tikzpicture}
                }
            }
            +
            \BasisE_{
                \scalebox{0.75}{
                    \begin{tikzpicture}[Centering,xscale=0.55,yscale=0.45]
                        \node[Leaf](1)at(0,0){};
                        \node[Node](2)at(0,-1){$\GenS_1$};
                        \node[Node](3)at(0,-2){$\GenS_0$};
                        \draw[Edge](1)--(2);
                        \draw[Edge](2)--(3);
                    \end{tikzpicture}
                }
            }
            \otimes
            \BasisE_\epsilon.
        \end{equation}
        \item If $\# \Signature \App 1 \geq 2$, let $\GenS, \GenS' \in \Signature \App 1$
        such that $\GenS \ne \GenS'$. We have
        \begin{equation}
            \Coproduct \App
            \BasisE_{
                \scalebox{0.75}{
                    \begin{tikzpicture}[Centering,xscale=0.55,yscale=0.45]
                        \node[Leaf](1)at(0,0){};
                        \node[Node](2)at(0,-1){$\GenS$};
                        \node[Node](3)at(0,-2){$\GenS'$};
                        \node[Leaf](4)at(0,-3){};
                        \draw[Edge](1)--(2);
                        \draw[Edge](2)--(3);
                        \draw[Edge](3)--(4);
                    \end{tikzpicture}
                }
            }
            =
            \BasisE_\epsilon
            \otimes
            \BasisE_{
                \scalebox{0.75}{
                    \begin{tikzpicture}[Centering,xscale=0.55,yscale=0.45]
                        \node[Leaf](1)at(0,0){};
                        \node[Node](2)at(0,-1){$\GenS$};
                        \node[Node](3)at(0,-2){$\GenS'$};
                        \node[Leaf](4)at(0,-3){};
                        \draw[Edge](1)--(2);
                        \draw[Edge](2)--(3);
                        \draw[Edge](3)--(4);
                    \end{tikzpicture}
                }
            }
            +
            \BasisE_{
                \scalebox{0.75}{
                    \begin{tikzpicture}[Centering,xscale=0.55,yscale=0.45]
                        \node[Leaf](1)at(0,0){};
                        \node[Node](2)at(0,-1){$\GenS$};
                        \node[Leaf](3)at(0,-2){};
                        \draw[Edge](1)--(2);
                        \draw[Edge](2)--(3);
                    \end{tikzpicture}
                }
            }
            \otimes
            \BasisE_{
                \scalebox{0.75}{
                    \begin{tikzpicture}[Centering,xscale=0.55,yscale=0.45]
                        \node[Leaf](1)at(0,0){};
                        \node[Node](2)at(0,-1){$\GenS'$};
                        \node[Leaf](3)at(0,-2){};
                        \draw[Edge](1)--(2);
                        \draw[Edge](2)--(3);
                    \end{tikzpicture}
                }
            }
            +
            \BasisE_{
                \scalebox{0.75}{
                    \begin{tikzpicture}[Centering,xscale=0.55,yscale=0.45]
                        \node[Leaf](1)at(0,0){};
                        \node[Node](2)at(0,-1){$\GenS$};
                        \node[Node](3)at(0,-2){$\GenS'$};
                        \node[Leaf](4)at(0,-3){};
                        \draw[Edge](1)--(2);
                        \draw[Edge](2)--(3);
                        \draw[Edge](3)--(4);
                    \end{tikzpicture}
                }
            }
            \otimes
            \BasisE_\epsilon.
        \end{equation}
        \item Otherwise, there is an $n \geq 2$ such that $\# \Signature \App n \geq 1$. By
        taking $\GenS \in \Signature \App n$, we have
        \begin{align}
            \Coproduct \App
            \BasisE_{
                \scalebox{0.75}{
                    \begin{tikzpicture}[Centering,xscale=0.5,yscale=0.5]
                        \node[Leaf](1)at(3.25,0){};
                        \node[Node](2)at(3.25,-1){$\GenS$};
                        \node[Node](3)at(1,-2){$\GenS$};
                        \node[Leaf](4)at(0.5,-3){};
                        \node[Leaf](5)at(1.5,-3){};
                        \node[Node](6)at(2.5,-2){$\GenS$};
                        \node[Leaf](7)at(2,-3){};
                        \node[Leaf](8)at(3,-3){};
                        \node[Leaf](9)at(4,-2){};
                        \node[Leaf](10)at(5,-2){};
                        \draw[Edge](1)--(2);
                        \draw[Edge](2)--(3);
                        \draw[Edge](2)--(6);
                        \draw[Edge](2)--(9);
                        \draw[Edge](2)--(10);
                        \draw[Edge](3)--(4);
                        \draw[Edge](3)--(5);
                        \draw[Edge](6)--(7);
                        \draw[Edge](6)--(8);
                        \node[right of=4,font=\scriptsize,node distance=0.275cm]{$\dots$};
                        \node[right of=7,font=\scriptsize,node distance=0.275cm]{$\dots$};
                        \node[right of=9,font=\scriptsize,node distance=0.275cm]{$\dots$};
                    \end{tikzpicture}
                }
            }
            & =
            \BasisE_\epsilon
            \otimes
            \BasisE_{
                \scalebox{0.75}{
                    \begin{tikzpicture}[Centering,xscale=0.5,yscale=0.5]
                        \node[Leaf](1)at(3.25,0){};
                        \node[Node](2)at(3.25,-1){$\GenS$};
                        \node[Node](3)at(1,-2){$\GenS$};
                        \node[Leaf](4)at(0.5,-3){};
                        \node[Leaf](5)at(1.5,-3){};
                        \node[Node](6)at(2.5,-2){$\GenS$};
                        \node[Leaf](7)at(2,-3){};
                        \node[Leaf](8)at(3,-3){};
                        \node[Leaf](9)at(4,-2){};
                        \node[Leaf](10)at(5,-2){};
                        \draw[Edge](1)--(2);
                        \draw[Edge](2)--(3);
                        \draw[Edge](2)--(6);
                        \draw[Edge](2)--(9);
                        \draw[Edge](2)--(10);
                        \draw[Edge](3)--(4);
                        \draw[Edge](3)--(5);
                        \draw[Edge](6)--(7);
                        \draw[Edge](6)--(8);
                        \node[right of=4,font=\scriptsize,node distance=0.275cm]{$\dots$};
                        \node[right of=7,font=\scriptsize,node distance=0.275cm]{$\dots$};
                        \node[right of=9,font=\scriptsize,node distance=0.275cm]{$\dots$};
                    \end{tikzpicture}
                }
            }
            +
            \BasisE_{
                \scalebox{0.75}{
                    \begin{tikzpicture}[Centering,xscale=0.55,yscale=0.45]
                        \node[Leaf](1)at(0.5,0){};
                        \node[Node](2)at(0.5,-1){$\GenS$};
                        \node[Leaf](3)at(0,-2){};
                        \node[Leaf](4)at(1,-2){};
                        \draw[Edge](1)--(2);
                        \draw[Edge](2)--(3);
                        \draw[Edge](2)--(4);
                        \node[right of=3,font=\scriptsize,node distance=0.275cm]{$\dots$};
                    \end{tikzpicture}
                }
            }
            \otimes
            \BasisE_{
                \scalebox{0.75}{
                    \begin{tikzpicture}[Centering,xscale=0.55,yscale=0.45]
                        \node[Leaf](2)at(0.5,-1){};
                        \node[Node](3)at(0.5,-2){$\GenS$};
                        \node[Leaf](4)at(0,-3){};
                        \node[Leaf](5)at(1,-3){};
                        \node[Leaf](6)at(2.5,-1){};
                        \node[Node](7)at(2.5,-2){$\GenS$};
                        \node[Leaf](8)at(2,-3){};
                        \node[Leaf](9)at(3,-3){};
                        \draw[Edge](2)--(3);
                        \draw[Edge](3)--(4);
                        \draw[Edge](3)--(5);
                        \draw[Edge](6)--(7);
                        \draw[Edge](7)--(8);
                        \draw[Edge](7)--(9);
                        \node[right of=4,font=\scriptsize,node distance=0.275cm]{$\dots$};
                        \node[right of=8,font=\scriptsize,node distance=0.275cm]{$\dots$};
                    \end{tikzpicture}
                }
            }
            +
            \BasisE_{
                \scalebox{0.75}{
                    \begin{tikzpicture}[Centering,xscale=0.5,yscale=0.5]
                        \node[Leaf](1)at(2,0){};
                        \node[Node](2)at(2,-1){$\GenS$};
                        \node[Node](3)at(1,-2){$\GenS$};
                        \node[Leaf](4)at(0.5,-3){};
                        \node[Leaf](5)at(1.5,-3){};
                        \node[Leaf](6)at(2,-2){};
                        \node[Leaf](7)at(3,-2){};
                        \draw[Edge](1)--(2);
                        \draw[Edge](2)--(3);
                        \draw[Edge](2)--(6);
                        \draw[Edge](2)--(7);
                        \draw[Edge](3)--(4);
                        \draw[Edge](3)--(5);
                        \node[right of=4,font=\scriptsize,node distance=0.275cm]{$\dots$};
                        \node[right of=6,font=\scriptsize,node distance=0.275cm]{$\dots$};
                    \end{tikzpicture}
                }
            }
            \otimes
            \BasisE_{
                \scalebox{0.75}{
                    \begin{tikzpicture}[Centering,xscale=0.55,yscale=0.45]
                        \node[Leaf](1)at(0.5,0){};
                        \node[Node](2)at(0.5,-1){$\GenS$};
                        \node[Leaf](3)at(0,-2){};
                        \node[Leaf](4)at(1,-2){};
                        \draw[Edge](1)--(2);
                        \draw[Edge](2)--(3);
                        \draw[Edge](2)--(4);
                        \node[right of=3,font=\scriptsize,node distance=0.275cm]{$\dots$};
                    \end{tikzpicture}
                }
            }
            \notag
            \\
            & \quad
            +
            \BasisE_{
                \scalebox{0.75}{
                    \begin{tikzpicture}[Centering,xscale=0.5,yscale=0.5]
                        \node[Leaf](1)at(2.25,0){};
                        \node[Node](2)at(2.25,-1){$\GenS$};
                        \node[Leaf](3)at(0,-2){};
                        \node[Node](4)at(1.5,-2){$\GenS$};
                        \node[Leaf](5)at(1,-3){};
                        \node[Leaf](6)at(2,-3){};
                        \node[Leaf](7)at(3,-2){};
                        \node[Leaf](8)at(4,-2){};
                        \draw[Edge](1)--(2);
                        \draw[Edge](2)--(3);
                        \draw[Edge](2)--(4);
                        \draw[Edge](2)--(7);
                        \draw[Edge](2)--(8);
                        \draw[Edge](4)--(5);
                        \draw[Edge](4)--(6);
                        \node[right of=5,font=\scriptsize,node distance=0.275cm]{$\dots$};
                        \node[right of=7,font=\scriptsize,node distance=0.275cm]{$\dots$};
                    \end{tikzpicture}
                }
            }
            \otimes
            \BasisE_{
                \scalebox{0.75}{
                    \begin{tikzpicture}[Centering,xscale=0.55,yscale=0.45]
                        \node[Leaf](1)at(0.5,0){};
                        \node[Node](2)at(0.5,-1){$\GenS$};
                        \node[Leaf](3)at(0,-2){};
                        \node[Leaf](4)at(1,-2){};
                        \draw[Edge](1)--(2);
                        \draw[Edge](2)--(3);
                        \draw[Edge](2)--(4);
                        \node[right of=3,font=\scriptsize,node distance=0.275cm]{$\dots$};
                    \end{tikzpicture}
                }
            }
            +
            \BasisE_{
                \scalebox{0.75}{
                    \begin{tikzpicture}[Centering,xscale=0.5,yscale=0.5]
                        \node[Leaf](1)at(3.25,0){};
                        \node[Node](2)at(3.25,-1){$\GenS$};
                        \node[Node](3)at(1,-2){$\GenS$};
                        \node[Leaf](4)at(0.5,-3){};
                        \node[Leaf](5)at(1.5,-3){};
                        \node[Node](6)at(2.5,-2){$\GenS$};
                        \node[Leaf](7)at(2,-3){};
                        \node[Leaf](8)at(3,-3){};
                        \node[Leaf](9)at(4,-2){};
                        \node[Leaf](10)at(5,-2){};
                        \draw[Edge](1)--(2);
                        \draw[Edge](2)--(3);
                        \draw[Edge](2)--(6);
                        \draw[Edge](2)--(9);
                        \draw[Edge](2)--(10);
                        \draw[Edge](3)--(4);
                        \draw[Edge](3)--(5);
                        \draw[Edge](6)--(7);
                        \draw[Edge](6)--(8);
                        \node[right of=4,font=\scriptsize,node distance=0.275cm]{$\dots$};
                        \node[right of=7,font=\scriptsize,node distance=0.275cm]{$\dots$};
                        \node[right of=9,font=\scriptsize,node distance=0.275cm]{$\dots$};
                    \end{tikzpicture}
                }
            }
            \otimes
            \BasisE_\epsilon.
        \end{align}
    \end{enumerate}
    In all these cases, we observe that $\Coproduct$ is not cocommutative,
    proving~\ref{item:commutativity_cocommutativity_natural_hopf_algebra_free_operad_2}.
\end{Proof}

By using
Proposition~\ref{prop:commutativity_cocommutativity_natural_hopf_algebra_free_operad}, we
have exactly the following possibilities:
\begin{enumerate}[label=({\sf \arabic*})]
    \item $\NaturalHopfAlgebra \App \SetTerms \App \Signature$ is commutative and
    cocommutative. This happens if and only if the profile of $\Signature$ is~$k0^\omega$
    with $k \leq 1$;
    \item $\NaturalHopfAlgebra \App \SetTerms \App \Signature$ is noncommutative and
    cocommutative. This happens if and only the profile of $\Signature$ is $k0^\omega$ with
    $k \geq 2$ or is~$010^\omega$;
    \item $\NaturalHopfAlgebra \App \SetTerms \App \Signature$ is
    noncommutative and non-cocommutative. This happens for all other possible profiles
    of $\Signature$.
\end{enumerate}
There is no natural Hopf algebra of a free operad which is commutative and
non-cocommutative.

\subsubsection{Some examples}
We give three classes of examples of natural Hopf algebras of free operads and present some
of their properties by using
Propositions~\ref{prop:hilbert_series_natural_hopf_algebra_free_operad}
and~\ref{prop:commutativity_cocommutativity_natural_hopf_algebra_free_operad}. Let
$\Signature$ be a finite signature of profile $w$.
\begin{enumerate}[label=({\sf \arabic*})]
    \item If $w = k0^\omega$, $k \geq 1$, then $\NaturalHopfAlgebra \App \SetTerms \App
    \Signature$ is the free associative algebra on $k$ generators, endowed with the
    unshuffling cocommutative coproduct.
    \item If $w = 0k0^\omega$, $k \geq 1$, then the bases of $\NaturalHopfAlgebra \App
    \SetTerms \App \Signature$ are indexed by the set of words whose letters are themselves
    nonempty words on $[k]$. For any $n \geq 1$, the dimension of the $n$-th homogeneous
    component of $\NaturalHopfAlgebra \App \SetTerms \App \Signature$ is $k^n 2^{n - 1}$.
    Moreover, when $k = 1$, $\NaturalHopfAlgebra \App \SetTerms \App \Signature$ is the Hopf
    algebra of noncommutative symmetric functions $\NCSF$~\cite{GKLLRT95}. Observe that
    $\NaturalHopfAlgebra \App \SetTerms \App \Signature$ is cocommutative only if $k = 1$.
    \item If $w = 00k0^\omega$, $k \geq 1$, then the bases of $\NaturalHopfAlgebra \App
    \SetTerms \App \Signature$ are indexed by the set of words whose letters are binary
    trees different from the leaf and such that internal nodes are decorated on $[k]$. For
    any $n \geq 1$, the dimension of the $n$-th homogeneous component of
    $\NaturalHopfAlgebra \App \SetTerms \App \Signature$ is $k^n \binom{2n - 1}{n }$. The
    particular case for $k = 1$ is studied in~\cite[Chapter 6, Section 3]{Gir11}.
\end{enumerate}

\subsubsection{Coproduct description} \label{subsubsec:coproduct_description}
We now introduce an alternative way to describe the coproduct of $\NaturalHopfAlgebra \App
\SetTerms \App \Signature$ for any signature $\Signature$.

Let $\ForestF$ be a reduced $\Signature$-forest. Given a set $I$ of nodes of $\ForestF$, the
\Def{restriction} $\ForestF \App I $of $\ForestF$ on $I$ is the reduced $\Signature$-forest
obtained by keeping only the nodes of
$\ForestF$ which are in $I$ and their adjacent edges. For instance,
\begin{equation} \label{equ:forest_restriction_example}
    \scalebox{0.75}{
        \begin{tikzpicture}[Centering,xscale=0.4,yscale=0.45]
            \node[Leaf](2)at(1.25,-1){};
            \node[Node,MarkB](3)at(1.25,-2){$\GenB$};
            \node[Node,MarkB](4)at(0,-3){$\GenA$};
            \node[Leaf](5)at(0,-4){};
            \node[Node](6)at(2.5,-3){$\GenC$};
            \node[Leaf](7)at(1,-4){};
            \node[Node,MarkB](8)at(2.5,-4){$\GenB$};
            \node[Leaf](9)at(2,-5){};
            \node[Leaf](10)at(3,-5){};
            \node[Leaf](11)at(4,-4){};
            \node[Leaf](12)at(6.5,-1){};
            \node[Node](13)at(6.5,-2){$\GenA$};
            \node[Node](14)at(6.5,-3){$\GenC$};
            \node[Leaf](15)at(5,-4){};
            \node[Leaf](16)at(6,-4){};
            \node[Node,MarkB](17)at(7.5,-4){$\GenB$};
            \node[Leaf](18)at(7,-5){};
            \node[Leaf](19)at(8,-5){};
            \draw[Edge](2)--(3);
            \draw[Edge](3)--(4);
            \draw[Edge](3)--(6);
            \draw[Edge](4)--(5);
            \draw[Edge](6)--(7);
            \draw[Edge](6)--(8);
            \draw[Edge](6)--(11);
            \draw[Edge](8)--(9);
            \draw[Edge](8)--(10);
            \draw[Edge](12)--(13);
            \draw[Edge](13)--(14);
            \draw[Edge](14)--(15);
            \draw[Edge](14)--(16);
            \draw[Edge](14)--(17);
            \draw[Edge](17)--(18);
            \draw[Edge](17)--(19);
        \end{tikzpicture}
    }
    \App \Bra{1, 2, 4, 7}
    =
    \scalebox{0.75}{
        \begin{tikzpicture}[Centering,scale=0.4]
            \node[Leaf](2)at(0.5,-1){};
            \node[Node,MarkB](3)at(0.5,-2){$\GenB$};
            \node[Node,MarkB](4)at(0,-3){$\GenA$};
            \node[Leaf](5)at(0,-4){};
            \node[Leaf](6)at(1,-3){};
            \node[Leaf](7)at(2.5,-1){};
            \node[Node,MarkB](8)at(2.5,-2){$\GenB$};
            \node[Leaf](9)at(2,-3){};
            \node[Leaf](10)at(3,-3){};
            \node[Leaf](11)at(4.5,-1){};
            \node[Node,MarkB](12)at(4.5,-2){$\GenB$};
            \node[Leaf](13)at(4,-3){};
            \node[Leaf](14)at(5,-3){};
            \draw[Edge](2)--(3);
            \draw[Edge](3)--(4);
            \draw[Edge](3)--(6);
            \draw[Edge](4)--(5);
            \draw[Edge](7)--(8);
            \draw[Edge](8)--(9);
            \draw[Edge](8)--(10);
            \draw[Edge](11)--(12);
            \draw[Edge](12)--(13);
            \draw[Edge](12)--(14);
        \end{tikzpicture}
    }.
\end{equation}
As a particular case, observe that $\ForestF \App \emptyset = \epsilon$. A pair $\Par{I_1,
I_2}$ of sets is \Def{$\ForestF$-admissible} if $I_1 \sqcup I_2 = [\Deg \App \ForestF]$, for
any $i_1 \in I_1$, all ancestors of the internal node $i_1$ of $\ForestF$ belong to $I_1$,
and for any $i_2 \in I_2$, all descendants of the internal node $i_2$ of $\ForestF$ belong
to $I_2$. This property is denoted by $\Par{I_1, I_2} \AdmissibleSet \ForestF$. For
instance, by denoting by $\ForestF$ the reduced forest of the left-hand side
of~\eqref{equ:forest_restriction_example}, the pair $\Par{\Bra{1, 2, 4, 7}, \Bra{3, 5, 6}}$
is not $\ForestF$-admissible, while the pair $\Par{\Bra{1, 3, 5}, \Bra{2, 4, 6, 7}}$ is.

\begin{Statement}{Proposition}{prop:admissible_sets_coproduct}
    For any signature $\Signature$ and any reduced $\Signature$-forest $\ForestF$, the
    coproduct of $\NaturalHopfAlgebra \App \SetTerms \App \Signature$ satisfies
    \begin{equation} \label{equ:admissible_sets_coproduct}
        \Coproduct \App \BasisE_\ForestF
        = \sum_{I_1, I_2 \subseteq [\Deg \App \ForestF]}
        \Iverson{\Par{I_1, I_2} \AdmissibleSet \ForestF} \;
        \BasisE_{\ForestF \App I_1} \otimes \BasisE_{\ForestF \App I_2}.
    \end{equation}
\end{Statement}
\begin{Proof}
    Let us denote by $\Coproduct'$ the coproduct defined
    by~\eqref{equ:admissible_sets_coproduct}. Let $\TermT \in \SetTerms \App \Signature$ and
    $\Par{I_1, I_2}$ be a $\TermT$-admissible pair of sets. By the fact that $I_1$ (resp.\
    $I_2$) is closed w.r.t.\ the ancestor (resp.\ descendant) relation of $\TermT$, and by
    the definition of the composition map of free operads, this last property is equivalent
    to the fact that $\TermT$ decomposes as $\TermT = \TermT' \Han{\TermT'_1, \dots,
    \TermT'_\ell}$ with $\ell \in \N$, $\TermT', \TermT'_1, \dots, \TermT'_\ell \in
    \SetTerms \App \Signature$, $\Reduced \App \TermT' = \TermT \App I_1$, and $\Reduced
    \App \Par{\TermT'_1 \Conc \dots \Conc \TermT'_\ell} = \TermT \App I_2$. This shows that
    $\Coproduct \App \BasisE_\TermT = \Coproduct' \App \BasisE_\TermT$ and hence, shows that
    $\Coproduct$ and $\Coproduct'$ coincide on the elements of the $\BasisE$-basis indexed
    by $\Signature$-terms.

    Now, let $\ForestF, \ForestF' \in \Reduced \App \SetForests \App \Signature$ and
    $\Par{I_1, I_2}$ be a pair of sets such that $\Par{I_1, I_2} \AdmissibleSet \ForestF
    \Conc \ForestF'$. This is equivalent to the fact that there exist a unique partition
    $\Bra{I'_1, I''_1}$ of $I_1$ and a unique partition $\Bra{I'_2, I''_2}$ of $I_2$ such
    that $\Par{I'_1, I'_2} \AdmissibleSet \ForestF$ and $\Par{I'''_1, I'''_2} \AdmissibleSet
    \ForestF'$, where $I'''_1$ and $I'''_2$ are the sets obtained by respectively
    decrementing by $\Deg \App \ForestF$ each element of $I''_1$ and $I''_2$. This
    observation leads to the fact that $\Coproduct'$ is a morphism of associative algebras.
    Moreover, as shown before, $\Coproduct$ and $\Coproduct'$ coincide on the elements of
    the $\BasisE$-basis indexed by $\Signature$-terms. Since these elements are the
    algebraic generators of $\NaturalHopfAlgebra \App \SetTerms \App \Signature$, the
    coproducts $\Coproduct$ and $\Coproduct'$ are the same.
\end{Proof}

\subsection{Quotients of free operads and Hopf subalgebras}
\label{subsec:quotients_hopf_subalgebras}
We show here that under some conditions, the natural Hopf algebra of an operad can be
realized as a Hopf subalgebra of the natural Hopf algebra of a free operad.

\subsubsection{Equivalence relations on reduced forests}
Let $\Signature$ be a signature and $\Equiv$ be an operad congruence of the free operad
$\SetTerms \App \Signature$. Let us denote by $\pi_{\Equiv} : \SetTerms \App \Signature \to
{\SetTerms \App \Signature} /_{\Equiv}$ the canonical projection map associated with
$\Equiv$. The map $\pi_{\Equiv}$ is extended as a map from $\SetForests \App \Signature$ to
${\SetTerms \App \Signature /_{\Equiv}}^*$ by setting $\pi_{\equiv} \App \ForestF :=
\Par{\pi_{\Equiv} \App \ForestF \App 1} \dots \Par{\pi_{\Equiv} \App \ForestF \App \Length
\App \ForestF}$ for any $\ForestF \in \SetForests \App \Signature$.

The following two properties on operad congruences play an important role here. The operad
congruence $\Equiv$ is \Def{compatible with the degree} if $\TermT \Equiv \TermT'$ implies
$\Deg \App \TermT = \Deg \App \TermT'$ for any $\TermT, \TermT' \in \SetTerms \App
\Signature$. The operad congruence $\Equiv$ is of \Def{finite type} if for any $\TermT \in
\SetTerms \App \Signature$, the $\Equiv$-equivalence class $[\TermT]_{\Equiv}$ of $\TermT$
is finite.

\begin{Statement}{Proposition}{prop:finitely_factorizable_graded_quotient}
    Let $\Signature$ be a signature and $\Equiv$ be an operad congruence of $\SetTerms \App
    \Signature$ compatible with the degree and of finite type. The quotient operad
    $\SetTerms \App \Signature /_{\Equiv}$ is finitely factorizable. Moreover, the map
    $\Deg$ sending each $\Equiv$-equivalence class $[\TermT]_{\Equiv}$ to the degree of any
    $\Signature$-term belonging to it is a grading of $\SetTerms \App \Signature
    /_{\Equiv}$.
\end{Statement}
\begin{Proof}
    Let $x \in \SetTerms \App \Signature /_{\Equiv}$. Since $\Equiv$ is of finite type,
    there are finitely many $\Signature$-terms $\TermT$ such that $\pi_{\Equiv} \App \TermT
    = x$. Moreover, due to the fact that $\SetTerms \App \Signature$ is finitely
    factorizable, the number of pairs $\Par{\TermT', \Par{\TermT'_1, \dots, \TermT'_{\Arity
    \App \TermT'}}} \in \SetTerms \App \Signature \times \Par{\SetTerms \App \Signature}^*$
    such that $\TermT = \TermT' \Han{\TermT'_1, \dots, \TermT'_{\Arity \App \TermT'}}$ is
    finite. Therefore, $x$ admits finitely many factorizations $x = y \Han{y'_1, \dots,
    y'_{\Arity \App y}}$ where $y, y'_1, \dots, y'_{\Arity \App y} \in \SetTerms \App
    \Signature /_{\Equiv}$. This shows that $\SetTerms \App \Signature /_{\Equiv}$ is
    finitely factorizable.

    Finally, since $\Equiv$ is compatible with the degree, for any $x \in \SetTerms \App
    \Signature /_{\Equiv}$, all $\TermT \in \SetTerms \App \Signature$ such that
    $\pi_{\Equiv} \App \TermT = x$ have the same degree as $\Signature$-terms. Since by
    definition of $\Deg$, $\Deg \App x = \Deg \App \TermT$ where $\TermT$ is any
    $\Signature$-term satisfying $\pi_{\Equiv} \App \TermT = x$, $\Deg$ is a grading of
    $\SetTerms \App \Signature /_{\Equiv}$.
\end{Proof}

\subsubsection{Hopf subalgebras of natural Hopf algebras of a free operad}
\label{subsubsec:hopf_subalgebras_free_operads}
Let $\Signature$ be a signature and $\Equiv$ be an operad congruence of $\SetTerms \App
\Signature$ of finite type. Let $\phi : \NaturalHopfAlgebra \App \Par{\SetTerms \App
\Signature /_{\Equiv}} \to \NaturalHopfAlgebra \App \SetTerms \App \Signature$ be the linear
map defined, for any $x \in \Reduced \App \Par{{\SetTerms \App \Signature /_{\Equiv}}^*}$,
by
\begin{equation} \label{equ:injection_from_quotient_operad}
    \phi \App \BasisE_x
    :=
    \sum_{\ForestF \in \Reduced \App \SetForests \App \Signature}
    \Iverson{\pi_{\Equiv} \App \ForestF = x} \; \BasisE_\ForestF.
\end{equation}
Due to the fact that $\Equiv$ is of finite type, $\phi$ is a well-defined linear map.
Moreover, observe that when $\Equiv$ is compatible with the degree, $\phi \App \BasisE_x$ is
a homogeneous element of $\NaturalHopfAlgebra \App \SetTerms \App \Signature$.

\begin{Statement}{Theorem}{thm:quotient_operad_hopf_morphism}
    Let $\Signature$ be a signature and $\Equiv$ be an operad congruence of $\SetTerms \App
    \Signature$ compatible with the degree and of finite type. The map $\phi$ is an
    injective Hopf algebra morphism.
\end{Statement}
\begin{Proof}
    First, since $\Equiv$ is compatible with the degrees and is of finite type, by
    Proposition~\ref{prop:finitely_factorizable_graded_quotient}, the operad $\SetTerms \App
    \Signature /_{\Equiv}$ is finitely factorizable and graded. Hence, $\NaturalHopfAlgebra
    \App \Par{\SetTerms \App \Signature /_{\Equiv}}$ is a well-defined Hopf algebra.

    Observe that by definition of the extension of $\pi_{\Equiv}$ on $\SetForests \App
    \Signature$, for any $\ForestF \in \SetForests \App \Signature$, $\Length \App \ForestF
    = \Length \App \pi_{\Equiv} \App \ForestF$. For this reason, $\phi \App \BasisE_\epsilon
    = \BasisE_\epsilon$. Moreover, for any $x_1, x_2 \in \Reduced \App \Par{{\SetTerms \App
    \Signature /_{\Equiv}}^*}$, the fact that
    \begin{math}
        \phi \App \Par{\BasisE_{x_1} \Product \BasisE_{x_2}}
        = \phi \App \BasisE_{x_1} \Product \phi \App \BasisE_{x_2}
    \end{math}
    follows from a straightforward computation. Therefore, $\phi$ is a unital associative
    algebra morphism.

    Let us show that $\phi$ is a coalgebra morphism. For any  $x \in \SetTerms \App
    \Signature /_{\Equiv}$, we have
    \begin{equation} \label{equ:quotient_operad_hopf_morphism_1}
        \Coproduct \App \phi \App \BasisE_x
        =
        \sum_{\substack{
            \TermT \in \SetTerms \App \Signature \\
            \TermT_1, \dots, \TermT_{\Arity \App \TermT} \in \SetTerms \App \Signature}}
        \Iverson{
            \pi_{\Equiv}
            \App \Par{\TermT \Han{\TermT_1, \dots, \TermT_{\Arity \App \TermT}}}
            = x
        }
        \;
        \BasisE_{\Reduced \App \TermT}
        \otimes
        \BasisE_{\Reduced \App \TermT_1 \dots \TermT_{\Arity \App \TermT}}.
    \end{equation}
    Since $\Equiv$ is an operad congruence, the right-hand side
    of~\eqref{equ:quotient_operad_hopf_morphism_1} rewrites as
    \begin{equation} \label{equ:quotient_operad_hopf_morphism_2}
        \sum_{\substack{
            y \in \SetTerms \App \Signature /_{\Equiv} \\
            y_1, \dots, y_{\Arity \App y} \in \SetTerms \App \Signature /_{\Equiv}}}
        \Iverson{y \Han{y_1, \dots, y_{\Arity \App y}} = x}
        \sum_{\ForestF, \ForestF' \in \SetForests \App \Signature} \;
        \Iverson{\pi_{\Equiv} \App \ForestF = y}
        \Iverson{\pi_{\Equiv} \App \ForestF' = y_i \dots y_{\Arity \App y}}
        \;
        \BasisE_{\Reduced \App \ForestF} \otimes \BasisE_{\Reduced \App \ForestF'}.
    \end{equation}
    Now, observe that since $\Equiv$ is compatible with the degree, $\Han{\Leaf}_{\Equiv} =
    \Bra{\Leaf}$. By the definition of the extension of $\pi_{\Equiv}$ on $\SetForests
    \App \Signature$, this leads to the fact that for any $\ForestF \in \SetForests \App
    \Signature$ and $z \in {\SetTerms \App \Signature /_{\Equiv}}^*$, $\pi_{\Equiv} \App
    \ForestF = z$ implies $\pi_{\Equiv} \App \Reduced \App \ForestF = \Reduced \App z$.
    Moreover, for the same reason, for any $\ForestF \in \Reduced \App \SetForests \App
    \Signature$ and $z \in {\SetTerms \App \Signature /_{\Equiv}}^*$, $\pi_{\Equiv} \App
    \ForestF = \Reduced \App z$ implies that there exists a unique $\ForestF' \in
    \SetForests \App \Signature$ such that $\Reduced \App \ForestF' = \ForestF$ and
    $\pi_{\Equiv} \App \ForestF' = z$. For these reasons,
    \eqref{equ:quotient_operad_hopf_morphism_2} is equal to
    \begin{equation} \label{equ:quotient_operad_hopf_morphism_3}
        \sum_{\substack{
            y \in \SetTerms \App \Signature /_{\Equiv} \\
            y_1, \dots, y_{\Arity \App y} \in \SetTerms \App \Signature /_{\Equiv}}}
        \Iverson{y \Han{y_1, \dots, y_{\Arity \App y}} = x}
        \sum_{\ForestF, \ForestF' \in \Reduced \App \SetForests \App \Signature}
        \Iverson{\pi_{\Equiv} \App \ForestF = \Reduced \App y}
        \Iverson{
            \pi_{\Equiv} \App \ForestF'
            =
            \Reduced \App y_1 \dots y_{\Arity \App y}
        }
        \;
        \BasisE_\ForestF \otimes \BasisE_{\ForestF'}.
    \end{equation}
    It follows now by a straightforward computation
    that~\eqref{equ:quotient_operad_hopf_morphism_3} is equal to $\Par{\phi \otimes \phi}
    \App \Coproduct \App \BasisE_x$. We have shown that for any $x \in \SetTerms \App
    \Signature /_{\Equiv}$, $\Coproduct \App \phi \App \BasisE_x = \Par{\phi \otimes \phi}
    \App \Coproduct \App \BasisE_x$. Now, the fact that $\Coproduct$ and $\phi$ are
    associative algebra morphisms implies that for any $z \in {\SetTerms \App \Signature
    /_{\Equiv}}^*$, $\Coproduct \App \phi \App \BasisE_z = \Par{\phi \otimes \phi} \App
    \Coproduct \App \BasisE_z$. This shows that $\phi$ is a coalgebra morphism.

    Finally, $\phi$ is injective because for any $\ForestF \in \Reduced \App \SetForests
    \App \Signature$, there is exactly one $x \in \Reduced \App \Par{{\SetTerms \App
    \Signature /_{\Equiv}}^*}$ such that $\BasisE_\ForestF$ appears in $\phi \App
    \BasisE_x$. This establishes the statement of the theorem.
\end{Proof}

A consequence of Theorem~\ref{thm:quotient_operad_hopf_morphism} is that, for any operad
congruence $\Equiv$ of $\SetTerms \App \Signature$ compatible with the degree and of finite
type, the Hopf algebra $\NaturalHopfAlgebra \App \Par{\SetTerms \App \Signature /_{\Equiv}}$
can be realized as a Hopf subalgebra of $\NaturalHopfAlgebra \App \SetTerms \App
\Signature$. This result is analogous to~\cite[Theorem 3.13]{BG16}, which is within the
context of pros~\cite{McL65} rather than operads.

Let us consider the following example. Recall from
Section~\ref{subsubsec:faa_di_bruno_hopf_algebra} that the noncommutative Faà di Bruno Hopf
algebra $\NCFdB$ can be built as the natural Hopf algebra of the associative operad $\As$.
This operad is isomorphic to the quotient of $\SetTerms \App \Signature$ by the operad
congruence $\Equiv$ satisfying $\GenA \circ_1 \GenA \Equiv \GenA \circ_2 \GenA$, where
$\Signature$ is the binary signature $\Bra{\GenA}$. The Hopf algebra $\NCFdB$ can be
realized as a Hopf subalgebra of $\NaturalHopfAlgebra \App \SetTerms \App \Signature$
through the injection $\phi$ satisfying, for any $n \geq 1$,
\begin{equation}
    \phi \App \BasisE_{\alpha_n}
    =
    \sum_{\TermT \in \SetTerms \App \Signature} \Iverson{\Deg \App \TermT = n}
    \; \BasisE_\TermT.
\end{equation}
This is due to the fact that the canonical projection map $\pi_{\Equiv}$ satisfies
$\pi_{\Equiv} \App \TermT = \alpha_{\Deg \App \TermT}$ for any $\TermT \in \SetTerms
\App \Signature$. For instance,
\begin{equation}
    \phi \App \BasisE_{\alpha_3}
    =
    \BasisE_{
        \scalebox{0.75}{
            \begin{tikzpicture}[Centering,scale=0.4]
                \node[Leaf](1)at(1.5,0){};
                \node[Node](2)at(1.5,-1){$\GenA$};
                \node[Node](3)at(1,-2){$\GenA$};
                \node[Node](4)at(0.5,-3){$\GenA$};
                \node[Leaf](5)at(0,-4){};
                \node[Leaf](6)at(1,-4){};
                \node[Leaf](7)at(2,-3){};
                \node[Leaf](8)at(2.5,-2){};
                \draw[Edge](1)--(2);
                \draw[Edge](2)--(3);
                \draw[Edge](2)--(8);
                \draw[Edge](3)--(4);
                \draw[Edge](3)--(7);
                \draw[Edge](4)--(5);
                \draw[Edge](4)--(6);
            \end{tikzpicture}
        }
    }
    +
    \BasisE_{
        \scalebox{0.75}{
            \begin{tikzpicture}[Centering,scale=0.4]
                \node[Leaf](1)at(1.5,0){};
                \node[Node](2)at(1.5,-1){$\GenA$};
                \node[Node](3)at(1,-2){$\GenA$};
                \node[Leaf](4)at(0,-3){};
                \node[Node](5)at(1.5,-3){$\GenA$};
                \node[Leaf](6)at(1,-4){};
                \node[Leaf](7)at(2,-4){};
                \node[Leaf](8)at(2.5,-2){};
                \draw[Edge](1)--(2);
                \draw[Edge](2)--(3);
                \draw[Edge](2)--(8);
                \draw[Edge](3)--(4);
                \draw[Edge](3)--(5);
                \draw[Edge](5)--(6);
                \draw[Edge](5)--(7);
            \end{tikzpicture}
        }
    }
    +
    \BasisE_{
        \scalebox{0.75}{
            \begin{tikzpicture}[Centering,xscale=0.35,yscale=0.4]
                \node[Leaf](1)at(1.5,0){};
                \node[Node](2)at(1.5,-1){$\GenA$};
                \node[Node](3)at(0.5,-2){$\GenA$};
                \node[Leaf](4)at(0,-3){};
                \node[Leaf](5)at(1,-3){};
                \node[Node](6)at(2.5,-2){$\GenA$};
                \node[Leaf](7)at(2,-3){};
                \node[Leaf](8)at(3,-3){};
                \draw[Edge](1)--(2);
                \draw[Edge](2)--(3);
                \draw[Edge](2)--(6);
                \draw[Edge](3)--(4);
                \draw[Edge](3)--(5);
                \draw[Edge](6)--(7);
                \draw[Edge](6)--(8);
            \end{tikzpicture}
        }
    }
    +
    \BasisE_{
        \scalebox{0.75}{
            \begin{tikzpicture}[Centering,scale=0.4]
                \node[Leaf](1)at(1.5,0){};
                \node[Node](2)at(1.5,-1){$\GenA$};
                \node[Leaf](3)at(0.5,-2){};
                \node[Node](4)at(2,-2){$\GenA$};
                \node[Node](5)at(1.5,-3){$\GenA$};
                \node[Leaf](6)at(1,-4){};
                \node[Leaf](7)at(2,-4){};
                \node[Leaf](8)at(3,-3){};
                \draw[Edge](1)--(2);
                \draw[Edge](2)--(3);
                \draw[Edge](2)--(4);
                \draw[Edge](4)--(5);
                \draw[Edge](4)--(8);
                \draw[Edge](5)--(6);
                \draw[Edge](5)--(7);
            \end{tikzpicture}
        }
    }
    +
    \BasisE_{
        \scalebox{0.75}{
            \begin{tikzpicture}[Centering,scale=0.4]
                \node[Leaf](1)at(1.5,0){};
                \node[Node](2)at(1.5,-1){$\GenA$};
                \node[Leaf](3)at(0.5,-2){};
                \node[Node](4)at(2,-2){$\GenA$};
                \node[Leaf](5)at(1,-3){};
                \node[Node](6)at(2.5,-3){$\GenA$};
                \node[Leaf](7)at(2,-4){};
                \node[Leaf](8)at(3,-4){};
                \draw[Edge](1)--(2);
                \draw[Edge](2)--(3);
                \draw[Edge](2)--(4);
                \draw[Edge](4)--(5);
                \draw[Edge](4)--(6);
                \draw[Edge](6)--(7);
                \draw[Edge](6)--(8);
            \end{tikzpicture}
        }
    }.
\end{equation}

\section{Forest-like alphabets and polynomial realization}
\label{sec:polynomial_realization}
This section is the central part of this work. We first introduce the class of forest-like
alphabets. This particular class of related alphabet is required to build a polynomial
realization of $\NaturalHopfAlgebra \App \SetTerms \App \Signature$ where $\Signature$ is
any signature. We also present two particular forest-like alphabets that are useful for
establishing connections between $\NaturalHopfAlgebra \App \SetTerms \App \Signature$ and
some known Hopf algebras in the subsequent section.

\subsection{Forest-like alphabets and realizing map}
We begin here by defining the realizing map $\Realization_A$ of our polynomial realization
of $\NaturalHopfAlgebra \App \SetTerms \App \Signature$. We then prove some initial
properties of this map.

\subsubsection{Forest-like alphabets}
Let $S$ be a set and let $\RelatedAlphabetSignature \App S$ be the related alphabet
signature containing
\begin{enumerate}[label=({\sf \roman*})]
    \item a unary element $\RootRelation$;
    \item for any $s \in S$, a unary element $\DecorationRelation_s$;
    \item for any $j \geq 1$, a binary element $\EdgeRelation{j}$.
\end{enumerate}
An \Def{$S$-forest-like alphabet} is an $\RelatedAlphabetSignature \App S$-related alphabet
$A$. We call $\RootRelation^A$ the \Def{root relation} of $A$, $\DecorationRelation_s^A$, $s
\in S$, the \Def{$s$-decoration relation} of $A$, and $\EdgeRelation{j}^A$, $j \geq 1$, the
\Def{$j$-edge relation} of $A$. No conditions are required on these relations. We denote by
$\ClassForestLikeAlphabets \App S$ the class of $S$-forest-like alphabets. In the sequel, we
shall mainly consider $\Signature$-forest-like alphabets where $\Signature$ is a signature.

\subsubsection{Realizing map}
Let $\Signature$ be a signature and $A$ be an $\Signature$-forest-like alphabet. Given a
reduced $\Signature$-forest $\ForestF$, a word $w \in A^*$ is \Def{$A$-compatible} with
$\ForestF$ if the following four assertions are satisfied:
\begin{enumerate}[label=({\sf C\arabic*})]
    \item \label{item:compatibility_1}
    $\Length \App w = \Deg \App \ForestF$;
    \item \label{item:compatibility_2}
    for any internal node $i$ of $\ForestF$, if $i$ is a root of $\ForestF$, then $w \App i
    \in \RootRelation^A$;
    \item \label{item:compatibility_3}
    for any internal node $i$ of $\ForestF$, $w \App i \in
    \DecorationRelation_{\Decoration_\ForestF \App i}^A$;
    \item \label{item:compatibility_4}
    for any internal nodes $i$ and $i'$ of $\ForestF$, if $i \Edge{\ForestF}{j} i'$ for a $j
    \geq 1$, then $w \App i \EdgeRelation{j}^A w \App i'$.
\end{enumerate}
This property is denoted by $w \CompatibleWord{A} \ForestF$.

For instance, by setting
\begin{equation}
    \ForestF
    :=
    \scalebox{0.75}{
        \begin{tikzpicture}[Centering,scale=0.45]
            \node[Leaf](2)at(2,-1){};
            \node[Node](3)at(2,-2){$\GenB$};
            \node[Node](4)at(1,-3){$\GenC$};
            \node[Leaf](5)at(0,-4){};
            \node[Leaf](6)at(1,-4){};
            \node[Node](7)at(2,-4){$\GenA$};
            \node[Leaf](8)at(2,-5){};
            \node[Node](9)at(3.5,-3){$\GenA$};
            \node[Node](10)at(3.5,-4){$\GenB$};
            \node[Leaf](11)at(3,-5){};
            \node[Leaf](12)at(4,-5){};
            \node[Leaf](13)at(6,-1){};
            \node[Node](14)at(6,-2){$\GenC$};
            \node[Leaf](15)at(5,-3){};
            \node[Node](16)at(6,-3){$\GenA$};
            \node[Leaf](17)at(6,-4){};
            \node[Node](18)at(7.5,-3){$\GenB$};
            \node[Leaf](19)at(7,-4){};
            \node[Leaf](20)at(8,-4){};
            \draw[Edge](2)--(3);
            \draw[Edge](3)--(4);
            \draw[Edge](3)--(9);
            \draw[Edge](4)--(5);
            \draw[Edge](4)--(6);
            \draw[Edge](4)--(7);
            \draw[Edge](7)--(8);
            \draw[Edge](9)--(10);
            \draw[Edge](10)--(11);
            \draw[Edge](10)--(12);
            \draw[Edge](13)--(14);
            \draw[Edge](14)--(15);
            \draw[Edge](14)--(16);
            \draw[Edge](14)--(18);
            \draw[Edge](16)--(17);
            \draw[Edge](18)--(19);
            \draw[Edge](18)--(20);
        \end{tikzpicture}
    },
\end{equation}
any $A$-compatible word $w \in A^*$ with $\ForestF$ satisfies $\Length \App w = 8$, $w \App
1, w \App 6 \in \RootRelation^A$, $w \App 3, w \App 4, w \App 7 \in
\DecorationRelation_\GenA^A$, $w \App 1, w \App 5, w \App 8 \in
\DecorationRelation_\GenB^A$, $w \App 2, w \App 6 \in \DecorationRelation_\GenC^A$, $w \App
1 \EdgeRelation{1}^A w \App 2$, $w \App 1 \EdgeRelation{2}^A w \App 4$, $w \App 2
\EdgeRelation{3}^A w \App 3$, $w \App 4 \EdgeRelation{1}^A w \App 5$, $w \App 6
\EdgeRelation{2}^A w \App 7$, and $w \App 6 \EdgeRelation{3}^A w \App 8$.

Let $\Realization_A : \NaturalHopfAlgebra \App \SetTerms \App \Signature \to \K \Angle{A}$
be the linear map defined for any reduced $\Signature$-forest $\ForestF$ by
\begin{equation} \label{equ:realization_forest}
    \Realization_A \App \BasisE_{\ForestF}
    := \sum_{w \in A^*} \Iverson{w \CompatibleWord{A} \ForestF} \; w.
\end{equation}
The $A$-polynomial $\Realization_A \App \BasisE_\ForestF$ is the \Def{$A$-realization} of
$\ForestF$ on the $\BasisE$-basis.

\begin{Statement}{Proposition}{prop:associative_algebra_morphism}
    For any signature $\Signature$ and any $\Signature$-forest-like alphabet $A$, the map
    $\Realization_A$ is a graded unital associative algebra morphism.
\end{Statement}
\begin{Proof}
    Let $\ForestF_1, \ForestF_2 \in \Reduced \App \SetForests \App \Signature$. Observe that
    all nodes of $\ForestF_1$ are visited before all nodes of $\ForestF_2$ according to the
    left to right preorder traversal of the $\Signature$-forest $\ForestF_1 \Conc
    \ForestF_2$. Therefore, from the definition of the $A$-compatibility, it follows that
    for any $w \in A^*$, $w \CompatibleWord{A} \ForestF_1 \Conc \ForestF_2$ if and only if
    by setting $w_1$ as the prefix of $w$ of length $\Deg \App \ForestF_1$ and $w_2$ as the
    suffix of $w$ of length $\Deg \App \ForestF_2$, $w_1 \CompatibleWord{A} \ForestF_1$ and
    $w_2 \CompatibleWord{A} \ForestF_2$. This shows that $\Realization_A$ is an associative
    algebra morphism. Since moreover $\Realization_A \App \BasisE_\epsilon = 1$,
    $\Realization_A$ is a unital associative algebra morphism. Finally,
    Condition~\ref{item:compatibility_1} of the definition of compatibility implies that
    this morphism is graded.
\end{Proof}

\subsubsection{Realizing maps on quotients of forest-like alphabets}
\label{subsubsec:realizing_maps_quotient_alphabets}
Let $\Signature$ be a signature, $A$ be an $\Signature$-forest-like alphabet, and $\Equiv$
be a related alphabet congruence of $A$. Let us denote by
\begin{math}
    \pi_{\Equiv} : A \to A /_{\Equiv}
\end{math}
the canonical projection map associated with $\Equiv$. The map $\pi_{\Equiv}$ is extended as
a linear map from $\K \Angle{A}$ to $\K \Angle{A /_{\Equiv}}$ by setting
\begin{math}
    \pi_{\Equiv} \App a_1 \dots a_n
    := \Par{\pi_{\Equiv} \App a_1} \dots \Par{\pi_{\Equiv} \App a_n}
\end{math}
for any word $a_1 \dots a_n$, $n \geq 0$, on $A$.

Let us state a result useful to compute the $A /_{\Equiv}$-realization of a reduced
$\Signature$-forest $\ForestF$ on the $\BasisE$-basis from its $A$-realization, where
$\Equiv$ is a related alphabet congruence of $A$.

\begin{Statement}{Proposition}{prop:factorization_realization_congruence}
    Let $\Signature$ be a signature, $A$ be an $\Signature$-forest-like alphabet, and
    $\Equiv$ be a related alphabet congruence of $A$. If for any reduced $\Signature$-forest
    $\ForestF$, the map $\pi_{\Equiv}$ is a bijection when restricted on the domain of the
    words on $A$ which are $A$-compatible with $\ForestF$ and on the codomain of the words
    of $A /_{\Equiv}$ which are $A /_{\Equiv}$-compatible with $\ForestF$, then
    \begin{math}
        \Realization_{A /_{\Equiv}} = \pi_{\Equiv} \circ \Realization_A.
    \end{math}
\end{Statement}
\begin{Proof}
    Let $\ForestF \in \Reduced \App \SetForests \App \Signature$ and let us assume that the
    map $\pi_{\Equiv}$ is a bijection when restricted on the domain $X := \Bra{w \in A^* : w
    \CompatibleWord{A} \ForestF}$ and on the codomain $X' := \Bra{w' \in {A /_{\Equiv}}^* :
    w' \CompatibleWord{A /_{\Equiv}} \ForestF}$. We have first
    \begin{align} \label{equ:factorization_realization_congruence_1}
        \pi_{\Equiv} \App \Realization_A \App \BasisE_\ForestF
        & =
        \sum_{w \in A^*}
        \Iverson{w \CompatibleWord{A} \ForestF} \; \pi_{\Equiv} \App w
        \\
        & =
        \sum_{w \in A^*} \Iverson{w \CompatibleWord{A} \ForestF} \;
        \sum_{w' \in {A /_{\Equiv}}^*} \Iverson{\pi_{\Equiv} \App w = w'} \; w'
        \notag
        \\
        & =
        \sum_{w' \in {A /_{\Equiv}}^*} \sum_{w \in A^*}
        \Iverson{w \CompatibleWord{A} \ForestF} \Iverson{\pi_{\Equiv} \App w = w'} \; w'.
        \notag
    \end{align}
    Now, since $\Equiv$ is a related alphabet congruence of $A$ and $\pi_{\Equiv} : X \to
    X'$ is a bijection, the last term of~\eqref{equ:factorization_realization_congruence_1}
    is equal to
    \begin{align} \label{equ:factorization_realization_congruence_2}
        \sum_{w' \in {A /_{\Equiv}}^*} \sum_{w \in A^*}
        \Iverson{w \CompatibleWord{A} \ForestF} \Iverson{\pi_{\Equiv} \App w = w'}
        &
        \Iverson{w' \CompatibleWord{A /_{\Equiv}} \ForestF} \; w'
        \\
        & =
        \sum_{w' \in {A /_{\Equiv}}^*}
        \Iverson{w' \CompatibleWord{A /_{\Equiv}} \ForestF} \;
        \Par{
            \sum_{w \in A^*}
            \Iverson{w \CompatibleWord{A} \ForestF} \Iverson{\pi_{\Equiv} \App w = w'}
        }
        w'
        \notag
        \\
        & =
        \sum_{w' \in {A /_{\Equiv}}^*}
        \Iverson{w' \CompatibleWord{A /_{\Equiv}} \ForestF} \;
        \# \Bra{w \in X : \pi_{\Equiv} \App w = w'} \; w'
        \notag
        \\
        & =
        \sum_{w' \in {A /_{\Equiv}}^*}
        \Iverson{w' \CompatibleWord{A /_{\Equiv}} \ForestF} \; w'.
        \notag
    \end{align}
    The last term of~\eqref{equ:factorization_realization_congruence_2} is equal
    to~$\Realization_{A /_{\Equiv}} \App \BasisE_\ForestF$, which shows the stated property.
\end{Proof}

It is possible to confer greater generality to the statement of
Proposition~\ref{prop:factorization_realization_congruence} so that it essentially works
for any polynomial realization whose realizing map admits an expression analogous
to~\eqref{equ:realization_forest}. We did not, however, write it in these terms since
this degree of generality is not necessary for this work.

\subsection{Compatibility with the coproduct}
We define now a disjoint sum operation on the class $\ClassForestLikeAlphabets \App
\Signature$ of $\Signature$-forest-like alphabets in order to prove that $\Realization_A$ is
compatible with the coproduct of $\NaturalHopfAlgebra \App \SetTerms \App \Signature$.

\subsubsection{Disjoint sum of forest-like alphabets}
Let $\Signature$ be a signature. The \Def{disjoint sum} of two $\Signature$-forest-like
alphabets $A_1$ and $A_2$ is the $\Signature$-forest-like alphabet $A_1 \AlphabetSum A_2$
defined as the set $A_1 \sqcup A_2$ and such that
\begin{enumerate}[label=({\sf \roman*})]
    \item $\RootRelation^{A_1 \AlphabetSum A_2} := \RootRelation^{A_1} \sqcup
    \RootRelation^{A_2}$;
    \item for any $\GenS \in \Signature$, $\DecorationRelation_\GenS^{A_1 \AlphabetSum A_2}
    := \DecorationRelation_\GenS^{A_1} \sqcup \DecorationRelation_\GenS^{A_2}$;
    \item for any $a, a' \in A$, $a \EdgeRelation{j}^{A_1 \AlphabetSum A_2} a'$ holds if
    \begin{enumerate}[label=({\sf \alph*})]
        \item $a, a' \in A_1$ and $a \EdgeRelation{j}^{A_1} a'$,
        \item or $a, a' \in A_2$ and $a \EdgeRelation{j}^{A_2} a'$,
        \item or $a \in A_1$, $a' \in A_2$, and $a \in \RootRelation^{A_2}$.
    \end{enumerate}
\end{enumerate}
This operation $\AlphabetSum$ on $\ClassForestLikeAlphabets$ is clearly associative and
admits the empty $\Signature$-forest-like alphabet $\emptyset$ as unit.

\subsubsection{Compatibility with the coproduct}
We shall prove here that the realizing map $\Realization_A$ is compatible with the coproduct
of $\NaturalHopfAlgebra \App \SetTerms \App \Signature$. To establish this property, we need
the following lemma.

\begin{Statement}{Lemma}{lem:compatibility_sum_alphabets}
    Let $\Signature$ be a signature, $A_1$ and $A_2$ be two $\Signature$-forest-like
    alphabets, $\ForestF$ be a reduced $\Signature$-forest, and $w$ be a word on $A_1
    \AlphabetSum A_2$. By setting $I_1 := \Occurrences{A_1} \App w$ and $I_2 :=
    \Occurrences{A_2} \App w$, the following two assertions are equivalent:
    \begin{enumerate}[label=({\sf \roman*})]
        \item \label{item:compatibility_sum_alphabets_1}
        the word $w$ is $A_1 \AlphabetSum A_2$-compatible with $\ForestF$;
        \item \label{item:compatibility_sum_alphabets_2}
        the pair $\Par{I_1, I_2}$ is $\ForestF$-admissible, the word $w_{|I_1}$ is
        $A_1$-compatible with $\ForestF \App I_1$, and the word $w_{|I_2}$ is
        $A_2$-compatible with $\ForestF \App I_2$.
    \end{enumerate}
\end{Statement}
\begin{Proof}
    For any $k \in [2]$, let $\pi^{(k)} : I_k \to \Han{\# I_k}$ be the map such that
    $\pi^{(k)} \App i$ is the position in $w_{|I_k}$ of the letter of position $i$ of $w$.

    Assume first that~\ref{item:compatibility_sum_alphabets_1} holds. Let $i' \in I_1$ such
    that there is an internal node $i$ of $\ForestF$ satisfying $i \Edge{\ForestF}{j} i'$
    for a $j \geq 1$. Since $w \CompatibleWord{A_1 \AlphabetSum A_2} \ForestF$, we have $w
    \App i \EdgeRelation{j}^{A_1 \AlphabetSum A_2} w \App i'$. Since $w \App i' \in A_1$, by
    definition of the operation $\AlphabetSum$, we necessarily have $w \App i
    \EdgeRelation{j}^{A_1} w \App i'$ with $w \App i \in A_1$. This shows that $i \in I_1$
    and that $I_1$ is closed w.r.t.\ the ancestor relation of $\ForestF$. Therefore, $I_2$
    is closed w.r.t.\ the descendant relation of $\ForestF$, which shows that $\Par{I_1,
    I_2} \AdmissibleSet \ForestF$. Now, by definition of the disjoint sum operation
    $\AlphabetSum$ on $\Signature$-forest-like alphabets, we have the following properties.
    \begin{enumerate}[label=({\sf \arabic*})]
        \item For any $k \in [2]$, since $\Length \App w_{|I_k} = \# I_k$ and $\Deg \App
        \ForestF \App I_k = \# I_k$, we have $\Length \App w_{|I_k} = \Deg \App \ForestF
        \App I_k$.
        \item Let $i \in I_k$, $k \in [2]$, such that $i$ is a root of $\ForestF \App I_k$.
        \begin{enumerate}[label=({\sf \alph*})]
            \item Assume that $k = 1$. Since $\Par{I_1, I_2} \AdmissibleSet \ForestF$, $i$
            is also a root of $\ForestF$. Therefore, since $w \App i \in \RootRelation^{A_1
            \AlphabetSum A_2}$, we have $w_{|I_1} \App \pi^{(1)} \App i \in
            \RootRelation^{A_1}$.
            \item Assume that $k = 2$. Since $\Par{I_1, I_2} \AdmissibleSet \ForestF$, we
            have two possibilities: either $i$ is a root of $\ForestF$, or $i$ is a child of
            an internal node $i'$ of $\ForestF$ such that $i' \in I_1$. In the first case,
            we have $w \App i \in \RootRelation^{A_1 \AlphabetSum A_2}$ and thus, $w_{|I_2}
            \App \pi^{(2)} \App i \in \RootRelation^{A_2}$. In the second case, we have $i'
            \Edge{\ForestF}{j} i$ for a $j \geq 1$, so that $w \App i' \EdgeRelation{j}^{A_1
            \AlphabetSum A_2} w \App i$. Since $w \App i' \in A_1$ and $w \App i \in A_2$,
            we have $w \App i \in \RootRelation^{A_2}$ and thus, $w_{|I_2} \App \pi^{(2)}
            \App i \in \RootRelation^{A_2}$.
        \end{enumerate}
        \item For any $k \in [2]$, let $i \in I_k$ such that the internal node $i$ of
        $\ForestF$ is decorated by $\GenS \in \Signature$. Since $w \App i \in
        \DecorationRelation_\GenS^{A_1 \AlphabetSum A_2}$, we have $w_{|I_k} \App \pi^{(k)}
        \App i \in \DecorationRelation_\GenS^{A_k}$.
        \item For any $k \in [2]$, let $i, i' \in I_k$ such that $i \Edge{\ForestF}{j} i'$
        for a $j \geq 1$. Since $w \App i \EdgeRelation{j}^{A_1 \AlphabetSum A_2} w \App
        i'$, we have $w_{|I_k} \App \pi^{(k)} \App i \EdgeRelation{j}^{A_k} w_{|I_k} \App
        \pi^{(k)} \App i'$.
    \end{enumerate}
    These properties together imply that~\ref{item:compatibility_sum_alphabets_2} holds.

    Assume conversely that~\ref{item:compatibility_sum_alphabets_2} holds. Again by
    definition of the disjoint sum operation $\AlphabetSum$ on $\Signature$-forest-like
    alphabets, we have the following properties.
    \begin{enumerate}[label=({\sf \arabic*})]
        \item Let $i$ be a root of $\ForestF$. For any $k \in [2]$, if $i \in I_k$, since
        $w_{|I_k} \CompatibleWord{A_k} \ForestF \App I_k$, we have $w_{|I_k} \App \pi^{(k)}
        \App i \in \RootRelation^{A_k}$, so that $w \App i \in \RootRelation^{A_k}$ and $w
        \App i \in \RootRelation^{A_1 \AlphabetSum A_2}$.
        \item For any $k \in [2]$, let $i$ be an internal node of $\ForestF$ decorated by
        $\GenS \in \Signature$. Since $w_{|I_k} \CompatibleWord{A_k} \ForestF \App I_k$, we
        have $w_{|I_k} \App \pi^{(k)} \App i \in \DecorationRelation_\GenS^{A_k}$, so that
        $w \App i \in \DecorationRelation_\GenS^{A_k}$ and $w \App i \in
        \DecorationRelation_\GenS^{A_1 \AlphabetSum A_2}$.
        \item Let $i$ and $i'$ be two internal nodes of $\ForestF$ such that $i
        \Edge{\ForestF}{j} i'$ for a $j \geq 1$. We have four possibilities, factored into
        three, depending on the set $I_1$ or $I_2$ to which each $i$ and $i'$ belong.
        \begin{enumerate}[label=({\sf \alph*})]
            \item If $i, i' \in I_k$ for a $k \in [2]$, since $w_{|I_k} \CompatibleWord{A_k}
            \ForestF \App I_k$, we have $w_{|I_k} \App \pi^{(k)} \App i
            \EdgeRelation{j}^{A_k} w_{|I_k} \App \pi^{(k)} \App i'$, so that $w \App i
            \EdgeRelation{j}^{A_k} w \App i'$ and $w \App i \EdgeRelation{j}^{A_1
            \AlphabetSum A_2} w \App i'$.
            \item If $i \in I_1$ and $i' \in I_2$, then, since $\Par{I_1, I_2}$ is
            $\ForestF$-admissible, $\pi^{(k)} \App i'$ is a root of $\ForestF \App I_2$.
            Therefore, since $w_{|I_2} \CompatibleWord{A_2} \ForestF \App I_2$, we have
            $w_{|I_2} \App \pi^{(k)} \App i' \in \RootRelation^{A_2}$ and $w \App i' \in
            \RootRelation^{A_2}$. Since moreover $w \App i \in A_1$, we have $w \App i
            \EdgeRelation{j}^{A_1 \AlphabetSum A_2} w \App i'$.
            \item If $i \in I_2$ and $i' \in I_1$, the node $i'$ of $\ForestF$ is such that
            its parent is not in $I_1$. This is contradictory to the fact that $\Par{I_1,
            I_2}$ is $\ForestF$-admissible so that this case cannot occur.
        \end{enumerate}
    \end{enumerate}
    These properties together imply that~\ref{item:compatibility_sum_alphabets_1} holds.
\end{Proof}

\begin{Statement}{Proposition}{prop:realization_alphabet_doubling}
    For any signature $\Signature$, any $\Signature$-forest-like alphabets $A_1$ and $A_2$,
    and any reduced $\Signature$-forest $\ForestF$,
    \begin{equation}
        \WordToTensor_{A_1, A_2}
        \App \Realization_{A_1 \AlphabetSum A_2}
        \App \BasisE_\ForestF
        =
        \Par{\Realization_{A_1} \otimes \Realization_{A_2}}
        \App \Coproduct \App \BasisE_\ForestF.
    \end{equation}
\end{Statement}
\begin{Proof}
    First, by definition of the map $\Realization_{A_1 \AlphabetSum A_2}$, we have
    \begin{equation} \label{equ:realization_alphabet_doubling_1}
        \WordToTensor_{A_1, A_2}
        \App \Realization_{A_1 \AlphabetSum A_2} \App \BasisE_\ForestF
        =
        \sum_{w \in \Par{A_1 \AlphabetSum A_2}^*}
        \Iverson{w \CompatibleWord{A_1 \AlphabetSum A_2} \ForestF}
        \; \WordToTensor_{A_1, A_2} \App w.
    \end{equation}
    By decomposing the sum intervening in the right-hand side
    of~\eqref{equ:realization_alphabet_doubling_1} following the occurrences of the letters
    of $w$ which belong to $A_1$ and to $A_2$, this polynomial is equal to
    \begin{equation} \label{equ:realization_alphabet_doubling_2}
        \sum_{I_1, I_2 \subseteq [\Deg \App \ForestF]}
        \enspace
        \sum_{w \in \Par{A_1 \AlphabetSum A_2}^*}
        \Iverson{\Occurrences{A_1} \App w = I_1}
        \Iverson{\Occurrences{A_2} \App w = I_2}
        \Iverson{w \CompatibleWord{A_1 \AlphabetSum A_2} \ForestF}
        \;
        \WordToTensor_{A_1, A_2} \App w.
    \end{equation}
    Now, by Lemma~\ref{lem:compatibility_sum_alphabets},
    \eqref{equ:realization_alphabet_doubling_2} is equal to
    \begin{multline} \label{equ:realization_alphabet_doubling_3}
        \sum_{I_1, I_2 \subseteq [\Deg \App \ForestF]}
        \Iverson{\Par{I_1, I_2} \AdmissibleSet \ForestF}
        \enspace
        \\
        \sum_{w \in \Par{A_1 \AlphabetSum A_2}^*}
        \Iverson{\Occurrences{A_1} \App w = I_1}
        \Iverson{\Occurrences{A_2} \App w = I_2}
        \Iverson{w_{|A_1} \CompatibleWord{A_1} \ForestF \App I_1}
        \Iverson{w_{|A_2} \CompatibleWord{A_2} \ForestF \App I_2}
        \;
        \WordToTensor_{A_1, A_2} \App w.
    \end{multline}
    By expressing the second sum of~\eqref{equ:realization_alphabet_doubling_3} by summing
    instead on $w_1$ and $w_2$ which are respectively the subwords $w_{|A_1}$ and $w_{|A_2}$
    of $w$ where, additionally, $I_1$ (resp.\ $I_2$) virtually specifies the positions of
    the letters of $w_1$ (resp.\ $w_2$) in $w$, this polynomial is equal to
    \begin{equation} \label{equ:realization_alphabet_doubling_4}
        \sum_{I_1, I_2 \subseteq [\Deg \App \ForestF]}
        \Iverson{\Par{I_1, I_2} \AdmissibleSet \ForestF}
        \enspace
        \sum_{\substack{
            w_1 \in A_1^* \\
            w_2 \in A_2^*
        }}
        \Iverson{w_1 \CompatibleWord{A_1} \ForestF \App I_1}
        \Iverson{w_2 \CompatibleWord{A_2} \ForestF \App I_2}
        \;
        w_1 \otimes w_2.
    \end{equation}
    Finally, \eqref{equ:realization_alphabet_doubling_4} rewrites as
    \begin{equation}
        \sum_{I_1, I_2 \subseteq [\Deg \App \ForestF]}
        \Iverson{\Par{I_1, I_2} \AdmissibleSet \ForestF}
        \enspace
        \Par{\sum_{w_1 \in A_1^*}
        \Iverson{w_1 \CompatibleWord{A_1} \ForestF \App I_1} \; w_1}
        \otimes
        \Par{\sum_{w_2 \in A_2^*}
        \Iverson{w_2 \CompatibleWord{A_2} \ForestF \App I_2} \; w_2}
    \end{equation}
    and, by Proposition~\ref{prop:admissible_sets_coproduct}, as
    \begin{equation}
        \sum_{I_1, I_2 \subseteq [\Deg \App \ForestF]}
        \Iverson{\Par{I_1, I_2} \AdmissibleSet \ForestF}
        \;
        \Realization_{A_1} \App \BasisE_{\ForestF \App I_1}
        \otimes
        \Realization_{A_2} \App \BasisE_{\ForestF \App I_2}
        =
        \Par{\Realization_{A_1} \otimes \Realization_{A_2}}
        \App \Coproduct \App \BasisE_\ForestF.
    \end{equation}
    This shows the expected property.
\end{Proof}

\subsection{Forest-like alphabet of positions and injectivity}
We introduce now a particular $\Signature$-forest-like alphabet $\PositionsAlphabet \App
\Signature$ for which the realizing map $\Realization_{\PositionsAlphabet \App \Signature}$
is injective. We also introduce and study a related alphabet quotient $\LengthsAlphabet \App
\Signature$ of this $\Signature$-forest-like alphabet.

\subsubsection{Forest-like alphabet of positions}
For any signature $\Signature$, the \Def{$\Signature$-forest-like alphabet of positions} is
the $\Signature$-forest-like alphabet
\begin{equation}
    \PositionsAlphabet \App \Signature
    :=
    \Bra{\Letter^\GenS_u : \GenS \in \Signature \text{ and } u \in \N^*}
\end{equation}
such that
\begin{enumerate}[label=({\sf \roman*})]
    \item the root relation $\RootRelation^{\PositionsAlphabet \App \Signature}$
    satisfies
    \begin{math}
        \Letter^\GenS_{0^\ell}
        \in \RootRelation^{\PositionsAlphabet \App \Signature}
    \end{math}
    for any $\ell \in \N$;
    \item for any $\GenS \in \Signature$, the $\GenS$-decoration relation
    $\DecorationRelation_\GenS^{\PositionsAlphabet \App \Signature}$
    satisfies
    \begin{math}
        \Letter^\GenS_u
        \in \DecorationRelation_\GenS^{\PositionsAlphabet \App \Signature}
    \end{math}
    for any $u \in \N^*$;
    \item for any $j \geq 1$, the $j$-edge relation $\EdgeRelation{j}^{\PositionsAlphabet
    \App \Signature}$ satisfies
    \begin{math}
        \Letter^\GenS_u
        \EdgeRelation{j}^{\PositionsAlphabet \App \Signature}
        \Letter^{\GenS'}_{u \Conc j \Conc 0^\ell}
    \end{math}
    for any $\GenS, \GenS' \in \Signature$, $u \in \N^*$, and $\ell \in \N$.
\end{enumerate}

Here is an example of the $\PositionsAlphabet \App \SignatureExample$-realization
of an $\SignatureExample$-forest on the $\BasisE$-basis:
\begin{equation}
    \Realization_{\PositionsAlphabet \App \Signature} \App
    \BasisE_{
        \scalebox{0.75}{
            \begin{tikzpicture}[Centering,xscale=0.35,yscale=0.45]
                \node[Leaf](2)at(2.5,-1){};
                \node[Node](3)at(2.5,-2){$\GenC$};
                \node[Node](4)at(1,-3){$\GenB$};
                \node[Leaf](5)at(0.5,-4){};
                \node[Leaf](6)at(1.5,-4){};
                \node[Leaf](7)at(2.5,-3){};
                \node[Node](8)at(4,-3){$\GenA$};
                \node[Node](9)at(4,-4){$\GenC$};
                \node[Leaf](10)at(3,-5){};
                \node[Leaf](11)at(4,-5){};
                \node[Leaf](12)at(5,-5){};
                \node[Leaf](13)at(6.5,-1){};
                \node[Node](14)at(6.5,-2){$\GenB$};
                \node[Node](15)at(6,-3){$\GenA$};
                \node[Leaf](16)at(6,-4){};
                \node[Leaf](17)at(7,-3){};
                \draw[Edge](2)--(3);
                \draw[Edge](3)--(4);
                \draw[Edge](3)--(7);
                \draw[Edge](3)--(8);
                \draw[Edge](4)--(5);
                \draw[Edge](4)--(6);
                \draw[Edge](8)--(9);
                \draw[Edge](9)--(10);
                \draw[Edge](9)--(11);
                \draw[Edge](9)--(12);
                \draw[Edge](13)--(14);
                \draw[Edge](14)--(15);
                \draw[Edge](14)--(17);
                \draw[Edge](15)--(16);
            \end{tikzpicture}
        }
    }
    =
    \sum_{\ell_1, \dots, \ell_6 \in \N}
    \Letter^\GenC_{0^{\ell_1}} \;
    \Letter^\GenB_{0^{\ell_1} 1 0^{\ell_2}} \;
    \Letter^\GenA_{0^{\ell_1} 3 0^{\ell_3}} \;
    \Letter^\GenC_{0^{\ell_1} 3 0^{\ell_3} 1 0^{\ell_4}} \;
    \Letter^\GenB_{0^{\ell_5}} \;
    \Letter^\GenA_{0^{\ell_5} 1^{\ell_6}}.
\end{equation}

\subsubsection{Injectivity}
Let $\PositionsEncoding : \SetForests \App \Signature \to \PositionsAlphabet \App
\Signature^*$ be the map defined for any $\Signature$-forest $\ForestF$ of degree $n
\geq 0$ by
\begin{equation}
    \PositionsEncoding \App \ForestF
    :=
    \Letter^{\Decoration_\ForestF \App 1}_{\Position_\ForestF \App 1}
    \dots
    \Letter^{\Decoration_\ForestF \App n}_{\Position_\ForestF \App n}.
\end{equation}
For instance,
\begin{equation} \label{equ:positions_encoding_example}
    \PositionsEncoding \App
    \scalebox{0.75}{
        \begin{tikzpicture}[Centering,scale=0.45]
            \node[Leaf](2)at(2.25,-1){};
            \node[Node](3)at(2.25,-2){$\GenB$};
            \node[Node](4)at(1,-3){$\GenC$};
            \node[Leaf](5)at(0,-4){};
            \node[Leaf](6)at(1,-4){};
            \node[Node](7)at(2,-4){$\GenA$};
            \node[Leaf](8)at(2,-5){};
            \node[Node](9)at(3.5,-3){$\GenA$};
            \node[Node](10)at(3.5,-4){$\GenB$};
            \node[Leaf](11)at(3,-5){};
            \node[Leaf](12)at(4,-5){};
            \node[Leaf](13)at(6,-1){};
            \node[Node](14)at(6,-2){$\GenC$};
            \node[Leaf](15)at(4.5,-3){};
            \node[Node](16)at(6,-3){$\GenA$};
            \node[Leaf](17)at(6,-4){};
            \node[Node](18)at(7.5,-3){$\GenB$};
            \node[Leaf](19)at(7,-4){};
            \node[Leaf](20)at(8,-4){};
            \draw[Edge](2)--(3);
            \draw[Edge](3)--(4);
            \draw[Edge](3)--(9);
            \draw[Edge](4)--(5);
            \draw[Edge](4)--(6);
            \draw[Edge](4)--(7);
            \draw[Edge](7)--(8);
            \draw[Edge](9)--(10);
            \draw[Edge](10)--(11);
            \draw[Edge](10)--(12);
            \draw[Edge](13)--(14);
            \draw[Edge](14)--(15);
            \draw[Edge](14)--(16);
            \draw[Edge](14)--(18);
            \draw[Edge](16)--(17);
            \draw[Edge](18)--(19);
            \draw[Edge](18)--(20);
        \end{tikzpicture}
    }
    =
    \Letter^\GenB_{\epsilon} \;
    \Letter^\GenC_{1} \;
    \Letter^\GenA_{13} \;
    \Letter^\GenA_{2} \;
    \Letter^\GenB_{21} \;
    \Letter^\GenC_{\epsilon} \;
    \Letter^\GenA_{2} \;
    \Letter^\GenB_{3}.
\end{equation}

\begin{Statement}{Lemma}{lem:positions_compatible}
    For any signature $\Signature$ and any reduced $\Signature$-forest $\ForestF$,
    the word $\PositionsEncoding \App \ForestF$ is $\PositionsAlphabet \App
    \Signature$-compatible with $\ForestF$.
\end{Statement}
\begin{Proof}
    Let us show that $w := \PositionsEncoding \App \ForestF$ satisfies the four conditions
    to be $\PositionsAlphabet \App \Signature$-compatible with $\ForestF$. First, since
    $\Length \App w = \Deg \App \ForestF$, Condition~\ref{item:compatibility_1} holds.
    Moreover, for any root $i$ of $\ForestF$, we have $w \App i =
    \Letter^{\Decoration_\ForestF \App i}_ \epsilon$. Therefore, $w \App i \in
    \RootRelation^{\PositionsAlphabet \App \Signature}$, showing
    that~\ref{item:compatibility_2} holds. Let $i$ be an internal node of $\ForestF$. Since
    $w \App i = \Letter^{\Decoration_\ForestF \App i}_{\Position_\ForestF \App i}$, we have
    that
    \begin{math}
        w \App i \in
        \DecorationRelation^{\PositionsAlphabet \App \Signature}_
            {\Decoration_\ForestF \App i}.
    \end{math}
    Hence, \ref{item:compatibility_3} checks out. Assume that $i$ and $i'$ are two internal
    nodes of $\ForestF$ such that $i \Edge{\ForestF}{j} i'$ for a $j \geq 1$. By definition
    of internal node positions, this implies that $\Position_\ForestF \App i' =
    \Position_\ForestF \App i \Conc j$. Now, since $w \App i = \Letter^{\Decoration_\ForestF
    \App i}_{\Position_\ForestF \App i}$ and $w \App i' = \Letter^{\Decoration_\ForestF \App
    i'}_{\Position_\ForestF \App i'}$, we have
    \begin{math}
        \Letter^{\Decoration_\ForestF \App i}_{\Position_\ForestF \App i}
        \EdgeRelation{j}
        \Letter^{\Decoration_\ForestF \App i'}_{\Position_\ForestF \App i'}.
    \end{math}
    Therefore, \ref{item:compatibility_4} holds, showing that $w
    \CompatibleWord{\PositionsAlphabet \App \Signature} \ForestF$.
\end{Proof}

Given a word $w := \Letter^{\GenS_1}_{u_1} \dots \Letter^{\GenS_n}_{u_n}$, $n \geq 0$, on
$\PositionsAlphabet \App \Signature$, the \Def{weight} $\Weight \App w$ of $w$ is $\Length
\App u_1 + \dots + \Length \App u_n$. For instance, the weight of the word appearing in the
right-hand side of~\eqref{equ:positions_encoding_example} is~$8$.

\begin{Statement}{Lemma}{lem:positions_triangular}
    For any signature $\Signature$ and any reduced $\Signature$-forest $\ForestF$,
    \begin{equation}
        \Realization_{\PositionsAlphabet \App \Signature} \App \BasisE_\ForestF
        = \PositionsEncoding \App \ForestF
        + \sum_{w \in \PositionsAlphabet \App \Signature^*}
        \Iverson{w \CompatibleWord{\PositionsAlphabet \App \Signature} \ForestF}
        \Iverson{\Weight \App w > \Weight \App \PositionsEncoding \App \ForestF} \;
        w.
    \end{equation}
\end{Statement}
\begin{Proof}
    Let $w := \Letter^{\GenS_1}_{u_1} \dots \Letter^{\GenS_n}_{u_n}$, $n \geq 0$, be a word
    on $\PositionsAlphabet \App \Signature$ such that $w \CompatibleWord{\PositionsAlphabet
    \App \Signature} \ForestF$. Since $w$ is $\PositionsAlphabet \App \Signature$-compatible
    with $\ForestF$, $n = \Deg \App \ForestF$ and for any internal node $i$ of $\ForestF$,
    $\GenS_i = \Decoration_\ForestF \App i$. Let us show that for any $i \in [n]$,
    ${u_i}_{|\N \setminus \{0\}} = \Position_\ForestF \App i$. Since $w$ is
    $\PositionsAlphabet \App \Signature$-compatible with $\ForestF$, if $i$ and $i'$ are two
    internal nodes of $\ForestF$ such that $i \Edge{\ForestF}{j} i'$ for a $j \geq 1$, then
    $u_{i'} = u_i \Conc j \Conc 0^\ell$ for an $\ell \in \N$. This implies that, for any
    internal node $i$ of $\ForestF$, by denoting by $i_0$ the root of the $\Signature$-term
    to which $i$ belongs and by $i_1, \dots, i_{k - 1}$ the internal nodes of $\ForestF$
    such that
    \begin{math}
        i_0 \Edge{\ForestF}{j_1} i_1 \Edge{\ForestF}{j_2}
        \dots
        \Edge{\ForestF}{j_{k - 1}} i_{k - 1} \Edge{\ForestF}{j_k} i
    \end{math}
    for some positive integers $j_1, j_2, \dots, j_{k - 1}, j_k$, we have
    \begin{math}
        u_i = 0^{\ell_0} j_1 0^{\ell_1} j_2 0^{\ell_2}
        \dots j_{k - 1} 0^{\ell_{k - 1}} j_k 0^{\ell_k}
    \end{math}
    for some $\ell_0, \ell_1, \ell_2, \dots, \ell_{k - 1}, \ell_k \in \N$. Therefore, since
    \begin{math}
        {u_i}_{|\N \setminus \{0\}}
        = j_1 j_2 \dots j_{k - 1} j_k
        = \Position_\ForestF \App i,
    \end{math}
    the expected property is established.

    This property together with Lemma~\ref{lem:positions_compatible} show that all monomials
    appearing in $\Realization_{\PositionsAlphabet \App \Signature} \App \BasisE_\ForestF$
    are of the form $w := \Letter^{\Decoration_\ForestF \App 1}_{u_1} \dots
    \Letter^{\Decoration_\ForestF \App n}_{u_n}$ where $n := \Deg \App \ForestF$ and
    ${u_i}_{|\N \setminus \{0\}} = \Position_\ForestF \App i$ for all $i \in [n]$. Hence,
    $\Weight \App \PositionsEncoding \App \ForestF \leq \Weight \App w$. The property of the
    statement follows.
\end{Proof}

We can now state the next property required to establish that $\Realization_A$ is a
polynomial realization of~$\NaturalHopfAlgebra \App \SetTerms \App \Signature$.

\begin{Statement}{Proposition}{prop:realization_injectivity}
    For any signature $\Signature$, the map $\Realization_{\PositionsAlphabet \App
    \Signature}$ is injective.
\end{Statement}
\begin{Proof}
    By Lemma~\ref{lem:positions_compatible}, $\Realization_{\PositionsAlphabet \App
    \Signature}$ sends any basis element $\BasisE_\ForestF$ of $\NaturalHopfAlgebra \App
    \SetTerms \App \Signature$ to an $\PositionsAlphabet \App \Signature$-polynomial
    in which the monomial $\PositionsEncoding \App \ForestF$ appears. This monomial encodes
    the decoration and the position of each internal node of $\ForestF$, from the first to
    the last one. It is thus possible to reconstruct $\ForestF$ from $\PositionsEncoding
    \App \ForestF$. Moreover, by Lemma~\ref{lem:positions_triangular}, it is possible to
    reconstruct $\ForestF$ from $\Realization_{\PositionsAlphabet \App \Signature} \App
    \BasisE_\ForestF$. This shows the stated property.
\end{Proof}

\subsubsection{Polynomial realizations of natural Hopf algebras of free operads}
We are now in position to state the main result of this work, namely a polynomial
realization of natural Hopf algebras of free operads.

\begin{Statement}{Theorem}{thm:realization}
    For any signature $\Signature$, the class of $\Signature$-forest-like alphabets,
    together with the alphabet disjoint sum operation $\AlphabetSum$, the map
    $\Realization_A$, and the alphabet $\PositionsAlphabet \App \Signature$, form a
    polynomial realization of the Hopf algebra $\NaturalHopfAlgebra \App \SetTerms \App
    \Signature$.
\end{Statement}
\begin{Proof}
    This is a consequence of Propositions~\ref{prop:associative_algebra_morphism},
    \ref{prop:realization_alphabet_doubling}, and~\ref{prop:realization_injectivity}.
\end{Proof}

\subsubsection{Polynomial realizations of natural Hopf algebras of operads}
Recall from Section~\ref{subsec:quotients_hopf_subalgebras} that when $\Operad$ is a
quotient of a free operad $\SetTerms \App \Signature$ where $\Signature$ is a signature, the
map $\phi$ defined by~\eqref{equ:injection_from_quotient_operad} is a Hopf algebra injection
from $\NaturalHopfAlgebra \App \Operad$ to $\NaturalHopfAlgebra \App \SetTerms \App
\Signature$. For any $\Signature$-forest-like alphabet $A$, let $\bar{\Realization}_A :
\NaturalHopfAlgebra \App \Operad \to \K \Angle{A}$ be the map $\Realization_A \circ \phi$.
The $A$-polynomial $\bar{\Realization}_A \App \BasisE_x$ is the \Def{$A$-realization} of $x
\in \Operad$ on the $\BasisE$-basis.

\begin{Statement}{Proposition}{prop:realization_quotient_operads}
    Let $\Signature$ be a signature and $\Operad$ be the quotient of the free operad
    $\SetTerms \App \Signature$ by an operad congruence $\Equiv$ which is compatible with
    the degree and of finite type. The class of $\Signature$-forest-like alphabets, together
    with the alphabet disjoint sum operation $\AlphabetSum$, the map $\bar{\Realization}_A$,
    and the alphabet $\PositionsAlphabet \App \Signature$, form a polynomial realization of
    the Hopf algebra~$\NaturalHopfAlgebra \App \Operad$.
\end{Statement}
\begin{Proof}
    This is a consequence of the fact that, by
    Theorem~\ref{thm:quotient_operad_hopf_morphism}, $\phi$ is an injective Hopf algebra
    morphism from $\NaturalHopfAlgebra \App \Operad$ to $\NaturalHopfAlgebra \App \SetTerms
    \App \Signature$ and of the fact that, by Theorem~\ref{thm:realization},
    $\Par{\ClassForestLikeAlphabets \App \Signature, \AlphabetSum, \Realization_A,
    \PositionsAlphabet \App \Signature}$ is a polynomial realization of $\NaturalHopfAlgebra
    \App \SetTerms \App \Signature$. From this, it follows straightforwardly that
    $\bar{\Realization}_A$ satisfies the required conditions to form a polynomial
    realization of~$\NaturalHopfAlgebra \App \Operad$.
\end{Proof}

Let us consider an example of application of
Proposition~\ref{prop:realization_quotient_operads} to build a polynomial realization of the
noncommutative Faà du Bruno Hopf algebra $\NCFdB$. The injection of $\NaturalHopfAlgebra
\App \As$, which is isomorphic to $\NCFdB$ (see
Section~\ref{subsubsec:faa_di_bruno_hopf_algebra}), is presented in
Section~\ref{subsubsec:hopf_subalgebras_free_operads}. By denoting by $\Signature$ the
binary signature $\Bra{\GenA}$, for any $\Signature$-forest-like alphabet $A$ and $n \geq
1$,
\begin{equation}
    \bar{\Realization}_A \App \BasisE_{\alpha_n}
    =
    \sum_{\TermT \in \SetTerms \App \Signature} \Iverson{\Deg \App \TermT = n}
    \; \Realization_A \App \BasisE_\TermT.
\end{equation}
In particular, we have
\begin{equation}
    \bar{\Realization}_{\PositionsAlphabet \App \Signature} \App \BasisE_{\alpha_1}
   =
   \Realization_{\PositionsAlphabet \App \Signature} \App
   \BasisE_{
        \scalebox{0.75}{
            \begin{tikzpicture}[Centering,scale=0.4]
                \node[Leaf](1)at(0.5,0){};
                \node[Node](2)at(0.5,-1){$\GenA$};
                \node[Leaf](3)at(0,-2){};
                \node[Leaf](4)at(1,-2){};
                \draw[Edge](1)--(2);
                \draw[Edge](2)--(3);
                \draw[Edge](2)--(4);
            \end{tikzpicture}
        }
   }
   =
   \sum_{\ell_1 \in \N}
   \Letter^\GenA_{0^{\ell_1}},
\end{equation}
\begin{equation}
    \bar{\Realization}_{\PositionsAlphabet \App \Signature} \App \BasisE_{\alpha_2}
    =
    \Realization_{\PositionsAlphabet \App \Signature} \App
    \BasisE_{
        \scalebox{0.75}{
            \begin{tikzpicture}[Centering,scale=0.4]
                \node[Leaf](1)at(1,0){};
                \node[Node](2)at(1,-1){$\GenA$};
                \node[Node](3)at(0.5,-2){$\GenA$};
                \node[Leaf](4)at(0,-3){};
                \node[Leaf](5)at(1,-3){};
                \node[Leaf](6)at(2,-2){};
                \draw[Edge](1)--(2);
                \draw[Edge](2)--(3);
                \draw[Edge](2)--(6);
                \draw[Edge](3)--(4);
                \draw[Edge](3)--(5);
            \end{tikzpicture}
        }
    }
    +
    \Realization_{\PositionsAlphabet \App \Signature} \App
    \BasisE_{
        \scalebox{0.75}{
            \begin{tikzpicture}[Centering,scale=0.4]
                \node[Leaf](1)at(1,0){};
                \node[Node](2)at(1,-1){$\GenA$};
                \node[Leaf](3)at(0,-2){};
                \node[Node](4)at(1.5,-2){$\GenA$};
                \node[Leaf](5)at(1,-3){};
                \node[Leaf](6)at(2,-3){};
                \draw[Edge](1)--(2);
                \draw[Edge](2)--(3);
                \draw[Edge](2)--(4);
                \draw[Edge](4)--(5);
                \draw[Edge](4)--(6);
            \end{tikzpicture}
        }
    }
    =
    \sum_{\ell_1, \ell_2 \in \N}
    \Letter^\GenA_{0^{\ell_1}} \; \Letter^\GenA_{0^{\ell_1} 1 0^{\ell_2}}
    +
    \sum_{\ell_1, \ell_2 \in \N}
    \Letter^\GenA_{0^{\ell_1}} \; \Letter^\GenA_{0^{\ell_1} 2 0^{\ell_2}},
\end{equation}
\begin{align}
    \bar{\Realization}_{\PositionsAlphabet \App \Signature} \App \BasisE_{\alpha_3}
    & =
    \Realization_{\PositionsAlphabet \App \Signature} \App
    \BasisE_{
        \scalebox{0.75}{
            \begin{tikzpicture}[Centering,scale=0.4]
                \node[Leaf](1)at(1.5,0){};
                \node[Node](2)at(1.5,-1){$\GenA$};
                \node[Node](3)at(1,-2){$\GenA$};
                \node[Node](4)at(0.5,-3){$\GenA$};
                \node[Leaf](5)at(0,-4){};
                \node[Leaf](6)at(1,-4){};
                \node[Leaf](7)at(2,-3){};
                \node[Leaf](8)at(2.5,-2){};
                \draw[Edge](1)--(2);
                \draw[Edge](2)--(3);
                \draw[Edge](2)--(8);
                \draw[Edge](3)--(4);
                \draw[Edge](3)--(7);
                \draw[Edge](4)--(5);
                \draw[Edge](4)--(6);
            \end{tikzpicture}
        }
    }
    +
    \Realization_{\PositionsAlphabet \App \Signature} \App
    \BasisE_{
        \scalebox{0.75}{
            \begin{tikzpicture}[Centering,scale=0.4]
                \node[Leaf](1)at(1.5,0){};
                \node[Node](2)at(1.5,-1){$\GenA$};
                \node[Node](3)at(1,-2){$\GenA$};
                \node[Leaf](4)at(0,-3){};
                \node[Node](5)at(1.5,-3){$\GenA$};
                \node[Leaf](6)at(1,-4){};
                \node[Leaf](7)at(2,-4){};
                \node[Leaf](8)at(2.5,-2){};
                \draw[Edge](1)--(2);
                \draw[Edge](2)--(3);
                \draw[Edge](2)--(8);
                \draw[Edge](3)--(4);
                \draw[Edge](3)--(5);
                \draw[Edge](5)--(6);
                \draw[Edge](5)--(7);
            \end{tikzpicture}
        }
    }
    +
    \Realization_{\PositionsAlphabet \App \Signature} \App
    \BasisE_{
        \scalebox{0.75}{
            \begin{tikzpicture}[Centering,xscale=0.35,yscale=0.4]
                \node[Leaf](1)at(1.5,0){};
                \node[Node](2)at(1.5,-1){$\GenA$};
                \node[Node](3)at(0.5,-2){$\GenA$};
                \node[Leaf](4)at(0,-3){};
                \node[Leaf](5)at(1,-3){};
                \node[Node](6)at(2.5,-2){$\GenA$};
                \node[Leaf](7)at(2,-3){};
                \node[Leaf](8)at(3,-3){};
                \draw[Edge](1)--(2);
                \draw[Edge](2)--(3);
                \draw[Edge](2)--(6);
                \draw[Edge](3)--(4);
                \draw[Edge](3)--(5);
                \draw[Edge](6)--(7);
                \draw[Edge](6)--(8);
            \end{tikzpicture}
        }
    }
    +
    \Realization_{\PositionsAlphabet \App \Signature} \App
    \BasisE_{
        \scalebox{0.75}{
            \begin{tikzpicture}[Centering,scale=0.4]
                \node[Leaf](1)at(1.5,0){};
                \node[Node](2)at(1.5,-1){$\GenA$};
                \node[Leaf](3)at(0.5,-2){};
                \node[Node](4)at(2,-2){$\GenA$};
                \node[Node](5)at(1.5,-3){$\GenA$};
                \node[Leaf](6)at(1,-4){};
                \node[Leaf](7)at(2,-4){};
                \node[Leaf](8)at(3,-3){};
                \draw[Edge](1)--(2);
                \draw[Edge](2)--(3);
                \draw[Edge](2)--(4);
                \draw[Edge](4)--(5);
                \draw[Edge](4)--(8);
                \draw[Edge](5)--(6);
                \draw[Edge](5)--(7);
            \end{tikzpicture}
        }
    }
    +
    \Realization_{\PositionsAlphabet \App \Signature} \App
    \BasisE_{
        \scalebox{0.75}{
            \begin{tikzpicture}[Centering,scale=0.4]
                \node[Leaf](1)at(1.5,0){};
                \node[Node](2)at(1.5,-1){$\GenA$};
                \node[Leaf](3)at(0.5,-2){};
                \node[Node](4)at(2,-2){$\GenA$};
                \node[Leaf](5)at(1,-3){};
                \node[Node](6)at(2.5,-3){$\GenA$};
                \node[Leaf](7)at(2,-4){};
                \node[Leaf](8)at(3,-4){};
                \draw[Edge](1)--(2);
                \draw[Edge](2)--(3);
                \draw[Edge](2)--(4);
                \draw[Edge](4)--(5);
                \draw[Edge](4)--(6);
                \draw[Edge](6)--(7);
                \draw[Edge](6)--(8);
            \end{tikzpicture}
        }
    }
    \notag
    \\
    & =
    \sum_{\ell_1, \ell_2, \ell_3 \in \N}
    \Letter^\GenA_{0^{\ell_1}} \; \Letter^\GenA_{0^{\ell_1} 1 0^{\ell_2}} \;
    \Letter^\GenA_{0^{\ell_1} 1 0^{\ell_2} 1 0^{\ell_3}}
    +
    \sum_{\ell_1, \ell_2, \ell_3 \in \N}
    \Letter^\GenA_{0^{\ell_1}} \; \Letter^\GenA_{0^{\ell_1} 1 0^{\ell_2}} \;
    \Letter^\GenA_{0^{\ell_1} 1 0^{\ell_2} 2 0^{\ell_3}}
    \notag
    \\
    &
    \quad +
    \sum_{\ell_1, \ell_2, \ell_3 \in \N}
    \Letter^\GenA_{0^{\ell_1}} \; \Letter^\GenA_{0^{\ell_1} 1 0^{\ell_2}} \;
    \Letter^\GenA_{0^{\ell_1} 2 0^{\ell_3}}
    +
    \sum_{\ell_1, \ell_2, \ell_3 \in \N}
    \Letter^\GenA_{0^{\ell_1}} \; \Letter^\GenA_{0^{\ell_1} 2 0^{\ell_2}} \;
    \Letter^\GenA_{0^{\ell_1} 2 0^{\ell_2} 1 0^{\ell_3}}
    \notag
    \\
    &
    \quad +
    \sum_{\ell_1, \ell_2, \ell_3 \in \N}
    \Letter^\GenA_{0^{\ell_1}} \; \Letter^\GenA_{0^{\ell_1} 2 0^{\ell_2}} \;
    \Letter^\GenA_{0^{\ell_1} 2 0^{\ell_2} 2 0^{\ell_3}}.
\end{align}

\subsubsection{Forest-like alphabet of lengths}
Let $S$ be a set and $\Equiv$ be the related alphabet congruence of $\PositionsAlphabet \App
S$ satisfying $\Letter^s_u \Equiv \Letter^s_{u'}$ if $\Length \App u = \Length \App u'$ for
all $s \in S$ and $u, u' \in \N^*$. The \Def{$S$-forest-like alphabet of lengths} is the
$S$-forest-like alphabet $\LengthsAlphabet \App S := \PositionsAlphabet \App S /_{\Equiv}$.
Since each $\Equiv$-equivalence class contains a unique letter $\Letter^s_{0^\ell}$ with
$\ell \in \N$, we identify $\LengthsAlphabet \App S$ with the set $\Bra{\Letter^s_\ell : s
\in S \text{ and } \ell \in \N}$. By construction of $\LengthsAlphabet \App S$, and by using
this identification,
\begin{enumerate}[label=({\sf \roman*})]
    \item the root relation $\RootRelation^{\LengthsAlphabet \App S}$ satisfies
    \begin{math}
        \RootRelation^{\LengthsAlphabet \App S} = \LengthsAlphabet \App S;
    \end{math}
    \item for any $s \in S$, the $s$-decoration relation
    $\DecorationRelation_s^{\LengthsAlphabet \App S}$ satisfies
    \begin{math}
        \Letter^s_\ell \in \DecorationRelation_s^{\LengthsAlphabet \App S}
    \end{math}
    for any $\ell \in \N$;
    \item for any $j \geq 1$, the $j$-edge relation $\EdgeRelation{j}^{\LengthsAlphabet \App
    S}$ satisfies
    \begin{math}
        \Letter^s_\ell \EdgeRelation{j} \Letter^{s'}_{\ell'}
    \end{math}
    for any $s, s' \in S$ and $\ell, \ell' \in \N$ whenever $\ell < \ell'$.
\end{enumerate}

Let $\ProjectionLength : \PositionsAlphabet \App S \to \LengthsAlphabet \App S$ be the
canonical projection map associated with $\Equiv$. Under the previous identification, this
map satisfies $\ProjectionLength \App \Letter^s_u = \Letter^s_{\Length \App u}$ for any
$\Letter^s_u \in \PositionsAlphabet \App S$. As explained in
Section~\ref{subsubsec:realizing_maps_quotient_alphabets}, this map is extended as linear
map from $\K \Angle{\PositionsAlphabet \App S}$ to $\K \Angle{\LengthsAlphabet \App S}$.

\begin{Statement}{Proposition}{prop:projection_lengths}
    For any signature $\Signature$,
    \begin{math}
        \Realization_{\LengthsAlphabet \App \Signature}
        = \ProjectionLength
        \circ \Realization_{\PositionsAlphabet \App \Signature}.
    \end{math}
\end{Statement}
\begin{Proof}
    Let $\ForestF \in \Reduced \App \SetForests \App \Signature$ and $n := \Deg \App
    \ForestF$. Let also the sets
    \begin{math}
        X
        := \Bra{w \in \PositionsAlphabet^* :
        w \CompatibleWord{\PositionsAlphabet \App \Signature} \ForestF}
    \end{math}
    and
    \begin{math}
        X'
        := \Bra{w' \in \LengthsAlphabet^* :
        w' \CompatibleWord{\LengthsAlphabet \App \Signature} \ForestF}.
    \end{math}
    We begin by proving that the map $\ProjectionLength$ with domain $X$ and codomain $X'$
    is a bijection.

    First, let $w \in X$. We have thus
    \begin{math}
        w \CompatibleWord{\PositionsAlphabet \App \Signature} \ForestF,
    \end{math}
    so that
    \begin{math}
        w = \Letter^{\Decoration_\ForestF \App 1}_{u_1}
        \dots
        \Letter^{\Decoration_\ForestF \App n}_{u_n}
    \end{math}
    where each $u_i$, $i \in [n]$, is a word on $\N$. Let
    \begin{math}
        w' :=
        \ProjectionLength \App w
        = \Letter^{\Decoration_\ForestF \App 1}_{\Length \App u_1}
        \dots \Letter^{\Decoration_\ForestF \App n}_{\Length \App u_n}.
    \end{math}
    Since all letters of $\LengthsAlphabet \App \Signature$ belong to
    $\RootRelation^{\LengthsAlphabet \App \Signature}$, if $i$ is a root of $\ForestF$, then
    $w' \App i \in \RootRelation^{\LengthsAlphabet \App \Signature}$. Moreover, for any
    internal node $i$ of $\ForestF$,
    \begin{math}
        w' \App i
        \in \DecorationRelation_{\Decoration_\ForestF \App i}^
        {\LengthsAlphabet \App \Signature}.
    \end{math}
    Finally, for all internal nodes $i$ and $i'$ of $\ForestF$, $i \Edge{\ForestF}{j} i'$
    for a $j \geq 1$ implies that
    \begin{math}
        \Letter^{\Decoration_\ForestF \App i}_{u_i} 
        \EdgeRelation{j}^{\PositionsAlphabet \App \Signature}
        \Letter^{\Decoration_\ForestF \App i'}_{u_{i'}}.
    \end{math}
    Hence, there is an $\ell \geq 0$ such that $u_{i'} = u_i \Conc j \Conc 0^\ell$, so that
    $\Length \App u_i < \Length \App u_{i'}$. This shows that
    \begin{math}
        w' \App i \EdgeRelation{j}^{\LengthsAlphabet \App \Signature} w' \App i'.
    \end{math}
    All these properties imply that
    \begin{math}
        w' \CompatibleWord{\LengthsAlphabet \App \Signature} \ForestF,
    \end{math}
    so that $\ProjectionLength$ is a well-defined map from $X$ to~$X'$.

    Let $w' \in X'$. We have thus $w' \CompatibleWord{\LengthsAlphabet \App \Signature}
    \ForestF$, so that $w' = \Letter^{\Decoration_\ForestF \App 1}_{\ell_1} \dots
    \Letter^{\Decoration_\ForestF \App n}_{\ell_n}$ where $\ell_i \in \N$ for any $i \in
    [n]$. In particular, for all internal nodes $i$ and $i'$ of $\ForestF$, $i
    \Edge{\ForestF}{j} i'$ implies that
    \begin{math}
        \Letter^{\Decoration_\ForestF \App i}_{\ell_i} 
        \EdgeRelation{j}^{\LengthsAlphabet \App \Signature}
        \Letter^{\Decoration_\ForestF \App i'}_{\ell_{i'}}
    \end{math}
    so that $\ell_i < \ell_{i'}$. Now, let $w := \Letter^{\Decoration_\ForestF \App 1}_{u_1}
    \dots \Letter^{\Decoration_\ForestF \App n}_{u_n}$ be the word on $\PositionsAlphabet
    \App \Signature$ defined as follows. For any $i' \in [n]$, if $i'$ is a root of
    $\ForestF$, then we set $u_{i'} := 0^{\ell_{i'}}$. Otherwise, there is an internal node
    $i$ of $\ForestF$ and a $j \geq 1$ such that $i \Edge{\ForestF}{j} i'$. In this case, we
    set $u_{i'} := u_i \Conc j \Conc 0^{\ell_{i'} - \ell_i - 1}$. Note that this is
    well-defined since as just established, $\ell_{i'} - \ell_i - 1 \geq 0$. By definition
    of the alphabet $\PositionsAlphabet \App \Signature$, we have $w
    \CompatibleWord{\PositionsAlphabet \App \Signature} \ForestF$. Moreover, we have also
    $\ProjectionLength \App w = w'$. Therefore, $\ProjectionLength$ as a map from $X$ to
    $X'$, is surjective. Observe also that by construction, $w$ is the only element of $X$
    having $w'$ as image by $\ProjectionLength$. Hence, $\ProjectionLength$ as a map from
    $X$ to $X'$ is injective.

    The statement of the proposition follows now by
    Proposition~\ref{prop:factorization_realization_congruence}.
\end{Proof}

\begin{Statement}{Lemma}{lem:lengths_triangular}
    For any signature $\Signature$ and any reduced $\Signature$-forest $\ForestF$,
    \begin{equation}
        \Realization_{\LengthsAlphabet \App \Signature} \App \BasisE_\ForestF
        = \ProjectionLength \App \PositionsEncoding \App \ForestF
        +
        \sum_{w \in \PositionsAlphabet \App \Signature^*}
        \Iverson{w \CompatibleWord{\PositionsAlphabet \App \Signature} \ForestF}
        \Iverson{\Weight \App w > \Weight \App \PositionsEncoding \App \ForestF} \;
        \ProjectionLength \App w.
    \end{equation}
\end{Statement}
\begin{Proof}
    This is a direct consequence of Lemma~\ref{lem:positions_triangular} and
    Proposition~\ref{prop:projection_lengths}.
\end{Proof}

Here is an example of the $\LengthsAlphabet \App \SignatureExample$-realization of an
$\SignatureExample$-forest on the $\BasisE$-basis:
\begin{equation}
    \Realization_{\LengthsAlphabet \App \Signature} \App
    \BasisE_{
        \scalebox{0.75}{
            \begin{tikzpicture}[Centering,xscale=0.35,yscale=0.45]
                \node[Leaf](2)at(2.5,-1){};
                \node[Node](3)at(2.5,-2){$\GenC$};
                \node[Node](4)at(1,-3){$\GenB$};
                \node[Leaf](5)at(0.5,-4){};
                \node[Leaf](6)at(1.5,-4){};
                \node[Leaf](7)at(2.5,-3){};
                \node[Node](8)at(4,-3){$\GenA$};
                \node[Node](9)at(4,-4){$\GenC$};
                \node[Leaf](10)at(3,-5){};
                \node[Leaf](11)at(4,-5){};
                \node[Leaf](12)at(5,-5){};
                \node[Leaf](13)at(6.5,-1){};
                \node[Node](14)at(6.5,-2){$\GenB$};
                \node[Node](15)at(6,-3){$\GenA$};
                \node[Leaf](16)at(6,-4){};
                \node[Leaf](17)at(7,-3){};
                \draw[Edge](2)--(3);
                \draw[Edge](3)--(4);
                \draw[Edge](3)--(7);
                \draw[Edge](3)--(8);
                \draw[Edge](4)--(5);
                \draw[Edge](4)--(6);
                \draw[Edge](8)--(9);
                \draw[Edge](9)--(10);
                \draw[Edge](9)--(11);
                \draw[Edge](9)--(12);
                \draw[Edge](13)--(14);
                \draw[Edge](14)--(15);
                \draw[Edge](14)--(17);
                \draw[Edge](15)--(16);
            \end{tikzpicture}
        }
    }
    =
    \sum_{\ell_1, \dots, \ell_6 \in \N}
    \Iverson{\ell_1 < \ell_2}
    \Iverson{\ell_1 < \ell_3 < \ell_4}
    \Iverson{\ell_5 < \ell_6} \;
    \Letter^\GenC_{\ell_1} \;
    \Letter^\GenB_{\ell_2} \;
    \Letter^\GenA_{\ell_3} \;
    \Letter^\GenC_{\ell_4} \;
    \Letter^\GenB_{\ell_5} \;
    \Letter^\GenA_{\ell_6}.
\end{equation}
Observe in particular that since
\begin{equation}
    \Realization_{\LengthsAlphabet \App \Signature} \App
    \BasisE_{
        \scalebox{0.75}{
            \begin{tikzpicture}[Centering,xscale=0.35,yscale=0.45]
                \node[Leaf](2)at(2.5,-1){};
                \node[Node](3)at(2.5,-2){$\GenC$};
                \node[Node](4)at(1,-3){$\GenB$};
                \node[Leaf](5)at(0.5,-4){};
                \node[Leaf](6)at(1.5,-4){};
                \node[Leaf](7)at(2.5,-3){};
                \node[Node](8)at(4,-3){$\GenA$};
                \node[Node](9)at(4,-4){$\GenC$};
                \node[Leaf](10)at(3,-5){};
                \node[Leaf](11)at(4,-5){};
                \node[Leaf](12)at(5,-5){};
                \node[Leaf](13)at(6.5,-1){};
                \node[Node](14)at(6.5,-2){$\GenB$};
                \node[Node](15)at(6,-3){$\GenA$};
                \node[Leaf](16)at(6,-4){};
                \node[Leaf](17)at(7,-3){};
                \draw[Edge](2)--(3);
                \draw[Edge](3)--(4);
                \draw[Edge](3)--(7);
                \draw[Edge](3)--(8);
                \draw[Edge](4)--(5);
                \draw[Edge](4)--(6);
                \draw[Edge](8)--(9);
                \draw[Edge](9)--(10);
                \draw[Edge](9)--(11);
                \draw[Edge](9)--(12);
                \draw[Edge](13)--(14);
                \draw[Edge](14)--(15);
                \draw[Edge](14)--(17);
                \draw[Edge](15)--(16);
            \end{tikzpicture}
        }
    }
    =
    \Realization_{\LengthsAlphabet \App \Signature} \App
    \BasisE_{
        \scalebox{0.75}{
            \begin{tikzpicture}[Centering,xscale=0.35,yscale=0.45]
                \node[Leaf](2)at(3,-1){};
                \node[Node](3)at(3,-2){$\GenC$};
                \node[Node](4)at(1,-3){$\GenB$};
                \node[Leaf](5)at(0.5,-4){};
                \node[Leaf](6)at(1.5,-4){};
                \node[Node](7)at(3,-3){$\GenA$};
                \node[Node](8)at(3,-4){$\GenC$};
                \node[Leaf](9)at(2,-5){};
                \node[Leaf](10)at(3,-5){};
                \node[Leaf](11)at(4,-5){};
                \node[Leaf](12)at(4.5,-3){};
                \node[Leaf](13)at(6.5,-1){};
                \node[Node](14)at(6.5,-2){$\GenB$};
                \node[Leaf](15)at(6,-3){};
                \node[Node](16)at(7,-3){$\GenA$};
                \node[Leaf](17)at(7,-4){};
                \draw[Edge](2)--(3);
                \draw[Edge](3)--(4);
                \draw[Edge](3)--(7);
                \draw[Edge](3)--(12);
                \draw[Edge](4)--(5);
                \draw[Edge](4)--(6);
                \draw[Edge](7)--(8);
                \draw[Edge](8)--(9);
                \draw[Edge](8)--(10);
                \draw[Edge](8)--(11);
                \draw[Edge](13)--(14);
                \draw[Edge](14)--(15);
                \draw[Edge](14)--(16);
                \draw[Edge](16)--(17);
            \end{tikzpicture}
        }
    },
\end{equation}
the map $\Realization_{\LengthsAlphabet \App \Signature}$ is not injective. Nevertheless,
this map has interesting properties. The quotient $\NaturalHopfAlgebra \App \SetTerms \App
\Signature / \Ker \App \Realization_{\LengthsAlphabet \App \Signature}$ will be studied in
Section~\ref{subsubsec:quotient_nck} and is linked with decorated versions of noncommutative
Connes-Kreimer Hopf algebras. In the same vein, we shall see in
Section~\ref{subsubsec:polynomial_realization_n_mas} that
$\bar{\Realization}_{\LengthsAlphabet \App \Signature}$ is a polynomial realization
of~$\NCFdB$.

\section{Links with other Hopf algebras} \label{sec:links_hopf_algebras}
In this last section, we establish links between natural Hopf algebras of some operads and
other Hopf algebras by using the polynomial realization of natural Hopf algebras of operads
introduced in Section~\ref{sec:polynomial_realization}.

\subsection{Decorated word quasi-symmetric functions}
\label{subsec:decorated_word_quasi_symmetric_functions}
We show in this part that the $\LengthsAlphabet \App \Signature$-realizations of reduced
$\Signature$-forests on the $\BasisE$-basis span certain generalized word quasi-symmetric
functions.

\subsubsection{Packed decorated words}
Let $D$ be a set. A \Def{$D$-decorated letter} is a pair $(k, \GenD)$, denoted by $k^\GenD$,
where $\GenD \in D$ and $k$ is a positive integer. We call $k$ (resp.\ $\GenD$) the
\Def{value} (resp.\ the \Def{decoration}) of $k^\GenD$. The set of $D$-decorated letters is
denoted by $\SetDecoratedLetters \App D$. Each word on $\SetDecoratedLetters \App D$ is a
\Def{$D$-decorated word}. Given a $D$-decorated word $u \in \SetDecoratedLetters \App D^*$,
for any $i \in [\Length \App u]$, we denote by $\Value_u \App i$ the value of $u \App i$ and
by $\Decoration_u \App i$ the decoration of $u \App i$. The \Def{packing} of $u$ is the
$D$-decorated word $\Pack \App u$ obtained by replacing each $D$-decorated letter
$k_1^\GenD$ of $u$ by $k_2^\GenD$ where $k_2$ is the number of different values less than or
equal to $k_1$ among the $D$-decorated letters of $u$. When $u$ is a $D$-decorated word such
that $\Pack \App u = u$, $u$ is \Def{packed}.

For instance, for $D := \{\GenA, \GenB, \GenC\}$, $u : = 4^\GenB 2^\GenB 3^\GenA 4^\GenB
4^\GenC 6^\GenC 3^\GenA$ is a $D$-decorated word of length $7$ satisfying $\Value_u \App 5 =
4$, $\Decoration_u \App 5 = \GenC$, and $\Pack \App u = 3^\GenB 1^\GenB 2^\GenA 3^\GenB
3^\GenC 4^\GenC 2^\GenA$.

\subsubsection{Decorated word quasi-symmetric functions}
Let $D$ be a set and $\Monomial : \SetDecoratedLetters \App D^* \to \LengthsAlphabet \App
D^*$ be the map defined for any $u \in \SetDecoratedLetters \App D^*$ by $\Monomial \App u :
= \Letter^{\Decoration_u \App 1}_{\Value_u \App 1} \dots \Letter^{\Decoration_u \App \Length
\App u}_{\Value_u \App \Length \App u}$. For instance, for $D := \{\GenA, \GenB, \GenC\}$,
we have $\Monomial \App 2^\GenA 1^\GenA 1^\GenB 4^\GenC 2^\GenB = \Letter^\GenA_2 \;
\Letter^\GenA_1 \; \Letter^\GenB_1 \; \Letter^\GenC_4 \; \Letter^\GenB_2$. Moreover, for any
packed $D$-decorated word $u$, let $\BasisM_u$ be the $\LengthsAlphabet \App D$-polynomial
defined by
\begin{equation} \label{equ:realization_wqsym}
    \BasisM_u
    := \sum_{v \in \SetDecoratedLetters \App D^*} \Iverson{\Pack \App v = u} \;
    \Monomial \App v.
\end{equation}
For instance, for the same set $D$ as in the previous example,
\begin{equation}
    \BasisM_{2^\GenB 1^\GenC 1^\GenC 3^\GenA}
    = \sum_{\ell_1, \ell_2, \ell_3, \ell_4 \in \N}
    \Iverson{\ell_2 = \ell_3 < \ell_1 < \ell_4} \;
    \Letter^\GenB_{\ell_1} \;
    \Letter^\GenC_{\ell_2} \;
    \Letter^\GenC_{\ell_3} \;
    \Letter^\GenA_{\ell_4}.
\end{equation}
Let $\WQSym \App D$ be the $\K$-linear span of the set $\Bra{\BasisM_u : u \in \Pack \App
\Par{\SetDecoratedLetters \App D^*}}$. We call \Def{$D$-decorated word quasi-symmetric
function} each element of $\WQSym \App D$. When $D$ is a singleton, $\WQSym \App D$ is
isomorphic to $\WQSym$, the Hopf algebra of word quasi-symmetric functions introduced
in~\cite{Hiv99,NT06}. This space $\WQSym \App D$, enriching $\WQSym$ with colors (called
``decorations'' here), is analogous to some similar constructions presented in~\cite{NT10},
enriching the Hopf algebra $\NCSF$ of noncommutative symmetric functions~\cite{GKLLRT95},
the Hopf algebra $\FQSym$ of free quasi-symmetric functions~\cite{MR95,DHT02}, and the Hopf
algebra $\PQSym$ of parking quasi-symmetric functions~\cite{NT07} in a similar way.

\begin{Statement}{Theorem}{thm:map_to_wqsym}
    For any signature $\Signature$ and any reduced $\Signature$-forest $\ForestF$,
    $\Realization_{\LengthsAlphabet \App \Signature} \App \BasisE_\ForestF$ is an
    $\Signature$-decorated word quasi-symmetric function. More precisely,
    \begin{equation} \label{equ:map_to_wqsym_1}
        \Realization_{\LengthsAlphabet \App \Signature}
        \App \BasisE_\ForestF
        =
        \sum_{u \in \Pack \App \Par{\SetDecoratedLetters \App \Signature^*}}
        \Iverson{
            \Monomial \App u
            \CompatibleWord{\LengthsAlphabet \App \Signature} \ForestF
        } \;
        \BasisM_u.
    \end{equation}
\end{Statement}
\begin{Proof}
    First of all, observe that, by the general definition~\eqref{equ:realization_forest} of
    the map $\Realization_A : \NaturalHopfAlgebra \App \SetTerms \App \Signature \to \K
    \Angle{A}$, each monomial appearing in the left-hand side of~\eqref{equ:map_to_wqsym_1}
    admits $1$ as coefficient. Besides, by the definition~\eqref{equ:realization_wqsym}, for
    any $u \in \Pack \App \Par{\SetDecoratedLetters \App \Signature^*}$, each monomial
    appearing in $\BasisM_u$ admits $1$ as coefficient. Moreover, it is clear that if $u$
    and $u'$ are different packed $D$-decorated words, then, by definition of the map
    $\Pack$, the supports of $\BasisM_u$ and $\BasisM_{u'}$ are disjoint, For these reasons,
    to establish~\eqref{equ:map_to_wqsym_1}, it is enough to prove that the supports of its
    left-hand side and of its right-hand side are equal.

    Assume that $w = \Letter^{\GenS_1}_{\ell_1} \dots \Letter^{\GenS_n}_{\ell_n}$ is a
    monomial appearing in the left-hand side of~\eqref{equ:map_to_wqsym_1}, where $n := \Deg
    \App \ForestF$ and for any $i \in [n]$, $\GenS_i \in \Signature$ and $\ell_i \in \N$.
    Let the packed $D$-decorated word $u := \Pack \App \ell_1^{\GenS_1} \dots
    \ell_n^{\GenS_n}$. Since $w \CompatibleWord{\LengthsAlphabet \App \Signature}
    \ForestF$, and since the map $\Pack$ preserves the decorations of the letters and
    preserves the relative order between the values of the letters, we have $\Monomial \App
    u \CompatibleWord{\LengthsAlphabet \App \Signature} \ForestF$. Now, by construction of
    $u$ and by definition of $\BasisM_u$, the monomial $w$ appears in $\BasisM_u$. Hence,
    $w$ appears in the right-hand side of~\eqref{equ:map_to_wqsym_1}. The inverse property,
    consisting in the fact that any monomial appearing in the right-hand side
    of~\eqref{equ:map_to_wqsym_1} appears also in its left-hand side is shown by a similar
    reasoning carried out in the opposite direction.
\end{Proof}

For instance, in $\NaturalHopfAlgebra \App \SetTerms \App \SignatureExample$,
\begin{equation}
    \Realization_{\LengthsAlphabet \App \SignatureExample}
    \App
    \BasisE_{
        \scalebox{0.75}{
            \begin{tikzpicture}[Centering,xscale=0.4,yscale=0.45]
                \node[Leaf](2)at(1.25,-1){};
                \node[Node](3)at(1.25,-2){$\GenC$};
                \node[Node](4)at(0,-3){$\GenA$};
                \node[Leaf](5)at(0,-4){};
                \node[Leaf](6)at(1.25,-3){};
                \node[Node](7)at(2.5,-3){$\GenB$};
                \node[Leaf](8)at(2,-4){};
                \node[Leaf](9)at(3,-4){};
                \draw[Edge](2)--(3);
                \draw[Edge](3)--(4);
                \draw[Edge](3)--(6);
                \draw[Edge](3)--(7);
                \draw[Edge](4)--(5);
                \draw[Edge](7)--(8);
                \draw[Edge](7)--(9);
            \end{tikzpicture}
        }
    }
    =
    \BasisM_{1^\GenC 2^\GenA 2^\GenB} + \BasisM_{1^\GenC 2^\GenA 3^\GenB}
    + \BasisM_{1^\GenC 3^\GenA 2^\GenB}.
\end{equation}

\subsection{Connes-Kreimer Hopf algebras of trimmed forests} \label{subsec:connes_kreimer}
We establish in this section a link between natural Hopf algebras $\NaturalHopfAlgebra \App
\SetTerms \App \Signature$ of free operads and decorated versions of noncommutative
Connes-Kreimer Hopf algebras.

\subsubsection{Trimmed forests}
For any signature $\Signature$, an \Def{$\Signature$-trimmed forest} is a word of nonempty
planar rooted trees such that each node is decorated by an element $\GenS$ of $\Signature$
and has at most $\Arity \App \GenS$ children. By definition, in such forests there are no
leaves and thus, each node is internal. For instance,
\begin{equation} \label{equ:trimmed_forest}
    \scalebox{.75}{
        \begin{tikzpicture}[Centering,xscale=0.45,yscale=0.5]
            \node[Node,MarkB](2)at(1,-1){$\GenC$};
            \node[Node,MarkB](3)at(0,-2){$\GenA$};
            \node[Node,MarkB](4)at(2,-2){$\GenC$};
            \node[Node,MarkB](5)at(1,-3){$\GenC$};
            \node[Node,MarkB](6)at(2,-3){$\GenB$};
            \node[Node,MarkB](7)at(3,-3){$\GenB$};
            \node[Node,MarkB](8)at(4,-1){$\GenB$};
            \node[Node,MarkB](9)at(5,-1){$\GenA$};
            \node[Node,MarkB](10)at(5,-2){$\GenB$};
            \draw[Edge](2)--(3);
            \draw[Edge](2)--(4);
            \draw[Edge](4)--(5);
            \draw[Edge](4)--(6);
            \draw[Edge](4)--(7);
            \draw[Edge](9)--(10);
        \end{tikzpicture}
    }
\end{equation}
is an $\SignatureExample$-trimmed forest. Moreover, for any $\Signature$-forest $\ForestF$,
let $\Trim \App \ForestF$ be the $\Signature$-trimmed forest obtained by removing the leaves
of $\ForestF$. For instance, on the signature $\SignatureExample$,
\begin{equation} \label{equ:trim}
    \Trim \; \App \;
    \scalebox{.75}{
        \begin{tikzpicture}[Centering,xscale=0.5,yscale=0.45]
            \node[Leaf](2)at(0.5,-1){};
            \node[Node](3)at(0.5,-2){$\GenB$};
            \node[Node](4)at(0,-3){$\GenA$};
            \node[Leaf](5)at(0,-4){};
            \node[Leaf](6)at(1,-3){};
            \node[Leaf](7)at(2,-1){};
            \node[Leaf](8)at(2,-2){};
            \node[Leaf](9)at(4,-1){};
            \node[Node](10)at(4,-2){$\GenC$};
            \node[Leaf](11)at(3,-3){};
            \node[Node](12)at(5,-3){$\GenB$};
            \node[Leaf](13)at(5,-4){};
            \node[Node](14)at(4,-4){$\GenA$};
            \node[Leaf](15)at(4,-5){};
            \node[Node](16)at(6,-4){$\GenC$};
            \node[Leaf](17)at(5.5,-5){};
            \node[Leaf](18)at(6,-5){};
            \node[Leaf](19)at(6.5,-5){};
            \node[Leaf](20)at(4,-3){};
            \draw[Edge](2)--(3);
            \draw[Edge](3)--(4);
            \draw[Edge](3)--(6);
            \draw[Edge](4)--(5);
            \draw[Edge](7)--(8);
            \draw[Edge](9)--(10);
            \draw[Edge](10)--(11);
            \draw[Edge](10)--(12);
            \draw[Edge](12)--(13);
            \draw[Edge](12)--(14);
            \draw[Edge](12)--(16);
            \draw[Edge](14)--(15);
            \draw[Edge](16)--(17);
            \draw[Edge](16)--(18);
            \draw[Edge](16)--(19);
            \draw[Edge](10)--(20);
        \end{tikzpicture}
    }
    \enspace = \enspace
    \scalebox{.75}{
        \begin{tikzpicture}[Centering,xscale=0.5,yscale=0.5]
            \node[Node,MarkB](2)at(0,-1){$\GenB$};
            \node[Node,MarkB](3)at(0,-2){$\GenA$};
            \node[Node,MarkB](4)at(1.5,-1){$\GenC$};
            \node[Node,MarkB](5)at(1.5,-2){$\GenB$};
            \node[Node,MarkB](6)at(1,-3){$\GenA$};
            \node[Node,MarkB](7)at(2,-3){$\GenC$};
            \draw[Edge](2)--(3);
            \draw[Edge](4)--(5);
            \draw[Edge](5)--(6);
            \draw[Edge](5)--(7);
        \end{tikzpicture}
    }.
\end{equation}
Most definitions and properties concerning $\Signature$-forests introduced mainly in
Sections~\ref{subsubsec:terms}, \ref{subsubsec:forests},
and~\ref{subsubsec:coproduct_description} apply to $\Signature$-trimmed forests by
stipulating that an $\Signature$-trimmed forest $\TrimmedForestF$ satisfies a property $P$
if all reduced $\Signature$-forests $\ForestF$ such that $\Trim \App \ForestF =
\TrimmedForestF$ satisfy $P$. For instance, the $\SignatureExample$-trimmed forest
of~\eqref{equ:trimmed_forest} has degree $9$, its node $1$ is decorated by $\GenC$, its node
$5$ is decorated by $\GenB$, and~$\Height_\TrimmedForestF \App 5 = 2$.

The \Def{charge} $\Charge \App \TrimmedForestF$ of an $\Signature$-trimmed forest
$\TrimmedForestF$ is the positive integer defined recursively as follows. If $\Length \App
\TrimmedForestF \ne 1$, then
\begin{equation} \label{equ:charge_1}
    \Charge \App \TrimmedForestF
    := \prod_{i \in [\Length \App \TrimmedForestF]} \Charge \App \TrimmedForestF \App i.
\end{equation}
Otherwise, $\TrimmedForestF$ decomposes into a single root decorated by $\GenS \in
\Signature$ which is attached to an $\Signature$-trimmed forest $\TrimmedForestF'$. In this
case,
\begin{equation} \label{equ:charge_2}
    \Charge \App \TrimmedForestF
    := \binom{\Arity \App \GenS}{\Length \App \TrimmedForestF'}
    \Charge \App \TrimmedForestF'.
\end{equation}
For instance, by denoting by $\TrimmedForestF_1$ (resp.\ $\TrimmedForestF_2$) the
$\Signature$-trimmed forest of~\eqref{equ:trimmed_forest} (resp.\ the right-hand side
of~\eqref{equ:trim}), we have $\Charge \App \TrimmedForestF_1 = 3$ (resp.\ $\Charge \App
\TrimmedForestF_2 = 6$).

\begin{Statement}{Lemma}{lem:number_trimmed}
    For any signature $\Signature$ and any $\Signature$-trimmed forest $\TrimmedForestF$,
    the charge of $\TrimmedForestF$ is the cardinality of the set of reduced
    $\Signature$-forests $\ForestF$ such that $\Trim \App \ForestF = \TrimmedForestF$.
\end{Statement}
\begin{Proof}
    We proceed by structural induction on the set of $\Signature$-trimmed forests. Let
    $\TrimmedForestF$ be an $\Signature$-trimmed forest and $\ForestF$ be a reduced
    $\Signature$-forest such that $\Trim \App \ForestF = \TrimmedForestF$.

    If $\Length \App \TrimmedForestF \ne 1$, by definition of the map $\Trim$, $\ForestF$ is
    such that $\Length \App \ForestF = \Length \App \TrimmedForestF$ and $\Trim \App
    \ForestF \App i = \TrimmedForestF \App i$ for any $i \in [\Length \App \ForestF]$. By
    induction hypothesis, there are exactly $\Charge \App \TrimmedForestF \App i$ different
    reduced $\Signature$-forests $\ForestF \App i$ satisfying the previous property. For
    this reason, by~\eqref{equ:charge_1}, the total number of such reduced
    $\Signature$-forests $\ForestF$ is~$\Charge \App \TrimmedForestF$.

    Otherwise, $\Length \App \TrimmedForestF = 1$. In this case, $\TrimmedForestF$
    decomposes into a single root decorated by $\GenS \in \Signature$ which is attached to
    an $\Signature$-trimmed forest $\TrimmedForestF'$. By definition of the map $\Trim$,
    $\ForestF$ is such that $\Length \App \ForestF = 1$ and $\ForestF$ decomposes into a
    single root decorated by $\GenS$ which is attached to an $\Signature$-forest $\ForestF'$
    such that $\Trim \App \ForestF' = \TrimmedForestF'$. Since $\ForestF'$ is not
    necessarily reduced, $\ForestF'$ is obtained by shuffling the reduced
    $\Signature$-forest $\Reduced \App \ForestF'$ with the $\Signature$-forest made of
    $\Length \App \ForestF' - \Length \App \Reduced \App \ForestF'$ occurrences of $\Leaf$.
    The number of such configurations is the binomial of $\Length \App \ForestF' = \Arity
    \App \GenS$ choose $\Length \App \ForestF' - \Length \App \Reduced \App \ForestF' =
    \Arity \App \GenS - \Length \App \TrimmedForestF'$. Moreover, by induction hypothesis,
    there are exactly $\Charge \App \TrimmedForestF'$ different possibilities for the
    $\Signature$-forest $\ForestF'$. For these reasons, by~\eqref{equ:charge_2}, the total
    number of such reduced $\Signature$-forests $\ForestF$ is~$\Charge \App
    \TrimmedForestF$.
\end{Proof}

\subsubsection{Connes-Kreimer Hopf algebras}
For any signature $\Signature$, let $\NCK \App \Signature$ be the $\K$-linear span of the
set $\Trim \App \SetForests \App \Signature$. The \Def{elementary basis} of $\NCK \App
\Signature$ is the set $\Bra{\BasisE_\TrimmedForestF : \TrimmedForestF \in \Trim \App
\SetForests \App \Signature}$. This vector space is endowed with the product $\Product$
defined, for any $\TrimmedForestF_1, \TrimmedForestF_2 \in \Trim \App \SetForests \App
\Signature$, by
\begin{math}
    \BasisE_{\TrimmedForestF_1} \Product \BasisE_{\TrimmedForestF_2}
    := \BasisE_{\TrimmedForestF_1 \Conc \TrimmedForestF_2}
\end{math}
and with the coproduct $\Coproduct$ defined, for any $\TrimmedForestF \in \Trim \App
\SetForests \App \Signature$ by
\begin{equation} \label{equ:coproduct_nck}
    \Coproduct \App \BasisE_\TrimmedForestF
    :=
    \sum_{I_1, I_2 \subseteq \Han{\Deg \App \TrimmedForestF}}
    \Iverson{\Par{I_1, I_2} \AdmissibleSet \TrimmedForestF}
    \;
    \BasisE_{\TrimmedForestF \App I_1} \otimes \BasisE_{\TrimmedForestF \App I_2}.
\end{equation}
Note that we use in~\eqref{equ:coproduct_nck} the notion of restriction of an
$\Signature$-trimmed forest on a set of nodes and the notion of $\TrimmedForestF$-admissible
pair of sets, both introduced in Section~\ref{subsubsec:coproduct_description}. For
instance,
\begin{equation}
    \Coproduct \App
    \BasisE_{
        \scalebox{.75}{
            \begin{tikzpicture}[Centering,xscale=0.5,yscale=0.5]
                \node[Node,MarkB](1)at(0.5,0){$\GenC$};
                \node[Node,MarkB](2)at(0,-1){$\GenB$};
                \node[Node,MarkB](3)at(1,-1){$\GenA$};
                \draw[Edge](1)--(2);
                \draw[Edge](1)--(3);
            \end{tikzpicture}
        }
    }
    =
    \BasisE_\epsilon
    \otimes
    \BasisE_{
        \scalebox{.75}{
            \begin{tikzpicture}[Centering,xscale=0.5,yscale=0.5]
                \node[Node,MarkB](1)at(0.5,0){$\GenC$};
                \node[Node,MarkB](2)at(0,-1){$\GenB$};
                \node[Node,MarkB](3)at(1,-1){$\GenA$};
                \draw[Edge](1)--(2);
                \draw[Edge](1)--(3);
            \end{tikzpicture}
        }
    }
    +
    \BasisE_{
        \scalebox{.75}{
            \begin{tikzpicture}[Centering,xscale=0.5,yscale=0.5]
                \node[Node,MarkB](1)at(0,0){$\GenC$};
            \end{tikzpicture}
        }
    }
    \otimes
    \BasisE_{
        \scalebox{.75}{
            \begin{tikzpicture}[Centering,xscale=0.5,yscale=0.5]
                \node[Node,MarkB](2)at(0,-1){$\GenB$};
                \node[Node,MarkB](3)at(1,-1){$\GenA$};
            \end{tikzpicture}
        }
    }
    +
    \BasisE_{
        \scalebox{.75}{
            \begin{tikzpicture}[Centering,xscale=0.5,yscale=0.5]
                \node[Node,MarkB](1)at(0,0){$\GenC$};
                \node[Node,MarkB](2)at(0,-1){$\GenB$};
                \draw[Edge](1)--(2);
            \end{tikzpicture}
        }
    }
    \otimes
    \BasisE_{
        \scalebox{.75}{
            \begin{tikzpicture}[Centering,xscale=0.5,yscale=0.5]
                \node[Node,MarkB](1)at(0,0){$\GenA$};
            \end{tikzpicture}
        }
    }
    +
    \BasisE_{
        \scalebox{.75}{
            \begin{tikzpicture}[Centering,xscale=0.5,yscale=0.5]
                \node[Node,MarkB](1)at(0,0){$\GenC$};
                \node[Node,MarkB](2)at(0,-1){$\GenA$};
                \draw[Edge](1)--(2);
            \end{tikzpicture}
        }
    }
    \otimes
    \BasisE_{
        \scalebox{.75}{
            \begin{tikzpicture}[Centering,xscale=0.5,yscale=0.5]
                \node[Node,MarkB](1)at(0,0){$\GenB$};
            \end{tikzpicture}
        }
    }
    +
    \BasisE_{
        \scalebox{.75}{
            \begin{tikzpicture}[Centering,xscale=0.5,yscale=0.5]
                \node[Node,MarkB](1)at(0.5,0){$\GenC$};
                \node[Node,MarkB](2)at(0,-1){$\GenB$};
                \node[Node,MarkB](3)at(1,-1){$\GenA$};
                \draw[Edge](1)--(2);
                \draw[Edge](1)--(3);
            \end{tikzpicture}
        }
    }
    \otimes
    \BasisE_\epsilon.
\end{equation}

This Hopf algebra $\NCK \App \Signature$ is, in fact, a Hopf subalgebra of the
noncommutative Connes-Kreimer Hopf algebra of decorated forests introduced in~\cite{Foi02a}
and~\cite{Foi02b}. In the latter Hopf algebra, there are no restrictions on the arities of
the internal nodes of the forests, unlike in $\NCK \App \Signature$. Note that the first
instance of such Hopf algebras, involving non-decorated forests, appears in~\cite{CK98}.

\subsubsection{A quotient of $\NaturalHopfAlgebra \App \SetTerms \App \Signature$}
\label{subsubsec:quotient_nck}
Let $\Signature$ be a signature. We now study the kernel of the map
$\Realization_{\LengthsAlphabet \App \Signature}$ on the domain $\NaturalHopfAlgebra \App
\SetTerms \App \Signature$. We shall show that the image of this map is isomorphic to~$\NCK
\App \Signature$.

\begin{Statement}{Lemma}{lem:trimmed_equivalence}
    For any signature $\Signature$ and any reduced $\Signature$-forests $\ForestF_1$ and
    $\ForestF_2$, we have
    \begin{math}
        \Realization_{\LengthsAlphabet \App \Signature} \App \BasisE_{\ForestF_1}
        =
        \Realization_{\LengthsAlphabet \App \Signature} \App \BasisE_{\ForestF_2}
    \end{math}
    if and only if $\Trim \App \ForestF_1 = \Trim \App \ForestF_2$.
\end{Statement}
\begin{Proof}
    Let us first prove the sufficient condition of the statement by induction on $n$, the
    common degree of both $\ForestF_1$ and $\ForestF_2$. If $n = 0$, the property is
    immediate since $\ForestF_1 = \epsilon = \ForestF_2$. Otherwise, $n \geq 1$, and let
    $\ForestF'_1$ (resp.\ $\ForestF'_2$) be the $\Signature$-forest obtained by replacing
    the greatest internal node of $\ForestF_1$ (resp.\ $\ForestF_2$) by a leaf. Since by
    hypothesis,
    \begin{math}
        \Realization_{\LengthsAlphabet \App \Signature} \App \BasisE_{\ForestF_1}
        =
        \Realization_{\LengthsAlphabet \App \Signature} \App \BasisE_{\ForestF_2},
    \end{math}
    by Lemma~\ref{lem:lengths_triangular}, we have in particular that
    \begin{math}
        \ProjectionLength \App \PositionsEncoding \App \ForestF_1
        = \ProjectionLength \App \PositionsEncoding \App \ForestF_2.
    \end{math}
    Therefore,
    \begin{equation}
        \ProjectionLength \App \PositionsEncoding \App \ForestF_1
        = \ProjectionLength \App \PositionsEncoding \App \ForestF_1'
        \; \Conc \;
        \Letter^{\Decoration_{\ForestF_1} \App n}_{\Height_{\ForestF_1} \App n}
        = \ProjectionLength \App \PositionsEncoding \App \ForestF_2'
        \; \Conc \;
        \Letter^{\Decoration_{\ForestF_2} \App n}_{\Height_{\ForestF_2} \App n}
        = \ProjectionLength \App \PositionsEncoding \App \ForestF_2.
    \end{equation}
    From these observations, we deduce that $\Decoration_{\ForestF_1} \App n =
    \Decoration_{\ForestF_2} \App n$ and $\Height_{\ForestF_1} \App n = \Height_{\ForestF_2}
    \App n$. Moreover, by induction hypothesis, $\Trim \App \ForestF_1' = \Trim \App
    \ForestF_2'$. Since the internal node $n$ is the last one of both $\ForestF_1$ and
    $\ForestF_2$, this implies that $\Trim \App \ForestF_1 = \Trim \App \ForestF_2$ as
    expected.

    Conversely, assume that $\Trim \App \ForestF_1 = \Trim \App \ForestF_2$. In particular,
    this implies that $\ForestF_1$ and $\ForestF_2$ have a same degree $n$, for any $i \in
    [n]$, $\Decoration_{\ForestF_1} \App i = \Decoration_{\ForestF_2} \App i$, and
    $\Height_{\ForestF_1} \App i = \Height_{\ForestF_2} \App i$. Therefore,
    \begin{math}
        \ProjectionLength \App \PositionsEncoding \App \ForestF_1
        = \ProjectionLength \App \PositionsEncoding \App \ForestF_2.
    \end{math}
    Now, again by Lemma~\ref{lem:lengths_triangular},
    \begin{math}
        \Realization_{\LengthsAlphabet \App \Signature} \App \BasisE_{\ForestF_1}
        =
        \Realization_{\LengthsAlphabet \App \Signature} \App \BasisE_{\ForestF_2}.
    \end{math}
\end{Proof}

\begin{Statement}{Lemma}{lem:kernel_nck}
    For any signature $\Signature$, the kernel of the map $\Realization_{\LengthsAlphabet
    \App \Signature}$ is the subspace of $\NaturalHopfAlgebra \App \SetTerms \App
    \Signature$ generated by the elements $\BasisE_{\ForestF_1} - \BasisE_{\ForestF_2}$ such
    that $\ForestF_1$ and $\ForestF_2$ are reduced $\Signature$-forests satisfying $\Trim
    \App \ForestF_1 = \Trim \App \ForestF_2$.
\end{Statement}
\begin{Proof}
    This is a direct consequence of Lemma~\ref{lem:trimmed_equivalence} and of the
    triangularity property of the map $\Realization_{\LengthsAlphabet \App \Signature}$
    exhibited by Lemma~\ref{lem:lengths_triangular}.
\end{Proof}

Let
\begin{math}
    \pi : \NaturalHopfAlgebra \App \SetTerms \App \Signature \to
    \NaturalHopfAlgebra \App \SetTerms \App \Signature
    / \Ker \App \Realization_{\LengthsAlphabet \App \Signature}
\end{math}
be the quotient map associated with $\Realization_{\LengthsAlphabet \App \Signature}$.

\begin{Statement}{Theorem}{thm:map_to_nck}
    For any signature $\Signature$, the linear map
    \begin{math}
        \phi :
        \NaturalHopfAlgebra \App \SetTerms \App \Signature
        / \Ker \App \Realization_{\LengthsAlphabet \App \Signature}
        \to \NCK \App \Signature
    \end{math}
    satisfying $\phi \App \pi \App \BasisE_\ForestF = \BasisE_{\Trim \App \ForestF}$ for any
    reduced $\Signature$-forest $\ForestF$ is a Hopf algebra isomorphism.
\end{Statement}
\begin{Proof}
    By Lemma~\ref{lem:kernel_nck}, as a vector space,
    \begin{math}
        \NaturalHopfAlgebra \App \SetTerms \App \Signature
        / \Ker \App \Realization_{\LengthsAlphabet \App \Signature}
    \end{math}
    is isomorphic to the linear span of the set $\Trim \App \SetForests \App \Signature$.
    Therefore, $\phi$ is an isomorphism of vector spaces. Moreover, due to the definition of
    the product of $\NaturalHopfAlgebra \App \SetTerms \App \Signature$ and the description
    of the coproduct of $\NaturalHopfAlgebra \App \SetTerms \App \Signature$ provided by
    Proposition~\ref{prop:admissible_sets_coproduct}, straightforward computations lead to
    the fact that $\phi$ is a Hopf algebra morphism from
    \begin{math}
        \NaturalHopfAlgebra \App \SetTerms \App \Signature
        / \Ker \App \Realization_{\LengthsAlphabet \App \Signature}
    \end{math}
    to $\NCK \App \Signature$.
\end{Proof}

Let the linear map
\begin{math}
    \Realization : \NCK \App \Signature \to \K \Angle{\LengthsAlphabet \App \Signature}
\end{math}
defined, for any $\Signature$-trimmed forest $\TrimmedForestF$, by
\begin{math}
    \Realization \App \BasisE_\TrimmedForestF
    :=
    \Realization_{\LengthsAlphabet \App \Signature} \App \BasisE_\ForestF
\end{math}
where $\ForestF$ is a reduced $\Signature$-forest such that $\Trim \App \ForestF =
\TrimmedForestF$. By Lemma~\ref{lem:trimmed_equivalence}, this map is well-defined. We call
$\Realization \App \BasisE_\TrimmedForestF$ the \Def{length polynomial} of $\TrimmedForestF$
on the $\BasisE$-basis. For instance, we have
\begin{equation}
    \Realization
    \App
    \BasisE_{
        \scalebox{.75}{
            \begin{tikzpicture}[Centering,xscale=0.4,yscale=0.5]
                \node[Node,MarkB](2)at(1.25,-1){$\GenC$};
                \node[Node,MarkB](3)at(0.5,-2){$\GenA$};
                \node[Node,MarkB](4)at(2,-2){$\GenC$};
                \node[Node,MarkB](5)at(1,-3){$\GenC$};
                \node[Node,MarkB](6)at(2,-3){$\GenB$};
                \node[Node,MarkB](7)at(3,-3){$\GenB$};
                \node[Node,MarkB](8)at(3.5,-1){$\GenB$};
                \draw[Edge](2)--(3);
                \draw[Edge](2)--(4);
                \draw[Edge](4)--(5);
                \draw[Edge](4)--(6);
                \draw[Edge](4)--(7);
            \end{tikzpicture}
        }
    }
    =
    \sum_{\ell_1, \dots, \ell_7 \in \N}
    \Iverson{\ell_1 < \ell_2}
    \Iverson{\ell_1 < \ell_3}
    \Iverson{\ell_3 < \ell_4}
    \Iverson{\ell_3 < \ell_5}
    \Iverson{\ell_3 < \ell_6} \;
    \Letter^\GenC_{\ell_1} \;
    \Letter^\GenA_{\ell_2} \;
    \Letter^\GenC_{\ell_3} \;
    \Letter^\GenC_{\ell_4} \;
    \Letter^\GenB_{\ell_5} \;
    \Letter^\GenB_{\ell_6} \;
    \Letter^\GenB_{\ell_7}.
\end{equation}

\subsection{Natural Hopf algebras of multiassociative operads}
\label{subsec:multiassociative_operads}
We study here polynomial realizations of natural Hopf algebras of multiassociative operads.
Such operads are parameterized by a signature $\Signature$ and, depending on $\Signature$,
they lead to known Hopf algebras. We obtain in this way polynomial realizations of the Hopf
algebra of noncommutative symmetric functions, of the Hopf algebra of noncommutative
multi-symmetric functions, and of the noncommutative Faà di Bruno Hopf algebra.

\subsubsection{Multiassociative operads} \label{subsubsec:multiassociative_operads}
Let $\Signature$ be a signature. An \Def{$\Signature$-multiset} is a multiset $\Bag{\GenS_1,
\dots, \GenS_\ell}$ of elements of $\Signature$. Given an $\Signature$-term $\TermT$ of
degree $n \geq 0$, the \Def{content} $\Content \App \TermT$ of $\TermT$ is the
$\Signature$-multiset $\Bag{\Decoration_\TermT \App 1, \dots, \Decoration_\TermT \App n}$.
For instance,
\begin{equation}
    \Content \App
    \scalebox{.75}{
        \begin{tikzpicture}[Centering,xscale=0.4,yscale=0.45]
            \node[Leaf](1)at(2.0,0){};
            \node[Node](2)at(2.0,-1){$\GenB$};
            \node[Node](3)at(0.5,-2){$\GenA$};
            \node[Node](4)at(0.5,-3){$\GenB$};
            \node[Leaf](5)at(0,-4){};
            \node[Leaf](6)at(1,-4){};
            \node[Node](7)at(3.5,-2){$\GenC$};
            \node[Node](8)at(2.5,-3){$\GenA$};
            \node[Leaf](9)at(2.5,-4){};
            \node[Leaf](10)at(3.5,-3){};
            \node[Node](11)at(4.5,-3){$\GenB$};
            \node[Leaf](12)at(4,-4){};
            \node[Leaf](13)at(5,-4){};
            \draw[Edge](1)--(2);
            \draw[Edge](2)--(3);
            \draw[Edge](2)--(7);
            \draw[Edge](3)--(4);
            \draw[Edge](4)--(5);
            \draw[Edge](4)--(6);
            \draw[Edge](7)--(8);
            \draw[Edge](7)--(10);
            \draw[Edge](7)--(11);
            \draw[Edge](8)--(9);
            \draw[Edge](11)--(12);
            \draw[Edge](11)--(13);
        \end{tikzpicture}
    }
    = \Bag{\GenA, \GenA, \GenB, \GenB, \GenB, \GenC}.
\end{equation}
Let $\Equiv_{\MAs \App \Signature}$ be the equivalence relation on $\SetTerms \App
\Signature$ satisfying $\TermT \Equiv_{\MAs \App \Signature} \TermT'$ for any
$\Signature$-terms $\TermT$ and $\TermT'$ such that $\Content \App \TermT = \Content \App
\TermT'$.

\begin{Statement}{Proposition}{prop:equivalence_mas}
    For any signature $\Signature$, the equivalence relation $\Equiv_{\MAs \App \Signature}$
    is an operad congruence of $\SetTerms \App \Signature$.
\end{Statement}
\begin{Proof}
    Directly from the definition of $\Equiv_{\MAs \App \Signature}$, for any $\TermT,
    \TermT' \in \SetTerms \App \Signature$, if $\TermT \Equiv_{\MAs \App \Signature}
    \TermT'$, then the number of internal nodes decorated by any $\GenS \in \Signature$ is
    the same in $\TermT$ and $\TermT'$. Therefore, $\Arity \App \TermT = \Arity \App
    \TermT'$. Besides, let $\TermT, \TermT', \TermS \in \SetTerms \App \Signature$ such that
    $\TermT \Equiv_{\MAs \App \Signature} \TermT'$. Hence, $\Content \App \TermT = \Content
    \App \TermT'$, so that, for any $i \in [\Arity \App \TermT]$, from the definition of the
    partial composition of $\SetTerms \App \Signature$, $\Content \App \Par{\TermT \circ_i
    \TermS} = \Content \App \Par{\TermT' \circ_i \TermS}$. For the same reasons, for any $i
    \in [\Arity \App \TermS]$, $\Content \App \Par{\TermS \circ_i \TermT} = \Content \App
    \Par{\TermS \circ_i \TermT'}$. Therefore, we have $\TermT \circ_i \TermS \Equiv_{\MAs
    \App \Signature} \TermT' \circ_i \TermS$ and $\TermS \circ_i \TermT \Equiv_{\MAs \App
    \Signature} \TermS \circ_i \TermT'$. This establishes the statement of the proposition.
\end{Proof}

By Proposition~\ref{prop:equivalence_mas}, the quotient of $\SetTerms \App \Signature$ by
$\Equiv_{\MAs \App \Signature}$ is a well-defined operad, denoted by $\MAs \App \Signature$.
This operad $\MAs \App \Signature$ is a generalization of the $\gamma$-multiassociative
operad, introduced in~\cite{Gir16} (see also ~\cite{Gir16c}), which contains originally only
$\gamma \geq 0$ binary generators. For this reason, we call $\MAs \App \Signature$ the
\Def{$\Signature$-multiassociative operad}. Observe in particular that when the profile of
$\Signature$ is $010^\omega$, $\MAs \App \Signature$ is the free operad on a single
generator of arity $1$ and when the profile of $\Signature$ is $0010^\omega$, $\MAs \App
\Signature$ is the associative operad~$\As$.

The equivalence relation $\Equiv_{\MAs \App \Signature}$ is, directly from its definition,
compatible with the degree. Moreover, since there are finitely many $\Signature$-terms
having a given content, $\Equiv_{\MAs \App \Signature}$ is of finite type. Therefore, by
Proposition~\ref{prop:finitely_factorizable_graded_quotient}, $\MAs \App \Signature$ is
finitely factorizable and graded by the map~$\Deg$.

To describe a combinatorial realization of $\MAs \App \Signature$, let us introduce some
additional definitions about $\Signature$-multisets. Let $\MultisetM := \Bag{\GenS_1, \dots,
\GenS_\ell}$ be an $\Signature$-multiset. The \Def{arity} $\Arity \App \MultisetM$ of
$\MultisetM$ is $\Arity \App \GenS_1 + \dots + \Arity \App \GenS_\ell - \ell + 1$. Note that
the empty $\Signature$-multiset $\emptyset$ has arity $1$. The \Def{degree} $\Deg \App
\MultisetM$ of $\MultisetM$ is $\ell$. The \Def{union} $\MultisetM \cup \MultisetM'$ of two
$\Signature$-multisets $\MultisetM$ and $\MultisetM'$ is the usual union of multisets.

\begin{Statement}{Proposition}{prop:realization_operad_mas}
    For any signature $\Signature$, the operad $\MAs \App \Signature$ admits the following
    combinatorial realization. For any $n \geq 0$, $\MAs \App \Signature \App n$ is the set
    of $\Signature$-multisets of arity $n$. Moreover, for any $\MultisetM, \MultisetM' \in
    \MAs \App \Signature$ and $i \in [\Arity \App \MultisetM]$, $\MultisetM \circ_i
    \MultisetM'$ is the union of $\MultisetM$ and $\MultisetM'$.
\end{Statement}
\begin{Proof}
    First of all, by definition of $\Equiv_{\MAs \App \Signature}$, the set $\Content \App
    \SetTerms \App \Signature$ is a system of representatives of the quotient operad
    $\SetTerms \App \Signature /_{\Equiv_{\MAs \App \Signature}}$. By
    Proposition~\ref{prop:equivalence_mas}, this quotient is well-defined and is the operad
    $\MAs \App \Signature$. Moreover, observe that for any $\TermT \in \SetTerms \App
    \Signature$, by definition of the arity of $\Signature$-multisets, we have $\Arity \App
    \TermT = \Arity \App \Content \App \TermT$. Therefore, for any $n \geq 0$, $\MAs \App
    \Signature \App n$ can be identified with the set of $\Signature$-multisets of arity
    $n$. Finally, since for any $\TermT, \TermT' \in \SetTerms \App \Signature$ and $i \in
    [\Arity \App \TermT]$, we have
    \begin{math}
        \Content \App \Par{\TermT \circ_i \TermT'}
        = \Content \App \TermT \cup \Content \App \TermT',
    \end{math}
    the rule for the partial composition of $\MAs \App \Signature$ given in the
    statement of the proposition follows.
\end{Proof}

In the sequel, we shall identify $\MAs \App \Signature$ with its combinatorial realization
provided by Proposition~\ref{prop:realization_operad_mas}.

By Proposition~\ref{prop:realization_operad_mas}, the $\SignatureExample$-multiset
$\Bag{\GenA, \GenA, \GenB, \GenC, \GenC, \GenC}$ is an element of arity $8$ of the operad
$\MAs \App \SignatureExample$. Moreover, in this operad, we have the partial composition
\begin{equation}
    \Bag{\GenA, \GenB, \GenB, \GenB, \GenC}
    \circ_4
    \Bag{\ColB{\GenB}, \ColB{\GenC}, \ColB{\GenC}}
    = \Bag{\GenA, \GenB, \GenB, \GenB, \ColB{\GenB}, \GenC, \ColB{\GenC}, \ColB{\GenC}}.
\end{equation}

\subsubsection{Natural Hopf algebras} \label{subsubsec:n_mas}
Since by Proposition~\ref{prop:finitely_factorizable_graded_quotient}, $\MAs \App
\Signature$ is finitely factorizable and graded by the map $\Deg$, $\NaturalHopfAlgebra \App
\MAs \App \Signature$ is a well-defined Hopf algebra. By construction of
$\NaturalHopfAlgebra \App \MAs \App \Signature$, the bases of this Hopf algebra are indexed
by reduced words on $\MAs \App \Signature$. Moreover, by construction, the coproduct of
$\NaturalHopfAlgebra \App \MAs \App \Signature$ satisfies, for any nonempty
$\Signature$-multiset $\MultisetM \in \MAs \App \Signature$,
\begin{equation} \label{equ:coproduct_n_mas}
    \Coproduct \App \BasisE_\MultisetM
    =
    \sum_{
        \substack{
            \MultisetM' \in \MAs \App \Signature \\
            u \in \MAs \App \Signature^*
        }
    }
    \Iverson{\Arity \App \MultisetM' = \Length \App u} \;
    \Iverson{
        \MultisetM =
        \MultisetM' \cup u \App 1 \cup \dots \cup u \App \Length \App u
    }
    \;
    \BasisE_{\Reduced \App \MultisetM'} \otimes \BasisE_{\Reduced \App u}.
\end{equation}
For instance, in $\NaturalHopfAlgebra \App \MAs \App \SignatureExample$, we have
\begin{align}
    \Coproduct \App \BasisE_{\Bag{\GenA, \GenB, \GenB}}
    & =
    \BasisE_{\emptyset} \otimes \BasisE_{\Bag{\GenA, \GenB, \GenB}}
    + \BasisE_{\Bag{\GenA}} \otimes \BasisE_{\Bag{\GenB, \GenB}}
    + 2 \BasisE_{\Bag{\GenB}} \otimes \BasisE_{\Bag{\GenA, \GenB}}
    + \BasisE_{\Bag{\GenB}} \otimes \BasisE_{\Bag{\GenA} \Bag{\GenB}}
    \\
    &
    \quad
    + \BasisE_{\Bag{\GenB}} \otimes \BasisE_{\Bag{\GenB} \Bag{\GenA}}
    + 2 \BasisE_{\Bag{\GenA, \GenB}} \otimes \BasisE_{\Bag{\GenB}}
    + 3 \BasisE_{\Bag{\GenB, \GenB}} \otimes \BasisE_{\Bag{\GenA}}
    + \BasisE_{\Bag{\GenA, \GenB, \GenB}} \otimes \BasisE_{\emptyset}.
    \notag
\end{align}

\subsubsection{Multi-symmetric functions and Faà di Bruno Hopf algebras}
\label{subsubsec:n_mas_fdb}
Let $\Signature$ be a signature of profile $0 0^r s 0^\omega$ where $r, s \in \N$ and let
$\NCFdB_r^{(s)} := \NaturalHopfAlgebra \App \MAs \App \Signature$. Let us consider the
following particular cases.

\begin{enumerate}[label=({\sf \arabic*})]
    \item When $r = 0$, $\NCFdB_r^{(s)}$ is the Hopf algebra $\NCSF^{(s)}$ of noncommutative
    multi-symmetric functions of level $s$~\cite{NT10}. In particular, $\NCSF^{(1)}$ is the
    noncommutative symmetric functions Hopf algebra $\NCSF$~\cite{GKLLRT95}. Let us denote
    by $\GenA_1$, \dots, $\GenA_s$ the elements of $\Signature$ and by the word $u$ of
    length $s$ on $\N$ the $\Signature$-multiset having $u \App i$ occurrences of $\GenA_i$
    for any $i \in [s]$. By~\eqref{equ:coproduct_n_mas}, for any $u \in \N^s$, we have
    \begin{equation}
        \Coproduct \App \BasisE_u
        = \sum_{u_1, u_2 \in \N^s}
        \Iverson{u_1 \App i + u_2 \App i = u \App i \text{ for all } i \in [s]} \;
        \BasisE_{u_1} \otimes \BasisE_{u_2}.
    \end{equation}
    For instance, for $s = 3$,
    \begin{equation}
        \Coproduct \App \BasisE_{120}
        = \BasisE_{000} \otimes \BasisE_{120}
        + \BasisE_{010} \otimes \BasisE_{110}
        + \BasisE_{020} \otimes \BasisE_{100}
        + \BasisE_{100} \otimes \BasisE_{020}
        + \BasisE_{110} \otimes \BasisE_{010}
        + \BasisE_{020} \otimes \BasisE_{000}.
    \end{equation}
    \item When $s = 0$, $\NCFdB_r^{(s)}$ is the $r$-deformation $\NCFdB_r$ of
    $\NCFdB$~\cite{Foi08} (see also~\cite{Bul11}). This is a consequence of a result
    of~\cite{BG16}, which provides a generalization of the construction
    $\NaturalHopfAlgebra$, but taking at input pros~\cite{McL65} instead of operads. For
    some pros and some operads, the two constructions coincide.  Let us denote by $\GenA$
    the unique element of $\Signature$ and by the integer $d$ the $\Signature$-multiset made
    of $d$ occurrences of $\GenA$. Observe that $\Arity \App d = dr + 1$.
    By~\eqref{equ:coproduct_n_mas}, for any $d \in \N$, we have
    \begin{equation}
        \Coproduct \App \BasisE_d
        =
        \sum_{d' \in \HanL{d}}
        \sum_{
            \substack{
                \ell \geq 0 \\
                d'_1, \dots, d'_\ell \geq 1
            }
        }
        \Iverson{d' + d'_1 + \dots + d'_\ell = d} \;
        \binom{d' r + 1}{\ell} \;
        \BasisE_{d'} \otimes \BasisE_{d'_1 \dots d'_\ell}.
    \end{equation}
    For instance, for $s = 2$, we have
    \begin{equation}
        \Coproduct \App \BasisE_3
        =
        \BasisE_\epsilon \otimes \BasisE_3
        + 3 \BasisE_1 \otimes \BasisE_2
        + 3 \BasisE_1 \otimes \BasisE_{11}
        + 5 \BasisE_2 \otimes \BasisE_1
        + \BasisE_3 \otimes \BasisE_\epsilon.
    \end{equation}
    Observe that the Hopf algebra $\NCFdB_1$ is the noncommutative Faà di Bruno Hopf algebra
    $\NCFdB$~\cite{FG05,BFK06,Foi08}. The construction of this Hopf algebra is detailed in
    Section~\ref{subsubsec:faa_di_bruno_hopf_algebra}.
\end{enumerate}
From these two particular cases, we call $\NCFdB_r^{(s)}$ the \Def{noncommutative multi-Faà
di Bruno Hopf algebra}. This Hopf algebra $\NCFdB_r^{(s)}$ is to $\NCFdB_r$ what
$\NCSF^{(s)}$ is to~$\NCSF$.

\subsubsection{Polynomial realization} \label{subsubsec:polynomial_realization_n_mas}
Let $A$ be an $\Signature$-forest-like alphabet and $\MultisetM \in \MAs \App \Signature$ be
a nonempty $\Signature$-multiset. By using Theorem~\ref{thm:quotient_operad_hopf_morphism},
we obtain that the $A$-realization of $\MultisetM$ on the $\BasisE$-basis of
$\NaturalHopfAlgebra \App \MAs \App \Signature$ satisfies
\begin{equation} \label{equ:realization_faa_du_bruno}
    \bar{\Realization}_A \App \BasisE_{\MultisetM}
    = \sum_{\TermT \in \SetTerms \App \Signature}
    \Iverson{\Content \App \TermT = \MultisetM} \;
    \Realization_A \App \BasisE_\TermT.
\end{equation}
By Proposition~\ref{prop:realization_quotient_operads},
$\bar{\Realization}_{\PositionsAlphabet \App \Signature}$ is a polynomial realization of
$\NaturalHopfAlgebra \App \MAs \App \Signature$. On the other hand, the map
$\bar{\Realization}_{\LengthsAlphabet \App \Signature}$ exhibits an interesting property,
as stated in the following result.

\begin{Statement}{Theorem}{thm:n_mas_length_realization}
    For any signature $\Signature$ and any nonempty $\Signature$-multiset $\MultisetM \in
    \MAs \App \Signature$, the map~\begin{math}
        \bar{\Realization}_{\LengthsAlphabet \App \Signature} :
        \NaturalHopfAlgebra \App \MAs \App \Signature
        \to \K \Angle{\LengthsAlphabet \App \Signature}
    \end{math}
    satisfies
    \begin{equation} \label{equ:n_mas_length_realization}
        \bar{\Realization}_{\LengthsAlphabet \App \Signature}
        \App \BasisE_{\MultisetM}
        =
        \sum_{\TrimmedTermT \in \Trim \App \SetTerms \App \Signature}
        \Iverson{\Content \App \TrimmedTermT = \MultisetM} \;
        \Charge \App \TrimmedTermT \;
        \Realization \App \BasisE_\TrimmedTermT.
    \end{equation}
    Moreover, this map $\bar{\Realization}_{\LengthsAlphabet \App \Signature}$ is injective.
\end{Statement}
\begin{Proof}
    By using successively~\eqref{equ:realization_faa_du_bruno},
    Lemma~\ref{lem:number_trimmed}, Proposition~\ref{prop:projection_lengths}, and the
    notion of length polynomial of $\Signature$-trimmed forests introduced at the end of
    Section~\ref{subsubsec:quotient_nck}, we have
    \begin{align}
        \bar{\Realization}_{\LengthsAlphabet \App \Signature} \App
        \BasisE_\MultisetM
        &
        = \sum_{\TermT \in \SetTerms \App \Signature}
        \Iverson{\Content \App \TermT = \MultisetM} \;
        \Realization_{\LengthsAlphabet \App \Signature} \App \BasisE_\TermT
        \\
        &
        = \sum_{\TrimmedTermT \in \Trim \App \SetTerms \App \Signature}
        \Iverson{\Content \App \TrimmedTermT = \MultisetM} \;
        \sum_{\TermT \in \SetTerms \App \Signature}
        \Iverson{\Trim \App \TermT = \TrimmedTermT} \;
        \Realization_{\LengthsAlphabet \App \Signature} \App \BasisE_\TermT
        \notag
        \\
        &
        = \sum_{\TrimmedTermT \in \Trim \App \SetTerms \App \Signature}
        \Iverson{\Content \App \TrimmedTermT = \MultisetM} \;
        \Charge \App \TrimmedTermT \;
        \Realization \App \BasisE_\TrimmedTermT.
        \notag
    \end{align}
    This establishes the first part of the statement.

    Let us prove the injectivity of $\bar{\Realization}_{\LengthsAlphabet \App \Signature}$.
    First of all, observe from~\eqref{equ:n_mas_length_realization} that for any nonempty
    $\Signature$-multiset $\MultisetM$ of $\MAs \App \Signature$ of degree $n \geq 1$, all
    monomials $w$ appearing in $\bar{\Realization}_{\LengthsAlphabet \App \Signature} \App
    \BasisE_\MultisetM$ are such that $w = \Letter^{\GenS_1}_{\ell_1} \dots
    \Letter^{\GenS_n}_{\ell_n}$ where $\ell_1, \dots, \ell_n \in \N$, $\GenS_1, \dots,
    \GenS_n \in \Signature$, and $\MultisetM = \Bag{\GenS_1, \dots, \GenS_n}$. Therefore,
    for any nonempty $\Signature$-multisets $\MultisetM_1, \MultisetM_2$ of $\MAs \App
        \Signature$, 
    \begin{math}
        \bar{\Realization}_{\LengthsAlphabet \App \Signature} \App \BasisE_{\MultisetM_1}
        =
        \bar{\Realization}_{\LengthsAlphabet \App \Signature} \App \BasisE_{\MultisetM_2}
    \end{math}
    implies $\MultisetM_1 = \MultisetM_2$. Now,
    Expression~\eqref{equ:n_mas_length_realization} together with the fact that, by
    Theorem~\ref{thm:quotient_operad_hopf_morphism}, $\bar{\Realization}_{\LengthsAlphabet
    \App \Signature}$ is an algebra morphism, lead to the following property. For any $x \in
    \Reduced \App \Par{\MAs \App \Signature^*}$, there exist $w \in \LengthsAlphabet \App
    \Signature^*$ and $i \in \Han{\Length \App w}$ such that $w$ appears in
    $\bar{\Realization}_{\LengthsAlphabet \App \Signature} \App \BasisE_x$ and $w \App i =
    \Letter^\GenS_0$, $\GenS \in \Signature$ if and only if the word $x$ decomposes as $x =
    x_1 \Conc x_2$ for some $x_1, x_2 \in \Reduced \App \Par{\MAs \App \Signature^*}$ such
    that $i = \Deg \App x_1 + 1$. These two properties imply together that for any $x_1, x_2
    \in \Reduced \App \Par{\MAs \App \Signature^*}$, if 
    \begin{math}
        \bar{\Realization}_{\LengthsAlphabet \App \Signature} \App \BasisE_{x_1}
        =
        \bar{\Realization}_{\LengthsAlphabet \App \Signature} \App \BasisE_{x_2},
    \end{math}
    then there are some nonempty $\Signature$-multisets $\MultisetM_1$, \dots,
    $\MultisetM_n$, $n \geq 0$, of $\MAs \App \Signature$ such that $x_1 = \MultisetM_1
    \Conc \dots \Conc \MultisetM_n = x_2$. Hence, $\bar{\Realization}_{\LengthsAlphabet \App
    \Signature}$ is injective.
\end{Proof}

By Theorem~\ref{thm:n_mas_length_realization}, $\bar{\Realization}_{\LengthsAlphabet \App
\Signature}$ is an additional polynomial realization of $\NaturalHopfAlgebra \App \MAs \App
\Signature$. This realization is simpler than the previous one since it uses the
$\Signature$-forest-like alphabet $\LengthsAlphabet \App \Signature$ which is a quotient
of~$\PositionsAlphabet \App \Signature$.

For instance, in $\NaturalHopfAlgebra \App \MAs \App \SignatureExample$, we have
\begin{align}
    \bar{\Realization}_{\LengthsAlphabet \App \Signature}
    \App \BasisE_{\Bag{\GenA, \GenB, \GenB}}
    & =
    \sum_{\ell_1, \ell_2, \ell_3 \in \N}
    \Iverson{\ell_1 < \ell_2 < \ell_3} \;
    2 \Letter^\GenA_{\ell_1} \; \Letter^\GenB_{\ell_2} \; \Letter^\GenB_{\ell_3}
    \\
    &
    \quad +
    \sum_{\ell_1, \ell_2, \ell_3 \in \N}
    \Iverson{\ell_1 < \ell_2 < \ell_3} \;
    2 \Letter^\GenB_{\ell_1} \; \Letter^\GenA_{\ell_2} \; \Letter^\GenB_{\ell_3}
    +
    \sum_{\ell_1, \ell_2, \ell_3 \in \N}
    \Iverson{\ell_1 < \ell_2} \Iverson{\ell_1 < \ell_3} \;
    \Letter^\GenB_{\ell_1} \; \Letter^\GenA_{\ell_2} \; \Letter^\GenB_{\ell_3}
    \notag
    \\
    &
    \quad +
    \sum_{\ell_1, \ell_2, \ell_3 \in \N}
    \Iverson{\ell_1 < \ell_2 < \ell_3} \;
    4 \Letter^\GenB_{\ell_1} \; \Letter^\GenB_{\ell_2} \; \Letter^\GenA_{\ell_3}
    +
    \sum_{\ell_1, \ell_2, \ell_3 \in \N}
    \Iverson{\ell_1 < \ell_2} \Iverson{\ell_1 < \ell_3} \;
    \Letter^\GenB_{\ell_1} \; \Letter^\GenB_{\ell_2} \; \Letter^\GenA_{\ell_3}
    \notag
    \\
    & =
    \sum_{\ell_1, \ell_2, \ell_3 \in \N}
    \Iverson{\ell_1 < \ell_2 < \ell_3} \;
    2 \Letter^\GenA_{\ell_1} \; \Letter^\GenB_{\ell_2} \; \Letter^\GenB_{\ell_3}
    \notag
    \\
    &
    \quad +
    \sum_{\ell_1, \ell_2, \ell_3 \in \N}
    \Par{
        3 \; \Iverson{\ell_1 < \ell_2 < \ell_3} \;
        +
        \Iverson{\ell_1 < \ell_3 \leq \ell_2}
    } \;
    \Letter^\GenB_{\ell_1} \; \Letter^\GenA_{\ell_2} \; \Letter^\GenB_{\ell_3}
    \notag
    \\
    &
    \quad +
    \sum_{\ell_1, \ell_2, \ell_3 \in \N}
    \Par{
        5 \; \Iverson{\ell_1 < \ell_2 < \ell_3} \;
        +
        \Iverson{\ell_1 < \ell_3 \leq \ell_2}
    } \;
    \Letter^\GenB_{\ell_1} \; \Letter^\GenB_{\ell_2} \; \Letter^\GenA_{\ell_3}.
    \notag
\end{align}

Remark that by Theorem~\ref{thm:n_mas_length_realization} and by using the construction of
the noncommutative multi-Faà du Bruno Hopf algebra $\NCFdB_r^{(s)}$ presented in
Section~\ref{subsubsec:n_mas_fdb}, $\bar{\Realization}_{\LengthsAlphabet \App \Signature}$
is a polynomial realization of this Hopf algebra and also of $\NCSF^{(s)}$ and $\NCFdB_r$.

\subsection{Natural Hopf algebras of interstice operads}
\label{subsec:interstice_operads}
We study finally here polynomial realizations of natural Hopf algebras of interstice
operads. These Hopf algebras are known as double tensor Hopf algebras.

\subsubsection{Interstice operads}
Let $\Signature$ be a binary signature. An \Def{$\Signature$-word} is a word on
$\Signature$. Given an $\Signature$-term $\TermT$ of degree $n \geq 0$, the \Def{infix
reading} $\Infix \App \TermT$ of $\TermT$ is the $\Signature$-word $u$ of length $n$
obtained by reading the decorations of the internal nodes of $\TermT$ according to its left
to right infix traversal. Note that this traversal is well-defined since, because
$\Signature$ is binary, each internal node of $\TermT$ has exactly two children. For
instance, if $\Signature$ is the binary signature containing $\GenA$, $\GenB$, and $\GenC$,
\begin{equation}
    \Infix \App
    \scalebox{.75}{
        \begin{tikzpicture}[Centering,xscale=0.38,yscale=0.45]
            \node[Leaf](1)at(3,0){};
            \node[Node](2)at(3,-1){$\GenB$};
            \node[Node](3)at(1,-2){$\GenA$};
            \node[Leaf](4)at(0,-3){};
            \node[Node](5)at(1.5,-3){$\GenB$};
            \node[Leaf](6)at(1,-4){};
            \node[Leaf](7)at(2,-4){};
            \node[Node](8)at(4.5,-2){$\GenC$};
            \node[Node](9)at(3.5,-3){$\GenC$};
            \node[Leaf](10)at(3,-4){};
            \node[Leaf](11)at(4,-4){};
            \node[Node](12)at(5.5,-3){$\GenA$};
            \node[Leaf](13)at(5,-4){};
            \node[Leaf](14)at(6,-4){};
            \draw[Edge](1)--(2);
            \draw[Edge](2)--(3);
            \draw[Edge](2)--(8);
            \draw[Edge](3)--(4);
            \draw[Edge](3)--(5);
            \draw[Edge](5)--(6);
            \draw[Edge](5)--(7);
            \draw[Edge](8)--(9);
            \draw[Edge](8)--(12);
            \draw[Edge](9)--(10);
            \draw[Edge](9)--(11);
            \draw[Edge](12)--(13);
            \draw[Edge](12)--(14);
        \end{tikzpicture}
    }
    = \GenA \GenB \GenB \GenC \GenC \GenA.
\end{equation}
Let $\Equiv_{\Int \App \Signature}$ be the equivalence relation on $\SetTerms \App
\Signature$ satisfying $\TermT \Equiv_{\Int \App \Signature} \TermT'$ for any
$\Signature$-terms $\TermT$ and $\TermT'$ such that $\Infix \App \TermT = \Infix \App
\TermT'$.

\begin{Statement}{Proposition}{prop:equivalence_int}
    For any binary signature $\Signature$, the equivalence relation $\Equiv_{\Int \App
    \Signature}$ is an operad congruence of $\SetTerms \App \Signature$.
\end{Statement}
\begin{Proof}
    Directly from the definition of $\Equiv_{\Int \App \Signature}$, for any $\TermT,
    \TermT' \in \SetTerms \App \Signature$, if $\TermT \Equiv_{\Int \App \Signature}
    \TermT'$, then the number of internal nodes decorated by any $\GenS \in \Signature$ is
    the same in $\TermT$ and $\TermT'$. Therefore, $\Arity \App \TermT = \Arity \App
    \TermT'$. Besides, let $\TermT, \TermT', \TermS \in \SetTerms \App \Signature$ such that
    $\TermT \Equiv_{\Int \App \Signature} \TermT'$. Hence, $\Infix \App \TermT = \Infix \App
    \TermT'$, so that, for any $i \in [\Arity \App \TermT]$, from the definition of the
    partial composition of $\SetTerms \App \Signature$, $\Infix \App \Par{\TermT \circ_i
    \TermS} = \Infix \App \Par{\TermT' \circ_i \TermS}$. For the same reasons, for any $i
    \in [\Arity \App \TermS]$, $\Infix \App \Par{\TermS \circ_i \TermT} = \Infix \App
    \Par{\TermS \circ_i \TermT'}$. Therefore, we have $\TermT \circ_i \TermS \Equiv_{\Int
    \App \Signature} \TermT' \circ_i \TermS$ and $\TermS \circ_i \TermT \Equiv_{\Int \App
    \Signature} \TermS \circ_i \TermT'$. This establishes the statement of the proposition.
\end{Proof}

By Proposition~\ref{prop:equivalence_int}, the quotient of $\SetTerms \App \Signature$ by
$\Equiv_{\Int \App \Signature}$ is a well-defined operad, denoted by $\Int \App \Signature$.
This operad $\Int \App \Signature$ is in fact the interstice operad introduced
in~\cite{CG22}. Observe in particular that when the profile of $\Signature$ is
$0010^\omega$, $\Int \App \Signature$ is the associative operad $\As$ and when the profile
of $\Signature$ is $0020^\omega$, $\Int \App \Signature$ is the operad whose algebras are
equipped with two associative and mutually associative operations~\cite{Pir03}.

The equivalence relation $\Equiv_{\Int \App \Signature}$ is, directly from its definition,
compatible with the degree. Moreover, since there are finitely many $\Signature$-terms
having a given infix reading, $\Equiv_{\Int \App \Signature}$ is of finite type. Therefore,
by Proposition~\ref{prop:finitely_factorizable_graded_quotient}, $\Int \App \Signature$ is
finitely factorizable and graded by the map~$\Deg$. Observe that since for any
$\Signature$-terms $\TermT$ and $\TermT'$, $\Infix \App \TermT = \Infix \TermT'$ implies
$\Content \App \TermT = \Content \App \TermT'$, the operad $\MAs \App \Signature$ is a
quotient of~$\Int \App \Signature$.

To describe a combinatorial realization of $\Int \App \Signature$, let us introduce some
additional definitions about $\Signature$-words. Let $u$ be an $\Signature$-word. The
\Def{arity} $\Arity \App u$ of $u$ is $\Length \App u + 1$. Note that the empty
$\Signature$-word $\epsilon$ is the unique $\Signature$-word of arity $1$ and that there is
no $\Signature$-word of arity $0$. The \Def{degree} $\Deg \App u$ of $u$ is $\Length \App
u$. Given $i_1, i_2 \in \Han{\Length \App u}$ such that $i_1 \leq i_2$, let $u_{\Par{i_1,
i_2}}$ be the factor $u \App i_1 \dots u \App i_2$ of~$u$.

The operad $\Int \App \Signature$ admits the following combinatorial realization described
in~\cite{CG22}. For any $n \geq 0$, $\Int \App \Signature \App n$ is the set of
$\Signature$-words on $\Signature$ of arity $n$. Moreover, for any $u, u' \in \Int \App
\Signature$ and $i \in [\Arity \App u]$, $u \circ_i u'$ is the $\Signature$-word $u_{(1, i -
1)} \Conc u' \Conc u_{(i, \Length \App u)}$. In the sequel, we shall identify $\Int \App
\Signature$ with this combinatorial realization. For instance, if $\Signature$ is the binary
signature containing $\GenA$ and $\GenB$,
\begin{equation}
    \GenB \GenB \GenA \GenB \GenA \circ_3 \ColB{\GenA \GenA \GenB} =
    \GenB \GenB
    \ColB{\GenA \GenA \GenB}
    \GenA \GenB \GenA.
\end{equation}

\subsubsection{Hopf algebra of phrases} \label{subsubsec:hopf_algebra_phrases}
Let $\Signature$ be a binary signature. Since by
Proposition~\ref{prop:finitely_factorizable_graded_quotient}, $\Int \App \Signature$ is
finitely factorizable and graded by the map $\Deg$, $\Phr \App \Signature :=
\NaturalHopfAlgebra \App \Int \App \Signature$ is a well-defined Hopf algebra. By
construction, the bases of $\Phr \App \Signature$ are indexed by reduced words on $\Int \App
\Signature$. Following the terminology of~\cite{Man97}, such elements are called
\Def{$\Signature$-phrases}. An $\Signature$-phrase is denoted by separating the
$\Signature$-words forming it by commas. For instance, if $\Signature$ is the binary
signature containing $\GenA$ and $\GenB$, then $\GenA \GenA, \GenB \GenA \GenB, \GenB,
\GenB$ is an $\Signature$-phrase made of the $\Signature$-words $\GenA \GenA$, $\GenB \GenA
\GenB$, $\GenB$, and $\GenB$. Moreover, by construction, the coproduct of $\Phr \App
\Signature$ satisfies, for any nonempty $\Signature$-word $u$,
\begin{equation}
    \Coproduct \App \BasisE_u =
    \sum_{
        \substack{
            v \in \Int \App \Signature \\
            w_1, \dots, w_{\Length \App v + 1} \in \Int \App \Signature
        }
    }
    \Iverson{
        w_1 \Conc v \App 1
        \Conc \dots \Conc
        w_{\Length \App v} \Conc v \App \Length \App v
        \Conc w_{\Length \App v + 1}
        = u
    }
    \;
    \BasisE_{\Reduced \App v}
    \otimes
    \BasisE_{\Reduced \App w_1 \dots w_{\Length \App v + 1}}.
\end{equation}
For instance, in $\Phr \App \Signature$ where $\Signature$ is the binary signature of the
previous example, we have
\begin{equation}
    \Coproduct \App \BasisE_{\GenA \GenA \GenB}
    =
    \BasisE_\epsilon \otimes \BasisE_{\GenA \GenA \GenB}
    + \BasisE_{\GenA} \otimes \BasisE_{\GenA, \GenB}
    + \BasisE_{\GenA} \otimes \BasisE_{\GenA \GenB}
    + \BasisE_{\GenB} \otimes \BasisE_{\GenA \GenA}
    + \BasisE_{\GenA \GenA} \otimes \BasisE_{\GenB}
    + 2 \BasisE_{\GenA \GenB} \otimes \BasisE_{\GenA}
    + \BasisE_{\GenA \GenA \GenB} \otimes \BasisE_\epsilon.
\end{equation}
As noticed in~\cite{CV23}, this Hopf algebra $\Phr \App \Signature$ is isomorphic to the
double tensor Hopf algebra built in~\cite{EP15}.

\subsubsection{Polynomial realization} \label{subsubsec:realization_phrase_hopf_algebra}
Let $\Signature$ be a binary signature, $A$ be an $\Signature$-forest-like alphabet, and $u$
be a nonempty $\Signature$-word. By using Theorem~\ref{thm:quotient_operad_hopf_morphism},
we obtain that the $A$-realization of $u$ on the $\BasisE$-basis of $\Phr \App \Signature$
satisfies
\begin{equation}
    \bar{\Realization}_A \App \BasisE_u
    =
    \sum_{\TermT \in \SetTerms \App \Signature}
    \Iverson{\Infix \App \TermT = u} \;
    \Realization_A \App \BasisE_\TermT.
\end{equation}
By Proposition~\ref{prop:realization_quotient_operads},
$\bar{\Realization}_{\PositionsAlphabet \App \Signature}$ is a polynomial realization of
$\Phr \App \Signature$.

For instance, in $\Phr \App \Signature$ where $\Signature$ is the binary signature of the
examples of Section~\ref{subsubsec:hopf_algebra_phrases}, we have
\begin{align}
    \bar{\Realization}_{\PositionsAlphabet \App \Signature}
    \App \BasisE_{\GenA \GenA \GenB}
    & =
    \sum_{\ell_1, \ell_2, \ell_3 \in \N}
    \Letter^\GenB_{0^{\ell_1}} \;
    \Letter^\GenA_{0^{\ell_1} 1 0^{\ell_2}} \;
    \Letter^\GenA_{0^{\ell_1} 1 0^{\ell_2} 1 0^{\ell_3}}
    +
    \sum_{\ell_1, \ell_2, \ell_3 \in \N}
    \Letter^\GenB_{0^{\ell_1}} \;
    \Letter^\GenA_{0^{\ell_1} 1 0^{\ell_2}} \;
    \Letter^\GenA_{0^{\ell_1} 1 0^{\ell_2} 2 0^{\ell_3}}
    \notag
    \\
    & \quad
    +
    \sum_{\ell_1, \ell_2, \ell_3 \in \N}
    \Letter^\GenA_{0^{\ell_1}} \;
    \Letter^\GenA_{0^{\ell_1} 1 0^{\ell_2}} \;
    \Letter^\GenB_{0^{\ell_1} 2 0^{\ell_3}}
    +
    \sum_{\ell_1, \ell_2, \ell_3 \in \N}
    \Letter^\GenA_{0^{\ell_1}} \;
    \Letter^\GenB_{0^{\ell_1} 2 0^{\ell_2}} \;
    \Letter^\GenA_{0^{\ell_1} 2 0^{\ell_2} 1 0^{\ell_3}}
    \notag
    \\
    & \quad
    +
    \sum_{\ell_1, \ell_2, \ell_3 \in \N}
    \Letter^\GenA_{0^{\ell_1}} \;
    \Letter^\GenA_{0^{\ell_1} 2 0^{\ell_2}} \;
    \Letter^\GenB_{0^{\ell_1} 2 0^{\ell_2} 2 0^{\ell_3}}.
\end{align}
Observe that the map $\bar{\Realization}_{\LengthsAlphabet \App \Signature}$ is not
injective. Indeed, we have for instance
\begin{equation}
    \bar{\Realization}_{\LengthsAlphabet \App \Signature}
    \App \BasisE_{\GenA \GenB}
    =
    \sum_{\ell_1, \ell_2 \in \N}
    \Iverson{\ell_1 < \ell_2} \;
    \Par{
        \Letter^\GenB_{\ell_1} \; \Letter^\GenA_{\ell_2}
        +
        \Letter^\GenA_{\ell_1} \; \Letter^\GenB_{\ell_2}
    }
    =
    \bar{\Realization}_{\LengthsAlphabet \App \Signature}
    \App \BasisE_{\GenB \GenA}.
\end{equation}
For this reason, unlike the case of the Hopf algebra $\NaturalHopfAlgebra \App \MAs \App
\Signature$ presented in Section~\ref{subsubsec:polynomial_realization_n_mas},
$\bar{\Realization}_{\LengthsAlphabet \App \Signature}$ is not a polynomial realization of
$\Phr \App \Signature$.

\section{Conclusion and future work}
We have introduced a polynomial realization of the natural Hopf algebra $\NaturalHopfAlgebra
\App \SetTerms \App \Signature$ of a free operad $\SetTerms \App \Signature$
(Theorem~\ref{thm:realization}) and, additionally, of the natural Hopf algebra
$\NaturalHopfAlgebra \App \Operad$ of an operad $\Operad$ in the case where $\Operad$ can be
described as a quotient of a free operad satisfying certain properties
(Proposition~\ref{prop:realization_quotient_operads}). At the heart of these realizations
lies the notion of the position of internal nodes of a forest. Another important tool is
that of related alphabets, which provides a framework for working with polynomial
realizations. Although this has not been developed in this work, related alphabets allow for
a unified treatment of already known polynomial realizations. Applications of this
polynomial realization of $\NaturalHopfAlgebra \App \SetTerms \App \Signature$ are proposed.
We have seen for instance that $\NaturalHopfAlgebra \App \SetTerms \App \Signature$ can be
sent to a space of a decorated version of word quasi-symmetric functions
(Theorem~\ref{thm:map_to_wqsym}) and that it contains a Hopf subalgebra of a decorated
version of the noncommutative Connes-Kreimer Hopf algebra (Theorem~\ref{thm:map_to_nck}). As
another consequence, we have also provided polynomial realizations of the noncommutative Faà
di Bruno Hopf algebra (Theorem~\ref{thm:n_mas_length_realization}) and of the double tensor
Hopf algebra (Section~\ref{subsubsec:realization_phrase_hopf_algebra}). Here are some open
questions in this context as well as some avenues for future research.

In~\cite{Gir15} and~\cite{Gir16b,Gir16}, a family of operads based on various families of
combinatorial objects is constructed. This family includes the operad $\FCat_m$ involving
$m$-Fuss-Catalan objects, the operad $\Schr$ involving Schröder trees, the operad $\Motz$
involving Motzkin paths, the operad $\Comp$ involving integer compositions, the operad $\DA$
involving directed animals, the operad $\SComp$ involving segmented integer compositions,
the pluriassociative operad $\Dias_\gamma$ involving some words of integers, and the
polydendriform $\Dendr_\gamma$ involving binary trees with decorated edges. These operads,
when considered as quotients of free operads by an operad congruence, satisfy the conditions
listed in Section~\ref{subsec:quotients_hopf_subalgebras}. Therefore, their natural Hopf
algebras can be seen as natural Hopf subalgebras of natural Hopf algebras of free operads
(see Theorem~\ref{thm:quotient_operad_hopf_morphism}). The question here consists in
applying the results of this work to obtain polynomial realizations of these Hopf algebras.
As a consequence, we can hope to obtain new families of polynomials, generalizing symmetric
functions.

A second question is the following. We have established the fact that a family of
polynomials on the alphabet $\PositionsAlphabet \App \Signature$ of positions provides a
polynomial realization of natural Hopf algebras of operads. However, a similar family of
polynomials on the alphabet $\LengthsAlphabet \App \Signature$ of lengths admits some
interesting properties. For instance, the image of the map $\Realization_{\LengthsAlphabet
\App \Signature}$ is linked with a decorated version of word quasi-symmetric functions (see
Section~\ref{subsec:decorated_word_quasi_symmetric_functions}) and with a decorated version
of the noncommutative Connes-Kreimer Hopf algebra (see Section~\ref{subsec:connes_kreimer}).
Moreover, this map remains a polynomial realization of natural Hopf algebras of
multiassociative operads, including the noncommutative Faà di Bruno Hopf algebra and the
Hopf algebra of multi-symmetric functions (see
Section~\ref{subsec:multiassociative_operads}). In contrast, this map is not injective in
the case of the natural Hopf algebra of interstice operads (see
Section~\ref{subsec:interstice_operads}). The question here is first to describe the kernel
of this map in the previous particular case. Next, one question is to look for necessary and
sufficient conditions on the operad $\Operad$ to ensure that the map
$\Realization_{\LengthsAlphabet \App \Signature}$ is a polynomial realization
of~$\NaturalHopfAlgebra \App \Operad$.

A last research focus addressed here involves the suitable definition of a Cartesian product
$\AlphabetProduct$ on the class of $\Signature$-forest-like alphabets, leading to the
definition of an internal coproduct on $\NaturalHopfAlgebra \App \SetTerms \App \Signature$.
When this alphabet product operation $\AlphabetProduct$ is associative, this would endow
$\NaturalHopfAlgebra \App \SetTerms \App \Signature$ with a different coalgebra structure
(see for instance~\cite{NT10,FNT14,Foi20} for examples of such constructions). The main
difficulty here is to propose a coherent way to define the root relation, the
$\GenS$-decoration relations, $\GenS \in \Signature$, and the $j$-edge relations, $j \geq
1$, of $A_1 \AlphabetProduct A_2$ in order to get a coalgebra that exhibits properties such
as coassociativity and results in a pair of bialgebras in interaction~\cite{Foi20}.

\MakeReferences

\end{document}